\documentclass[a4paper]{article}

\usepackage[cp1251]{inputenc}
\usepackage[T2A]{fontenc}
\usepackage[russian]{babel}

\usepackage[a4paper,top=1.5cm,bottom=1.5cm,left=3cm,right=3cm,marginparwidth=1.75cm]{geometry}

\usepackage{amsmath}
\usepackage{graphicx}
\usepackage{comment}
\usepackage{amsfonts}

\title{Групповой анализ одномерного уравнения Больцмана. IV. Условие сохранения числа частиц}
\author{Боровских А.В., Платонова К.С.}
\date{}

\begin{document}
\maketitle
\setcounter{tocdepth}{2}
\tableofcontents
\section{Одномерное уравнение Больцмана и сопровождающие физические условия}

Работа является четвертой в серии, посвященной проблеме переноса группы симметрий с уравнения Больцмана на моментные величины, вычисления инвариантов этой группы в терминах моментных величин и решения в терминах инвариантов проблемы замыкания (с помощью так называемых уравнений состояния) моментной системы.

В ней, как и в предыдущих работах [1-3], рассматривается одномерный вариант уравнения, причем основное внимание сосредотачивается на том, чтобы явно в математической форме описать дополнительные ограничения физического характера на замены переменных, позволяющих получить содержательный и физически осмысленный результат.

Классическому уравнению Больцмана
$$f_t+cf_x+\mathcal{F}f_c=J,\eqno(1)$$
где $t$ -- время, $x$ -- координата, $c$ -- скорость, $f$ -- фазовая (по координатам и скоростям) плотность распределения количества частиц, $\mathcal{F}$ -- внешняя сила, $J$ -- интеграл столкновений, посвящено довольно много работ и целый ряд монографий (см., напр, [4-8]). Групповые свойства уравнения Больцмана исследовались ранее только в плане вычисления группы симметрий для уравнения с $\mathcal{F}=0$ (см. [9], [10]).

Имея в виду, в конечном счете, именно это уравнение, мы в настоящей работе рассматриваем упрощенный -- одномерный вариант уравнения с целью выяснения принципиальной возможности осуществить перенос группы симметрий с уравнения Больцмана на моментные величины.
Одномерный случай резко упрощает ситуацию, поскольку в этом случае интеграл столкновений оказывается нулевым (частицы при столкновении просто обмениваются импульсами, а это не меняет распределения $f$). Тем не менее сама задача группового анализа (см., напр. [11-13]) вносит целый ряд сложностей.

Во-первых, для эффективности группового анализа необходимо взять максимально широкий класс замен переменных (в нашем случае -- переменных $t, x, c, f$, а в случае группы эквивалентности -- и $\mathcal{F}$). В идеале все новые переменные могут быть функциями от всех старых, то есть речь идет о заменах общего вида.

Во-вторых, как следствие, для эффективности группового анализа необходимо все функциональные параметры, которые могут меняться, считаться произвольными функциями, и классификации производить по этому произвольному параметру. Поэтому функцию $\mathcal{F}$ следует считать произвольной функцией от $(t,x,c)$, а еще лучше -- и от $f$ (такого рода зависимости встречаются, например, в уравнении Власова [6]).

Эти два эвристических принципа группового анализа, будучи очень важными с математической точки зрения, тем не менее, с физической точки зрения могут привести к абсурду.
Так, совершенно понятно, что замена переменных, в которой система координат в $x$-ах (то есть в координатах) поворачивается, а в скоростях -- нет, физического смысла не имеет, и такие замены надо исключать. Поэтому, помимо самого уравнения, и даже, строго говоря, до рассмотрения самого уравнения, необходимо задать чисто математически ряд дополнительных условий, которые отражают сохранение того или иного физического смысла отношений между величинами, фигурирующими в уравнениях.

В предыдущих работах нами были сформулированы следующие ограничения, отражающие этот физический смысл. Это, во-первых, сохранение при заменах соотношений
$$dx=cdt,\qquad dc=\mathcal{F}dt,\eqno(2)$$ которые эквивалентны второму закону Ньютона и выражают собой тот факт, что $c$ -- это скорость в пространстве $x$-ов, а $\mathcal{F}$ -- это сила, приложенная к тому или иному объекту (в нашем случае -- к частице) в том же пространстве.

Во-вторых, нам понадобилось сохранение неизменными соотношений $$dt=dx=0,\eqno(3)$$ выражающих в аналитической форме тот факт, что при заменах переменных прямые $t={\rm const}$, $x={\rm const}$ перейдут в такие же прямые (а поскольку моментные величины вычисляются как интегралы по этим прямым, мы при замене переменных не потеряем смысла этих величин). К слову говоря, интегрирование по $c$ в трехмерном случае возникает уже в самом уравнении в интеграле столкновений, и поэтому там это условие потребуется как необходимое уже для инвариантности уравнения.

И, наконец, в-третьих, нам понадобилось условие сохранения при замене переменных числа частиц. Здесь имеется в виду, что при замене переменных (на физическом языке -- при переходе к другой системе отсчета) физический объект, а именно совокупность частиц в данном объеме, не меняется. Учитывая, что $f$ в уравнении Больцмана интерпретируется как фазовая плотность числа частиц, само количество (находящееся в бесконечно малом фазовом объеме $dxdc$) выражается величиной $f\,dxdc$, которая должна при заменах оставаться инвариантной.

Классификации одномерного уравнения Больцмана, учитывающие первое и второе условия, нами были выполнены в [1-3]. Что же касается третьего условия, то оказалось, что естественное, на первый взгляд, обобщение одномерного варианта уравнения (1) как уравнения
$$f_t+cf_x+\mathcal{F}(t,x,c)f_c=0$$
(и именно так записывают во всех справочниках вид левой части уравнения Больцмана) этому условию не соответствует: вычисление числа частиц, расположенного в некотором объеме фазового пространства $(x,c)$ и в образе этого объема по прошествии некоторого времени $t$ дает разные результаты. Поэтому учет третьего условия требует пересмотра самого уравнения (на что в свое время обратил наше внимание А.Н.Голубятников). Оказалось, что для выполнения требуемого условия уравнение должно иметь вид
$$f_t+cf_x+(\mathcal{F}(t,x,c)f)_c=0,\eqno(4)$$
то есть иметь дивергентную форму, которая в случае $\mathcal{F}=\mathcal{F}(t,x,c)$ раскрывается в виде
$$f_t+cf_x+\mathcal{F}(t,x,c)f_c+\mathcal{F}_c(t,x,c)f=0,\eqno(5)$$
а если $\mathcal{F}$ зависит не только от $(t,x,c)$, но и от $f$, то в виде
$$f_t+cf_x+\mathcal{F}(t,x,c,f)f_c+\mathcal{F}_c(t,x,c,f)f+\mathcal{F}_f(t,x,c,f)ff_c=
0,\eqno(5)$$
которые мы далее и будем рассматривать. Отметим, что для уравнений с $\mathcal{F}$, зависящей только от $t$ и $x$ (что обычно и подразумевается в литературе, посвященной уравнению Больцмана) все эти обобщения совпадают между собой, однако введение зависимости от $c$ и от $f$ резко меняет ситуацию.

Помимо того, что уравнение (4) является следствием закона сохранения величины $f$ (см., напр., [8]), оно обеспечивает этой величине еще одно важное свойство: "частицу"\, оказывается возможным отождествить с той или иной траекторией системы (2), так что величина $f$ уже считает плотность не неких "физических" \, частиц, а математических объектов -- траекторий, придавая каждой из них некий вес.

При этом удалось сделать еще один шаг -- модифицировать третье условие. Дело в том, что в третьем условии фигурирует фактически интегрирование по плоскости $t={\rm const}$, а при замене переменных такие плоскости переходят в поверхности произвольного вида. Требование же того, что плоскость переходит в плоскость представлялось нам излишним и обременительным. Теперь же, использовав соображения сохранения фазовых потоков, следующее из уравнений (4) или (5), удалось модифицировать третье условие, перейдя от поверхностей $t={\rm const}$ к поверхностям произвольным, задаваемым соотношением $t=\theta(x,c)$, и при этом инвариантной при замене переменных должна оставаться величина
$$(1-c\theta_x-\mathcal{F}\theta_c)f\,dx\,dc.\eqno(6)$$
Форма условия (6) обусловлена введением меры на множестве непрерывных траекторий, плотность которой вычисляется как плотность потока векторного поля $(1,c,\mathcal{F})\cdot f$ через поверхность $t=\theta(x,c)$.

В настоящей работе мы для уравнений (4) и (5) вычислим алгебру эквивалентности и алгебру симметрий как для случая $\mathcal{F}(t,x,c)$, так и для случая $\mathcal{F}(t,x,c,f)$, фактически закрыв вопрос о групповой классификации одномерных уравнений Больцмана.

Следует подчеркнуть специфику нашего подхода к групповому анализу, отличающуюся от довольно распространенного (и естественного, на первый взгляд), -- когда в центр внимания ставится уравнение или система уравнений, а группа симметрий рассматривается как подчиненный объект, оставляющий это уравнение или систему инвариантным. В нашем случае в центре внимания стоит именно группа (группа симметрий или группа эквивалентности), которая должна удовлетворять набору условий, каждое из которых может быть сформулировано в своей форме. Это может быть условие инвариантности какого-то уравнения, может быть условие инвариантности метрической или дифференциальной формы, может быть условие инвариантности тех или иных условий на поверхностях или многообразиях, и др.
Решение задачи групповой классификации с такой постановкой состоит в последовательном применении условий инвариантности того или иного объекта, получении из них условия на исходную алгебру Ли, и переходе к другому условию.

Обратим внимание, что такой "набор условий инвариантности различных соотношений" \ не равносилен "условию инвариантности объекта, удовлетворяющего всем этим соотношениям" \ даже в случае, когда такой объект существует. Желающие проверить это могут сравнить условие сохранения группой, действующей в плоскости переменных $(x,y)$, семейства прямых $x={\rm const}$ отдельно и семейства прямых $y={\rm const}$ отдельно (это группа замен $\bar x=\phi(x)$, $\bar y=\phi(y)$) с условием сохранения группой семейства объектов, для которых одновременно $x={\rm const}$ и $y={\rm const}$, то есть семейства точек (этому условию удовлетворяет группа всех диффеоморфизмов плоскости). В нашем же случае объекта, удовлетворяющего всем соотношениям (2)-(3), (6) и (4)/(5) просто не существует, поэтому попытка интерпретировать нашу постановку задачи в традиционном понимании бессмысленна.

\section{Формулировка основных результатов}

{\bf Теорема 1.} {\it
Алгебра Ли группы эквивалентности как уравнения (4) (с $\mathcal{F}(t,x,c)$), так и уравнения (5) (с $\mathcal{F}(t,x,c,f)$), сохраняющая дифференциальные соотношения (2), прямые (3) и форму (6), имеет вид:
$$\Xi=\tau\partial_t+\xi\partial_x+\alpha\partial_c+\eta\partial_f+\phi\partial_{\mathcal{F}},$$
где $\tau=\tau(t,x)$ и $\xi=\xi(t,x)$ -- произвольные функции,
$$\alpha = \xi_t + c\xi_x - c(\tau_t + c\tau_x), \quad \eta = f(3c\tau_x - 2\xi_x + \tau_t),\quad
\phi = \alpha_t + c\alpha_x + \mathcal{F}(\xi_x - 2\tau_t - 3c\tau_x).\eqno(7)$$

Эта алгебра порождает, вместе с отражением $\bar t=x$, $\bar x=t$, группу диффеоморфизмов пространства переменных $(t,x)$ и связанных с ними преобразований $c$, $f$ и $\mathcal{F}$:
$$\bar t=\varphi(t,x),\qquad \bar x=\psi(t,x),\qquad \bar c=\frac{\psi_t+c\psi_x}{\varphi_t+c\varphi_x},\qquad \bar f=f\frac{(\varphi_t+c\varphi_x)^3}{(\psi_x\varphi_t-\varphi_x\psi_t)^2},$$
$$\bar{\mathcal{F}}=\mathcal{F}\frac{\psi_x\varphi_t-\varphi_x\psi_t}{(\varphi_t+c\varphi_x)^3}+\frac1{\varphi_t+c\varphi_x}\left[ \left(\frac{\psi_t+c\psi_x}{\varphi_t+c\varphi_x}\right)_t+c\left(\frac{\psi_t+c\psi_x}{\varphi_t+c\varphi_x}\right)_x\right].\eqno(8)$$}

{\bf Следствие 1.} Классификацию уравнений Больцмана имеет смысл осуществлять с точностью до группы диффеоморфизмов пространства $(t,x)$.

{\bf Следствие 2.} В случае, когда $\mathcal{F}$ зависит от $f$, преобразования группы эквивалентности не меняют вид этой зависимости, поэтому вид зависимости от $f$ является классифицирующим признаком.

{\bf Следствие 3.} Поскольку скорость $c$ преобразуется дробно-рациональным образом, и в формулах (8) присутствуют только дробно-линейные функции от $c$, оказывается, что рациональная, иррациональная и различные типы трансцендентной зависимости $\mathcal{F}$ от $c$ несводимы друг к другу и поэтому тоже являются классифицирующим признаком.\\

{\bf Теорема 2.} {\it
1) Алгебра Ли группы симметрий уравнения (5) (с $\mathcal{F}(t,x,c,f)$), сохраняющей (2)-(3) и (6), имеет вид:
$\Xi=\tau\partial_t+\xi\partial_x+\alpha\partial_c+\eta\partial_f$,
где функции
$$\tau=\tau(t,x), \quad \xi=\xi(t,x),\quad \alpha = \xi_t + c\xi_x - c\tau_t - c^2\tau_x,\quad \eta = f(3c\tau_x - 2\xi_x + \tau_t)\eqno(9)$$
удовлетворяют уравнению
$$\alpha_t + c\alpha_x + \mathcal{F}\alpha_c - \mathcal{F}(\tau_t + c\tau_x) = \tau\mathcal{F}_t + \xi\mathcal{F}_x + \alpha\mathcal{F}_c + \eta\mathcal{F}_f.\eqno(10)$$

2) Алгебра Ли группы симметрий уравнения (4) (с $\mathcal{F}(t,x,c)$), сохраняющей (2)-(3) и (6), имеет вид
$\Xi=\tau\partial_t+\xi\partial_x+\alpha\partial_c+\eta\partial_f$,
где функции $\tau$, $\xi$, $\alpha$ и $\eta$ определяются теми же соотношениями (9) и удовлетворяют уравнению
$$\alpha_t + c\alpha_x + \mathcal{F}\alpha_c - \mathcal{F}(\tau_t + c\tau_x) = \tau\mathcal{F}_t + \xi\mathcal{F}_x + \alpha\mathcal{F}_c.\eqno(10')$$}

{\bf Замечание 1.} Как будет видно из доказательства, соотношения (9)-(10)/(10$'$) являются следствиями условий (2)-(3), (6), а уравнение (4)/(5) при этом сохраняется автоматически. Это не является неожиданным, поскольку уравнение (4)/(5) фактически представляет собой уравнение Лиувилля сохранения меры $fdxdc$ вдоль фазового потока, порожденного (2)-(3).

{\bf Замечание 2.} Это не означает "ненужности"\ условий (2)-(3): они, выражая сохранение физического смысла отношений между величинами, не зависят от того, какое именно уравнение мы рассматриваем. Могло рассматриваться и другое уравнение, и тогда условие инвариантности этого уравнения наложило бы дополнительное ограничение на группу симметрий.

{\bf Замечание 3.} Для трехмерного уравнения Больцмана аналогичное рассмотрение приводит не к тому, что уравнение Больцмана остается автоматически инвариантным, а к тому, что условие инвариантности уравнения дает распадающиеся классифицирующие уравнения -- одно для $\mathcal{F}$, а другое -- для интеграла столкновений, что оказывается важным для последующего анализа.\\

{\bf Теорема 3.} {\it
Алгебра $\Xi=\tau\partial_t+\xi\partial_x+\alpha\partial_c+\eta\partial_f$, определяемая соотношениями (9) и удовлетворяющая условию (10$\,'$), конечномерна для любой функции $\mathcal{F}(t,x,c)$.

Конечномерные нетривиальные алгебры (9), удовлетворяющие (10$\,'$) имеются для
функций $\mathcal{F}(t,x,c)$ классов (с точностью до преобразований из  группы, порожденной группой диффеоморфизмов (8)), представители которых приведены в таблице 1. В каждом случае I--IV предполагается, что функция $\mathcal{F}$ не лежит в предыдущем классе.

Для остальных $\mathcal{F}(t,x,c)$ уравнение (10$\,'$) не имеет других решений $\tau(t,x)$, $\xi(t,x)$, кроме тривиальных.}

\begin{figure}[h]
Таблица 1. Представители классов функций $\mathcal{F}$, для которых ($10'$) имеет решением нетривиальную алгебру ($A$ -- произвольная константа, $G(\cdot)$ и $\Phi(\cdot,\cdot)$ -- произвольные функции)
\begin{center}
\begin{tabular}{|c||c|l|}
\hline
& \text{          }$\mathcal{F}(t,x,c)$ & \text{          }Базис алгебры \\
\hline
\hline
I & $\mathcal{F} = 0$ & $\Xi_1= \partial_t$, $\Xi_2= \partial_x$, $\Xi_3 = t\partial_t$, $\Xi_4 = x\partial_x$, $\Xi_5 = x\partial_t$, $\Xi_6 = t\partial_x$, \\
& & $\Xi_7 = t^2\partial_t + tx\partial_x$, $\Xi_8 = tx\partial_t + x^2\partial_x$ \\
\hline
\hline
II.1 & $\mathcal{F}\!\!=\!\!\displaystyle A\exp\!\!{\int\!\!\frac{3c + a}{c^2 + bc + d}dc}\!\!$ & $\Xi_1 = \partial_t$, $\Xi_2 = \partial_x$, $\Xi_3 = (x+(a-b)t)\partial_t + ((a-2b)x-dt)\partial_x$\\
& $\mathcal{F} = Ac^a$ & $\Xi_1 = \partial_t$, $\Xi_2 = \partial_x$, $\Xi_3 = t\partial_t + \displaystyle\frac{a-2}{a-1}x\partial_x$\\
& $\mathcal{F} = A\exp(ac)$ & $\Xi_1 = \partial_t$, $\Xi_2 = \partial_x$, $\Xi_3 = t\partial_t + \displaystyle \left(x - \frac ta \right)\partial_x$\\
\hline
II.2 & $\mathcal{F} = \displaystyle A/x^{3}$ & $\Xi_1= \partial_t$, $\Xi_2 = 2t\partial_t + x\partial_x$, $\Xi_3=t^2\partial_t + tx\partial_x$\\
\hline
II.3 & $\mathcal{F} = \displaystyle A\left(1 - \frac{(t - ac)^2}{t^2 - 2ax}\right)^{\frac 32} $ & $\Xi_1 = a\partial_t + t\partial_x$, $\Xi_2 = t\partial_t + 2x\partial_x$,\\
&($a\ne 0$)& $\Xi_3 = (t^2-ax)\partial_t + tx\partial_x$\\
\hline
II.4 &$\!\!\mathcal{F}\!\! =\!\!\displaystyle A\left (\!\frac{(x-ct)^2+c^2+1}{x^2+t^2+1}\!\right)^{\frac 32}\!\!$& $\Xi_1 = (t^2+1)\partial_t+tx\partial_x$, $\Xi_2 = tx\partial_t+(x^2+1)\partial_x$,\\
& & $\Xi_3 = -x\partial_t+t\partial_x.$\\
\hline
\hline
\!\!III.1\!\!& $\mathcal{F} = G(c)$ & $\Xi_1=\partial_t$,\,\, $\Xi_2=\partial_x$\\
\hline
\!\!III.2\!\!& $\mathcal{F} = \displaystyle G(c)/t$ & $\Xi_1=\partial_x$, \,\,$\Xi_2=t\partial_t+x\partial_x$\\
\hline
\hline
IV& $\mathcal{F} = \Phi(x,c)$ & $\Xi_1=\partial_t$\\
 & $\mathcal{F} = \Phi(t,c)$ & $\Xi_1=\partial_x$\\
\hline
\end{tabular}
\end{center}
\end{figure}

{\bf Замечание 1.} Вторая и третья функции в классе II.1 являются частным случаем первой (и при замене переменных $\bar{t} = x$, $\bar{x} = t$ приводятся к первой при коэффициентах $b=d=0$ и $b = 1$, $d=0$ соответственно), но они являются наиболее простыми, после $\mathcal{F}=0$, поэтому мы их выделяем как частных, вспомогательных представителей этого семейства.

{\bf Замечание 2.} В теореме 3 фактически приводится уже окончательная классификация уравнений (4). Она является по существу перефразировкой результата из [3], где рассматривалось уравнение $f_t+cf_x+\mathcal{F}f_c=0$, отличающееся от (4) тем, что по $c$ дифференцируется только $f$. Групповая классификация этого уравнения с двумя лишь условиями (2)-(3) привела в точности к такому же уравнению (10$'$), и в [3] эти решения были перечислены. Поэтому мы доказательства теоремы 3 приводить не будем, сославшись на результат [3].\\

{\bf Теорема 4.} {\it
Алгебра $\Xi=\tau\partial_t+\xi\partial_x+\alpha\partial_c+\eta\partial_f$, определяемая соотношениями (9) и удовлетворяющая (10), конечномерна для любой функции $\mathcal{F}(t,x,c,f)$. 1
1
1

Конечномерные нетривиальные алгебры, удовлетворяющие (10), имеются для функций $\mathcal{F}(t,x,c,f)$ классов (с точностью до преобразований из  группы, порожденной группой диффеоморфизмов (8)), представители которых приведены в таблицах 2-5.

Каждая таблица представляет собой класс (классы занумерованы римскими цифрами) функций с одной и той же размерностью алгебры симметрий. Каждый класс делится на семейства (занумерованные уже арабскими цифрами), различающиеся типом зависимости от $f$, в каждом семействе имеются подсемейства, которые занумерованы тоже арабскими цифрами.В каждом подсемействе может быть указано несколько представителей -- канонический (для которого в третьей колонке указана алгебра) и вспомогательный (для которых указаны замены, приводящие их к каноническому виду). Вспомогательные представители удобны для уяснения связей с представителями следующих классов.

Для каждого класса I-V предполагается, что функция $\mathcal{F}$ не лежит в предыдущем классе (исключения и указания на соответствующие классы также приведены в таблицах) и что $\mathcal{F}_f\ne 0$.

Для остальных $\mathcal{F}(t,x,c,f)$ уравнение (10) не имеет других решений $\tau(t,x)$, $\xi(t,x)$ кроме тривиальных.}

\begin{figure}[h]
Таблица 2. Представители классов функций, для которых (10) имеет решением пятимерную алгебру ($P\ne 0$ -- произвольная константа)
\begin{center}
\begin{tabular}{|c||c|l|}
\hline
& \text{          }$\mathcal{F}(t,x,c,f)$ & \text{          }Базис алгебры \\
\hline
\hline
I & $\mathcal{F} = Pf^{-1}$ & $\Xi_1=\partial_t $, $\Xi_2 =\partial_x$, $\Xi_3=t\partial_t-x\partial_x$, $\Xi_4=t\partial_x $, $\Xi_5 =x\partial_t$\\
 & $\mathcal{F} = P(t^3f)^{-1}$ &$\bar t=-1/t$, $\bar x=x/t$\\
 & $\mathcal{F} = Pf^{-1}+Q$ &$\bar t=t$, $\bar x=x-Qt^2/2$\\
\hline
\end{tabular}
\end{center}
\end{figure}

\begin{figure}[h]
Таблица 3. Представители классов функций, для которых (10) имеет решением четырехмерную алгебру ($P\ne 0$, $Q$, $R$, $k$ -- произвольные константы)
\begin{center}
\begin{tabular}{|c||c|l|}
\hline
& \text{          }$\mathcal{F}(t,x,c,f)$ & \text{          }Базис алгебры / замена переменных\\
\hline
\hline
II.1 & $\mathcal{F} = Pf^j$, $j\ne 0$, $j\ne -1$& $\Xi_1= \partial_t$,  $\Xi_3 = (2j+1)t\partial_t+(j+2)x\partial_x$,\\
&& $\Xi_2= \partial_x$, $\Xi_4 =t\partial_x$\\
&$\mathcal{F}=Pt^{-3}f^j$&$\bar t=-1/t$, $\bar x=x/t$\\
& $\mathcal{F} = Pf+Qc+R$ ($Q\ne 0$)&  $\bar t=e^{-Qt}$, $\bar x=e^{-Qt}(x+Rt/Q)$ \\
& $\mathcal{F} = Pf+R$ &  $\bar t=t$, $\bar x=x-\frac 12 Rt^2$\\
&$\mathcal{F}=Pc^6f+Qc^3+Rc^2$&$\bar t=e^{Rx}$, $\bar x=e^{Rx}(t+Qx/R)$\\
& $\mathcal{F} = Pf^{-2}+Rc+S$ ($R\ne 0$)  & $\bar t=e^{Rt}$, $\bar x=x+St/R$   \\
&$\mathcal{F} = Pf^{-2}+S$&$\bar t=t$, $\bar x=x-St^2/2$ \\
& $\mathcal{F} = Pf^{-1/2}+Qc^2$ &$\bar t=t$, $\bar x=e^{-Qx}$\\
& $\mathcal{F} = Pf^{-1/2}+S$ & $\bar t=t$, $\bar x=x-St^2/2$\\
& $\mathcal{F} = Pt^{-3}f^{-1/2}$ & $\bar t=1/t$, $\bar x=x/t$\\
&$\mathcal{F}=\frac 1t[u_1c^{3(j+1)}(tf)^{j} + c\frac{1-j}{2j+1}]$&$\bar t=x$, $\bar x=t^{\frac{j+2}{2j+1}}$\\
\hline
II.2& $\mathcal{F} = P(c^2\pm k^2)^3f$, $k\ne 0$ & $\Xi_1= \partial_t$, $\Xi_2= \partial_x$, $\Xi_3 =t\partial_t+x\partial_x$, \\
&&$\Xi_4 = x\partial_t\mp k^2t\partial_x$\\
&$\mathcal{F} = P(c^2+k^2)^3f+Q(c^2+k^2)$&$\bar t=e^{Qx}\cos kQt$, $\bar x=e^{Qx}\sin kQt$  \\
&$\mathcal{F} = P(c^2-k^2)^3f+Q(c^2-k^2)$&$\bar t=e^{Qx}{\rm ch\,}kQt$, $\bar x=e^{Qx}{\rm sh\,}kQt$  \\
& $\mathcal{F} = Pc^3f$ & $\bar t=x+t$, $\bar x=x-t$\\
& $\mathcal{F} = Pc^3f+Qc^2+Rc\to \mathcal{F}=Pc^3f$,  & $\bar t=e^{Rt}$ (при $R\ne 0$),  $\bar x=e^{-Qx}$ (при $Q\ne 0$)\\
\hline
II.3&$\mathcal{F}=Pf^{-2}x^{-3}$&$\Xi_1= \partial_t$, $\Xi_2= x\partial_x$, $\Xi_3 = t\partial_t $,\\ &&$\Xi_4=t^2\partial_t+tx\partial_x $\\
& $\mathcal{F} = Pf^{-2}+Q(c^2+k^2)$, $Q\ne 0$ & $\bar t={\rm tg\,}kQt$, $\bar x=e^{-Qx}/\cos kQt$\\
& $\mathcal{F} = Pf^{-2}+Q(c^2-k^2)$, $Q\ne 0$ & $\bar t=e^{-2Qkt}$, $\bar x=e^{-Q(x+kt)}$\\
& $\mathcal{F} = Pf^{-2}+Qc^2$, $Q\ne 0$ & $\bar t=t$, $\bar x=e^{-Qx}$\\
& $\mathcal{F} = Pc^{-3}f^{-2}+Qc$, $Q\ne 0$ & $\bar t=x$, $\bar x=e^{Qt}$\\
& $\mathcal{F} = Pf^{-2}(x-ct)^{-3}$,  & $\bar t=x/t$, $\bar x=-1/t$,\\
\hline
II.4& $\mathcal{F} = P(t^2\pm k^2)^{-3/2}f^{-1/2}$ & $\Xi_1= \partial_x$,  $\Xi_2= t\partial_x$, $\Xi_3 =x\partial_x $, \\
&&$\Xi_4 = (t^2\pm k^2)\partial_t+tx\partial_x$\\
& $\mathcal{F} = Pf^{-1/2}+Q(c^2+k^2)$, $Q\ne 0$ & $\bar t={\rm tg\,} kQt$, $\bar x=e^{-Qx}/\cos kQt$\\
& $\mathcal{F} = Pf^{-1/2}+Q(c^2-k^2)$, $Q\ne 0$ & $\bar t={\rm cth\,} kQt$, $\bar x=e^{-Qx}/{\rm sh\,} kQt$\\
& $\mathcal{F} = Pf^{-1/2}+Rc+S$, $R\ne 0$ & $\bar t=\frac{1+e^{Rt}}{1-e^{Rt}}$, $\bar x=2\frac{x+St/R}{1+e^{Rt}}$\\
& $\mathcal{F} = Pt^{-3/2}f^{-1/2}$ &$\bar t=\frac{1+t}{1-t}$, $\bar x=\frac x{1-t}$ \\
\hline
II.5 & $\mathcal{F}=P\ln |f|+Q$ &  $\Xi_1= \partial_x$, $\Xi_2 =t\partial_x $,  $\Xi_3= \partial_t$, \\
&&$\Xi_4 = t\partial_t+(2x-\frac 32 Pt^2)\partial_x$\\
& $\mathcal{F}=t^{-3}(P\ln |f|+Q)$ &  $\bar t=-1/t$, $\bar x=x/t$ \\
\hline
\end{tabular}\\
\end{center}
\end{figure}

\begin{figure}
Таблица 4. Представители классов функций, для которых уравнение (10) имеет решением трехмерную алгебру ($G(\cdot)$ -- произвольная функция, для исключений указаны ссылки на таблицы 2 и 3)
\begin{center}
\begin{tabular}{|c||c|l|}
\hline
& &\\
& \text{          }$\mathcal{F}(t,x,c,f)$ & \text{          }Базис алгебры \\
\hline
\hline
& &$\Xi_1=-x\partial_t+t\partial_x$,\\
III.1 & $\mathcal{F}\!\!=\!\!\displaystyle\left(\frac {(x-ct)^2 + 1+c^2}{1+x^2+t^2}\right)^{\frac32} \!\!\!\!G\left(((x-ct)^2 + 1+c^2)^{\frac32}f\right)\!\!\!\!$ &  $\Xi_2 = tx\partial_t + (x^2+1)\partial_x$\\
&&$\Xi_3=(t^2+1)\partial_t + tx\partial_x$\\
\hline
III.2.1 & $\mathcal{F} = G(f)$ & $\Xi_1 = \partial_t$, $\Xi_2 = \partial_x$, $\Xi_3=t\partial_x$\\
& $G(f)\ne Pf^j+Q \to II.1, $ &\\
&$\mathcal{F}=c^3G(c^3f)$&$\bar t=x$, $\bar x=t$\\
\hline
III.2.2&$\mathcal{F}=\frac {1}{x^3}G(f)$&$\Xi_1=\partial_t$, $\Xi_2=2t\partial_t+x\partial_x$,\\ &&$\Xi_3=t^2\partial_t+tx\partial_x$\\
& $\mathcal{F} = \frac{c^3}{t}G(tc^3f)-\frac{c}{2t}$ &$\bar  t=x$, $\bar x=t^{1/2}$\\
\hline
III.2.3 & $\mathcal{F} = G(c)f $,  & $\Xi_1 = \partial_t$, $\Xi_2 = \partial_x$, \\
&$G(c)\ne (Pc^2+Qc+R)^3\to II.1, II.2$&$\Xi_3 =t\partial_t+x\partial_x$ \\
\hline
III.2.4 & $\mathcal{F} = G(t)f^{-1/2}$, & $\Xi_1 = \partial_x$, $\Xi_2 = t\partial_x$, $\Xi_3 =x\partial_x $\\
&$G(t)\ne (At^2+Bt+C)^{-3/2}\to II.1, II.4$&\\
\hline
III.3.1 & $\mathcal{F} = (c^2+n)^{3/2}G((c^2+n)^{3/2}f) $, & $\Xi_1 = \partial_t$, $\Xi_2 = \partial_x$, \\
&$G(z)\ne 1/z \to I$, $G(z)\ne Pz\to II.2$&$\Xi_3 =x\partial_t-nt\partial_x $\\
\hline
III.3.2 & $\mathcal{F} = c^mG(c^{3-m}f)$ & $\Xi_1 = \partial_t$, $\Xi_2 = \partial_x$, \\
&&$\Xi_3 =(m-1)t\partial_t +(m-2)x\partial_x$\\
&$\mathcal{F}\ne Pc^{-3}f^{-2}+Qc\to II.3$&\\
&$\mathcal{F}\ne Pf^{-1/2}+Qc\to II.4$&\\
&$\mathcal{F}\ne Pc^3f+Qc \to II.2$,&\\
&$\mathcal{F}\ne Pc^{3/2}f^{-1/2}+Qc^2 \to II.3$,&\\
&$\mathcal{F}\ne Pf^{-2}+Qc^2 \to II.3$&\\
& $\mathcal{F} = c^2G(cf)-\lambda c$, & $\bar t=e^{-\lambda t}$, $\bar x=x$,\\
&$\mathcal{F} = \frac {c^m}tG(tc^{3-m}f) + \frac c{t(m-2)}$&$\bar t=t^{\frac{m-1}{m-2}}$, $\bar x=x$\\
&$\mathcal{F} = \frac {c}tG(tc^2f)$&$\bar t=\ln|t|$, $\bar x=x$\\
\hline
III.3.3 & $\mathcal{F} = e^{mc}G(e^{-mc}f)$ & $\Xi_1 = \partial_t$, $\Xi_2 = \partial_x$, \\
&$G(z)\ne z \to II.1$ &$\Xi_3 =mt\partial_t+(mx-t)\partial_x$\\
&$\mathcal{F}=\frac 1t e^{mc}G(te^{-mc}f) +  \frac 1{mt}$&$\bar{t} = t$, $\bar{x} = x + (t-t\ln t)/m$\\
\hline
III.4.1 & $\mathcal{F}\!\!=\!\!(t^2\pm k^2)^{-3/2}e^{\lambda h(t)}G(fe^{2\lambda h(t)})$,  $h(t)\!\!=\!\!\int\!\!\frac{dt}{t^2\pm k^2}\!$, & $\Xi_1 = \partial_x$, $\Xi_2 = t\partial_x$, \\
& $G(z)\ne Pz^{-1/2}+Q\to II.3.1-II.3.3$ &$\Xi_3 =(t^2\pm k^2)\partial_t+(tx+\lambda x)\partial_x $\\
\hline
III.4.2 & $\mathcal{F} = t^{k-2}G(t^{2k-1}f)$, & $\Xi_1 = \partial_x$, $\Xi_2 = t\partial_x$, \\
&$G(z)\ne Pz^{-1/2}+Q\to II.4$, $\mathcal{F}\ne \frac Pf,\frac P{t^3f} \to I$,&$\Xi_3 =t\partial_t+kx\partial_x $\\
&$\mathcal{F}\ne P\ln |f|+Q,\frac {P\ln|f|+Q}{t^3}\to II.5$&\\
&$\mathcal{F}\ne Pf^j+Q, t^{-3}(Pf^j+Q) \to II.1$ &\\
& $\mathcal{F} = G(f)+Qc+R$, $Q\ne 0$, & $\bar t=e^{-Qt}$, $\bar x=x+\frac RQt$ \\
&$\mathcal{F}=\frac 1tG(tf)+\lambda\frac ct$, $\lambda\ne -1$&$\bar t=t^{\lambda+1}$, $\bar x=x$\\
\hline
III.4.3 & $\mathcal{F} = e^{kt}G(e^{2kt}f)$, $k\ne 0$ & $\Xi_1 = \partial_x$, $\Xi_2 = t\partial_x$,\\
&$G(z)\ne Pz^{-1/2}+Q\to II.1$& $\Xi_3 =\partial_t+kx\partial_x $\\
&$\mathcal{F}=\frac 1tG(tf)-\frac ct$&$\bar t=\ln|t|$, $\bar x=x$\\
\hline
III.5 &$\mathcal{F}=\left(1\pm \frac{(c\mp t)^2}{2x\mp t^2}\right)^{3/2}G\left((2x\mp t^2\pm(c\mp t)^2)^{3/2}f\right)$&$\Xi_1=\partial_t\pm t\partial_x$,$\Xi_2=t\partial_t+2x\partial_x$, \\
&&$\Xi_3=(t^2\mp x)\partial_t+tx\partial_x$\\
& $\mathcal{F} = \frac{(c^2\pm k^2)^{3/2}}{t}G\left(t(c^2\pm k^2)^{3/2}f\right)\pm\frac{c(c^2\pm k^2)}{k^2t}$, $k\ne 0$ & $\bar t=x$, $\bar x=\frac {k^2t^2\pm x^2}2$\\
\hline
\end{tabular} \\
\end{center}
\end{figure}
\begin{figure}
Таблица 5. Представители классов функций, для которых уравнение (10) имеет решением двух- или одномерную алгебру ($\Phi(\cdot,\cdot,\cdot)$ и $\Psi(\cdot,\cdot)$ -- произвольные функции)
\begin{center}
\begin{tabular}{|c||c|l|}
\hline
IV.1 & $\mathcal{F} = \Psi(c,f)$ & $\Xi_1=\partial_t$,\,\, $\Xi_2=\partial_x$\\
\hline
IV.2 & $\mathcal{F} = \displaystyle \frac {1}{t}\Psi(c,tf)$ & $\Xi_1=\partial_x$, \,\,$\Xi_2=t\partial_t+x\partial_x$\\
\hline
IV.3 & $\mathcal{F} = c\Psi(t,c^2f)$ & $\Xi_1=\partial_x$,\,\, $\Xi_2=x\partial_x$\\
\hline
IV.4 & $\mathcal{F} = \Psi(t,f)$  & $\Xi_1=\partial_x$,\,\, $\Xi_2=t\partial_x$\\
\hline
\hline
V & $\mathcal{F} = \Phi(t,c,f)$ & $\Xi_1=\partial_x$\\
& $\mathcal{F} = \Phi(x,c,f)$ & $\Xi_1=\partial_t$\\
\hline
\end{tabular}
\end{center}
\end{figure}

{\bf Замечание 1.} Отметим, что формат классификации в формулировке теоремы 4 определяется устройством группы эквивалентности (8): в таблицах 1-4 первой по иерархии классифицирующей характеристикой является размерность алгебры, второй -- тип зависимости от $f$, третьей -- тип зависимости от $c$.

{\bf Замечание 2.} {\em Каноническими} мы называем тех представителей, которым отвечают алгебры симметрий уравнения, являющиеся подалгебрами восьмимерной алгебры проективных преобразований в $\mathbb{R}^2$
введенной в работе С.Ли [14] (см. ткж. [15]) и состоящей из операторов
$$\partial_t,\quad \partial_x,\quad t\partial_t,\quad x\partial_t,\quad t\partial_x,\quad x\partial_x,\quad  t^2\partial_t+tx\partial_x,\quad tx\partial_t+x^2\partial_x.$$
С точки зрения анализа уравнения эта алгебра играет центральную роль: во-первых, для нее $$\xi_{tt}=2\xi_{tx}-\tau_{tt}=\xi_{xx}-2\tau_{tx}=\tau_{xx}=0,$$
и уравнение (10) оказывается однородным (а его решение -- соответствующий представитель класса -- не содержит, тем самым, лишних параметров). Во-вторых,
именно эта алгебра является алгеброй симметрий максимальной размерности для всех рассматриваемых уравнений, которая реализуется для уравнения с $\mathcal{F}=0$ (она и указана в формулировке теоремы 3). Сформулированный выше способ выбора представителей позволяет параллельно с иерархией семейств и подсемейств функций $\mathcal{F}$ выстроить и иерархию включений из соответствующих алгебр, что мы рассматриваем как идеальную форму такого рода результатов.

\section{Доказательство теорем 1 и 2.}
Доказательства обоих теорем мы будем проводить параллельно, так как оба результата устанавливаются по единой схеме: сначала выводятся условия на алгебры эквивалентности и симметрий, сохраняющие соотношения (2), затем - на алгебры, сохраняющие прямые (3), потом -- сохраняющие форму (6) и, наконец, оставляющие инвариантными уравнения (4) или (5).

\subsection{Условие сохранения соотношений (2)}
Как было сказано выше, физический смысл переменных $t,x,c$ не позволяет нам произвольно менять переменные, не учитывая связи между ними. Эти связи математически выражаются соотношениями $dx = cdt$ и $dc=\mathcal{F}dt$, которые при заменах переменных должны сохраняться. Найдем условия сохранения алгеброй эквивалентности и алгеброй симметрий данных соотношений для $\mathcal{F}(t,x,c,f)$ и $\mathcal{F}(t,x,c)$.

\subsubsection{Алгебра эквивалентности для $\mathcal{F}(t,x,c,f)$}

Мы будем рассматривать произвольные замены в пространстве переменных $(t,x,c,f,\mathcal{F})$, при которых одновременно сохраняются тождества $dx=cdt$ и $dc=\mathcal{F}dt$.

Поскольку в новых координатах (обозначим их теми же буквами с чертой) должно быть выполнено $d\bar{x}-\bar{c}d\bar{t}=0$, мы, подставляя в это равенство асимптотические (вблизи $a=0$) представления замены $\bar t\sim t+a\tau$, $\bar x\sim x+a\xi$, $\bar c\sim c+a\alpha$, $\bar f\sim f+a\eta$, $\bar{\mathcal{F}}\sim\mathcal{F}+a\phi$, получаем
$$dx +a\xi_tdt + a\xi_xdx + a\xi_cdc + a\xi_fdf + a\xi_{\mathcal{F}}d\mathcal{F} - (c + a\alpha)(dt + a\tau_tdt + a\tau_xdx + a\tau_cdc + a\tau_fdf+ a\tau_{\mathcal{F}}d\mathcal{F}) \sim 0$$
с точностью до бесконечно малых первого порядка по $a$. Представляя полученное выражение в виде многочлена по $a$, выпишем слагаемые первого порядка по $a$, что дает соотношение
$$(c\tau_t + \alpha - \xi_t)dt + (c\tau_x - \xi_x)dx + (c\tau_c - \xi_c)dc + (c\tau_f - \xi_f)df + (c\tau_{\mathcal{F}} - \xi_{\mathcal{F}})d\mathcal{F} = 0.$$

Подставим $dx = cdt$ и $dc=\mathcal{F}dt$ и, воспользовавшись тем, что приращения $dt$, $df$ и $d\mathcal{F}$ независимы, получаем:
$$\xi_t + c\xi_x + \mathcal{F}\xi_c - c(\tau_t + c\tau_x + \mathcal{F}\tau_c) = \alpha, \quad c\tau_f - \xi_f = 0 , \quad c\tau_{\mathcal{F}} - \xi_{\mathcal{F}} = 0.\eqno(11)$$
Рассуждая аналогично, найдем условия сохранения соотношения $d{c}-\mathcal{F}d{t}=0$:
$$\alpha_t + c\alpha_x + \mathcal{F}\alpha_c - \mathcal{F}(\tau_t + c\tau_x + \mathcal{F}\tau_c) = \phi, \,\, \alpha_f - \mathcal{F}\tau_f = 0,\,\,  \alpha_{\mathcal{F}} - \mathcal{F}\tau_{\mathcal{F}} = 0.\eqno(12)$$
Из совместности условий $\alpha_f = \mathcal{F}\tau_f$ и $\alpha_{\mathcal{F}} = \mathcal{F}\tau_{\mathcal{F}}$ следует, что $\tau_f =0$, а значит, и $\xi_f = \alpha_f = \phi_f = 0$. Далее из $\alpha_{\mathcal{F}} = \mathcal{F}\tau_{\mathcal{F}}$ следует, что $\xi_c = c\tau_c$. Условие совместности $\xi_c = c\tau_c$ и
$\xi_{\mathcal{F}} = c\tau_{\mathcal{F}}$ дает $\tau_{\mathcal F}=0$, откуда $\xi_{\mathcal{F}}=\alpha_{\mathcal{F}}=0$.

Таким образом, в случае $\mathcal{F}=\mathcal{F}(t,x,c,f)$ алгебра эквивалентности имеет вид
$$\Xi=\tau(t,x,c)\partial_t+\xi(t,x,c)\partial_x+\alpha(t,x,c)\partial_c+\eta(t,x,c,f,\mathcal{F})\partial_f+\phi(t,x,c,\mathcal{F})\partial_{\mathcal{F}},$$
где $\tau(t,x,c)$ и $\xi(t,x,c)$ -- произвольные функции, удовлетворяющие равенству $\xi_c-c\tau_c=0$,
$$\alpha = \xi_t + c\xi_x - c(\tau_t + c\tau_x),\quad
\phi = \alpha_t + c\alpha_x + \mathcal{F}\alpha_c - \mathcal{F}(\tau_t + c\tau_x + \mathcal{F}\tau_c),\eqno(13)$$
а $\eta(t,x,c,f,\mathcal{F})$ -- пока что произвольная функция.

Остается заметить, что от условия $\xi_c-c\tau_c=0$ можно избавиться, если ввести функцию $\beta(t,x,c)=\xi-c\tau$. В этом случае мы получим, что $\tau=-\beta_c$, $\xi=\beta+c\tau=\beta-c\beta_c$, и никаких дополнительных условий уже не требуется. Таким образом, алгебра определяется одной произвольной функцией $\beta(t,x,c)$ и одной произвольной функцией $\eta(t,x,c,f,\mathcal{F})$.

\subsubsection{Алгебра эквивалентности для $\mathcal{F}(t,x,c)$}

В этом случае нам понадобится дополнительно потребовать, чтобы при заменах переменных сохранялось условие $\mathcal{F}_f=0$, что дает условие $\tau_f=\xi_f=\alpha_f=\phi_f=0$.

Остальные рассуждения аналогичны предыдущему случаю: условие сохранения $dx=cdt$ дает равенство
$$a\xi_tdt + dx + a\xi_xdx + a\xi_cdc + a\xi_fdf + a\xi_{\mathcal{F}}d\mathcal{F} -$$
$$-(c + a\alpha)(dt + a\tau_tdt + a\tau_xdx + a\tau_cdc + a\tau_fdf+ a\tau_{\mathcal{F}}d\mathcal{F}) = 0,$$
с точностью до бесконечно малых первого порядка по $a$, расщепление которого дает систему уравнений (11), а условия сохранения соотношения $d{c}-\mathcal{F}(t,x,c)d{t}=0$ -- систему (12).

Единственное отличие от предыдущего случая -- что равенства $\tau_f = \xi_f = \alpha_f = \phi_f = 0$ мы получаем не из условий совместности, а они у нас изначально равны нулю в силу сохранения равенства $\mathcal{F}_f=0$.

Таким образом, алгебра эквивалентности для этого случая оказалась той же самой, что и в предыдущем пункте, что означает не только то, что алгебра эквивалентности предыдущего пункта переводит в себя подсемейство функций $\mathcal{F}(t,x,c)$, не зависящих от $f$, но и то, что в этом подсемействе нет дополнительных преобразований, переводящих его в себя.

\subsubsection{Алгебра симметрий для $\mathcal{F}(t,x,c,f)$}

Рассуждения здесь аналогичны рассуждениям в предыдущих пунктах, с той только разницей, что $\mathcal{F} = \mathcal{F}(t,x,c,f)$ предполагается теперь фиксированной функцией, а не самостоятельной переменной, и поэтому, с одной стороны, все компоненты алгебры уже зависят только от четырех переменных $(t,x,c,f)$, а с другой -- $\Xi\mathcal{F}$ дает не компоненту $\phi$, а обычное дифференцирование функции:
$$\Xi\mathcal{F}=\mathcal{F}_t\tau+\mathcal{F}_x\xi+\mathcal{F}_c\alpha+\mathcal{F}_f\eta.$$

Условие сохранения $dx=cdt$ дает нам тогда равенство
$$(c\tau_t-\xi_t +\alpha)dt + (c\tau_x - \xi_x)dx +(c\tau_c- \xi_c)dc + (c\tau_f- \xi_f)df=0,$$
а условие сохранения соотношения $dc=\mathcal{F}dt$ -- равенство
$$(\mathcal{F}\tau_x - \alpha_x)dx + (\mathcal{F}\tau_c - \alpha_c)dc + (\mathcal{F}\tau_f - \alpha_f)df + (\mathcal{F}\tau_t - \alpha_t +\tau\mathcal{F}_t+\xi\mathcal{F}_x+ \alpha\mathcal{F}_c+\eta\mathcal{F}_f)dt = 0.$$

Соответственно место соотношений (11) и (12) занимают
$$\xi_f - c\tau_f = \alpha_f - \mathcal{F}\tau_f = 0,\eqno(14)$$
$$\alpha = \xi_t + c\xi_x + \mathcal{F}\xi_c - c\tau_t - c^2\tau_x - c\mathcal{F}\tau_c,\eqno(15)$$
$$\alpha_t + c\alpha_x + \mathcal{F}\alpha_c - \mathcal{F}\tau_t - c\mathcal{F}\tau_x - \mathcal{F}^2\tau_c = \tau\mathcal{F}_t + \xi\mathcal{F}_x + \alpha\mathcal{F}_c + \eta\mathcal{F}_f.\eqno(16)$$
Подставим $\alpha$ из (15) во второе соотношение $(14)$ и получим $\mathcal{F}_f(\xi_c-c\tau_c)=0$,
так что возможны два варианта:\\
1) $\mathcal{F} = \mathcal{F}(t,x,c)$, что будет рассмотрено в следующем пункте;\\
2) $\xi_c - c\tau_c = 0$, тогда из совместности $\xi_c = c\tau_c$ и $\xi_f = c\tau_f$ получаем $\tau_f = \xi_f = \alpha_f = 0$.

Таким образом, алгебра симметрий имеет вид
$$\Xi=\tau(t,x,c)\partial_t+\xi(t,x,c)\partial_x+\alpha(t,x,c)\partial_c+\eta(t,x,c,f)\partial_f,$$
где $\tau(t,x,c)$ и $\xi(t,x,c)$ -- произвольные функции, удовлетворяющие равенству $\xi_c-c\tau_c=0$ (как мы уже выше отмечали, это означает, что они выражаются через одну произвольную функцию $\beta(t,x,c)$: $\tau=-\beta_c$, $\xi=\beta-c\beta_c$),
$$\alpha=\xi_t - c\tau_t + c\xi_x - c^2\tau_x,\eqno(17)$$ а $\eta$ -- пока что произвольная функция, и эти четыре компоненты удовлетворяют уравнению
$$\alpha_t + c\alpha_x + \mathcal{F}\alpha_c - \mathcal{F}\tau_t - c\mathcal{F}\tau_x - \mathcal{F}^2\tau_c = \tau\mathcal{F}_t + \xi\mathcal{F}_x + \alpha\mathcal{F}_c + \eta\mathcal{F}_f.\eqno(18)$$

\subsubsection{Алгебра симметрий для $\mathcal{F}(t,x,c)$}
Она вычисляется аналогично предыдущему случаю, с тем только отличием, что всюду $\mathcal{F}_f=0$. Это дает те же самые соотношения (14)-(16), из которых, впрочем, не следует никаких дополнительных условий, поскольку второе соотношение в (14) оказывается тождеством. В итоге алгебра симметрий имеет вид
$$\Xi=\tau(t,x,c,f)\partial_t+\xi(t,x,c,f)\partial_x+\alpha(t,x,c,f)\partial_c+\eta(t,x,c,f)\partial_f,$$
где $\tau(t,x,c,f)$ и $\xi(t,x,c,f)$ -- произвольные функции, удовлетворяющие условию $\xi_f - c\tau_f =0$, $\alpha = \xi_t + c\xi_x + \mathcal{F}\xi_c - c\tau_t - c^2\tau_x - c\mathcal{F}\tau_c$, а $\eta$ --  произвольная функция, и эти четыре компоненты удовлетворяют уравнению
$$\alpha_t + c\alpha_x + \mathcal{F}\alpha_c - \mathcal{F}\tau_t - c\mathcal{F}\tau_x - \mathcal{F}^2\tau_c = \tau\mathcal{F}_t + \xi\mathcal{F}_x + \alpha\mathcal{F}_c.\eqno(19)$$

\subsection{Условие сохранения прямых (3)}
Для сохранения интегрирования по прямым при вычислении моментных функций при заменах переменных необходимо, чтобы прямые $t = const$, $x = const$ перешли в прямые того же вида. Для алгебры эквивалентности и алгебры симметрий это означает сохранение соотношений $dt = 0$ и $dx = 0$, то есть выполнение следующих условий: $\tau_f = \xi_f = \tau_c = \xi_c = 0$. Это избавляет нас от квадратичных по $\mathcal{F}$ членов в (13), (18) и (19), а также от условий (14), которые теперь выполнены автоматически.

\subsection{Условие сохранения числа частиц (6)}
Последнее условие, которое необходимо наложить на алгебры симметрий и эквивалентности, -- сохранение числа частиц (6).

\subsubsection{Алгебра эквивалентности для $\mathcal{F}(t,x,c,f)$}
Напомним, что по итогам предыдущих рассуждений алгебра эквивалентности уравнения с $\mathcal{F}(t,x,c,f)$ имеет вид:
$$\tau = \tau(t,x), \quad \xi = \xi(t,x), \quad \alpha = \xi_t + c\xi_x - c(\tau_t + c\tau_x), \quad \phi = \alpha_t + c\alpha_x + \mathcal{F}(\alpha_c - \tau_t - c\tau_x),\eqno(20)$$
$$\eta = \eta(t,x,c,f,\mathcal{F}).$$

Выразим условие (6) в новых переменных и совершим обратный переход к $t,x,c,f$  и $\mathcal{F}$, учитывая, что $dt = \theta_xdx + \theta_cdc$, так как интегрирование происходит по поверхности $t = \theta(x,c)$:
$$(f+ a\eta)(dx + a(\xi_t\theta_x + \xi_x)dx + a\xi_t\theta_cdc)(dc + a(\alpha_t\theta_x + \alpha_x)dx + a(\alpha_t\theta_c +  \alpha_c)dc)\cdot$$
$$\cdot (1 - c\theta_x - \mathcal{F}\theta_c - ac\Theta^x - a\alpha\theta_x - a\mathcal{F}\Theta^c - a\phi\theta_c) \sim f(1-c\theta_x - \mathcal{F}\theta_c)dxdc,\eqno(21)$$
где $\Theta^x$ и $\Theta^c$ -- компоненты векторного поля, отвечающие за преобразование $\theta_x$ и $\theta_c$ соответственно. Слагаемые нулевого порядка по $a$ в левой части (21) совпадают с $f(1-c\theta_x - \mathcal{F}\theta_c)dxdc$. Слагаемые первого порядка по $a$ в таком случае должны равняться нулю:
$$f[(\xi_t\theta_x + \xi_x)(1 - c\theta_x - \mathcal{F}\theta_c) + (\alpha_t\theta_c + \alpha_c)(1 - c\theta_x - \mathcal{F}\theta_c) - (c\Theta^x + \alpha\theta_x + \mathcal{F}\Theta^c + \phi\theta_c)] + \eta(1 - c\theta_x - \mathcal{F}\theta_c) = 0.$$

Воспользуемся формулами для $\Theta^x$ и $\Theta^c$:
$$\Theta^x = \tau_x + \tau_t\theta_x - \xi_x\theta_x - \xi_t\theta_x^2 - \alpha_x\theta_c - \alpha_t\theta_c\theta_x, \quad \Theta^c = \tau_t\theta_c - \xi_t\theta_x\theta_c - \alpha_c\theta_c - \alpha_t\theta_c^2,$$
подстановка их в полученное равенство дает
$$f(\xi_t\theta_x + \alpha_t\theta_c + \xi_x + \alpha_c)(1 - c\theta_x - \mathcal{F}\theta_c) + \eta(1 - c\theta_x - \mathcal{F}\theta_c) -$$
$$- cf(\tau_x + \tau_t\theta_x - \xi_x\theta_x - \xi_t\theta_x^2 - \alpha_x\theta_c - \alpha_t\theta_c\theta_x) - f\alpha\theta_x - f\mathcal{F}(\tau_t\theta_c - \xi_t\theta_x\theta_c - \alpha_c\theta_c - \alpha_t\theta_c^2) - f\phi\theta_c = 0.$$

Так как условие (6) должно выполняться для любой функции $\theta(x,c)$, величины $\theta_x$ и $\theta_c$ независимы, а значит коэффициенты при разных мономах от этих величин должны быть нулевыми:
$$f(\xi_t - c\tau_t - c\alpha_c - \alpha) = c\eta, \quad f(\alpha_t + c\alpha_x - \mathcal{F}(\xi_x + \tau_t) - \phi) = \mathcal{F}\eta, \quad f(c\tau_x - \xi_x - \alpha_c) = \eta.$$
Первое и второе соотношения является следствием третьего в силу определения $\alpha$ и $\phi$ из (20). Остается выражение
$$\eta = f(c\tau_x - \xi_x - \alpha_c),\eqno(22)$$
которое дополняет (20).

\subsubsection{Алгебра эквивалентности для $\mathcal{F}(t,x,c)$}
Вычисление условий сохранения соотношения (6) не использует зависимость/независимость функции $\mathcal{F}$ от $f$. Поэтому результаты в этом случае аналогичны результатам п. 3.3.1.

\subsubsection{Алгебра симметрий для $\mathcal{F}(t,x,c,f)$}
Для сохранения соотношения (6) должно быть выполнено:
$$f(1 + c\theta_x + \mathcal{F}\theta_c)dxdc  = (f+ a\eta)(a\xi_tdt + a\xi_xdx + dx)(a\alpha_tdt + a\alpha_xdx + a\alpha_cdc + dc)\cdot$$
$$\cdot (1 - (c+a\alpha)(\theta_x + a\Theta^x) - (\mathcal{F}+ a\mathcal{F}_t\tau + a\mathcal{F}_x\xi + a\mathcal{F}_c\alpha + a\mathcal{F}_f\eta)(\theta_c + a\Theta^c)).$$
Повторяя рассуждения п. 3.3.1, получим соотношения
$$c\eta = f(\xi_t - c\alpha_c - c\tau_t - \alpha), \quad \mathcal{F}\eta = f(\alpha_t + c\alpha_x - \mathcal{F}\tau_t - \mathcal{F}\xi_x - \mathcal{F}_t\tau - \mathcal{F}_x\xi - \mathcal{F}_c\alpha - \mathcal{F}_f\eta),$$
$$\eta = f(c\tau_x-\xi_x - \alpha_c ). $$
первое из которых является следствием третьего в силу (17), а второе -- в силу уравнения (18). То есть в итоге имеем:
$$\tau=\tau(t,x), \quad  \xi = \xi(t,x), \quad \alpha = \xi_t + c\xi_x - c\tau_t - c^2\tau_x, \quad \eta = f(3c\tau_x - 2\xi_x + \tau_t)\eqno(23)$$
и уравнение (18).

\subsubsection{Алгебра симметрий для $\mathcal{F}(t,x,c)$}
Так как вычисление условий сохранения числа частиц аналогично для $\mathcal{F}(t,x,c)$ и $\mathcal{F}(t,x,c,f)$ и, в силу сохранения физического смысла моментов, $\tau = \tau(t,x)$, $\xi = \xi(t,x)$, то мы получим те же соотношения (23) и уравнение (19).


\subsection{Условие сохранения уравнения}
Теперь рассмотрим главное -- условия инвариантности уравнения (5) опять же в различных версиях: с функцией $\mathcal{F}$, зависящей от $(t,x,c)$ и зависящей от $(t,x,c,f)$.

\subsubsection{Алгебра эквивалентности для $\mathcal{F}(t,x,c,f)$}
Условие инвариантности уравнения (5) имеет вид:
$$\eta^t + c\eta^x + \mathcal{F} \eta^c + \alpha f_x + \phi f_c + \mathcal{F}_c\eta + \phi^c f + \mathcal{F}_ff_c\eta + \mathcal{F}_ff\eta^c + \phi^f f_cf= 0\eqno(24)$$
на поверхности, заданной уравнением (5). Здесь подразумевается, что, с точностью до второго порядка по групповому параметру, $\bar f_{\bar t} \sim f_t + a\eta^t$, $\bar f_{\bar x} \sim f_x + a\eta^x$, $\bar f_{\bar c} \sim f_c + a\eta^c$ и $\bar{\mathcal{F}}_{\bar{c}} \sim \mathcal{F}_c + a\phi^c$, $\bar{\mathcal{F}}_{\bar f} \sim \mathcal{F}_f + a\phi^f$.

Специфика ситуации состоит в том, что функция $\mathcal{F}$ сама зависит от $t,x,c,f$, то есть мы имеем дело со сложной функцией. И если продолжение на производные $\mathcal{F}$ (с учетом соотношений, полученных в п. 3.3.1) осуществляется по тем же формулам, что и обычно
$$\phi^c = \phi_c + \phi_{\mathcal{F}} \mathcal{F}_c - \mathcal{F}_t\tau_c - \mathcal{F}_x\xi_c - \mathcal{F}_c\alpha_c -\mathcal{F}_f(\eta_c + \eta_{\mathcal{F}}\mathcal{F}_c), \quad \phi^f = \phi_{\mathcal{F}} \mathcal{F}_f -\mathcal{F}_f(\eta_f + \eta_{\mathcal{F}}\mathcal{F}_f),$$
поскольку это продолжение на производные $\mathcal{F}$ по своим (формальным) аргументам, так что здесь не важно, что подставлено вместо этих аргументов, то для продолжения на производные $f$ приходится пользоваться несколько другими формулами
$$\eta^t = \eta_t + \eta_ff_t + \eta_{\mathcal{F}}\mathcal{F}_t + \eta_{\mathcal{F}}\mathcal{F}_ff_t - f_t\tau_t - f_x\xi_t - f_c\alpha_t,$$
$$\eta^x = \eta_x + \eta_ff_x + \eta_{\mathcal{F}}\mathcal{F}_c + \eta_{\mathcal{F}}\mathcal{F}_ff_x - f_t\tau_x - f_x\xi_x - f_c\alpha_x,$$
$$\eta^c = \eta_c + \eta_ff_c + \eta_{\mathcal{F}}\mathcal{F}_c + \eta_{\mathcal{F}}\mathcal{F}_ff_c - f_t\tau_c - f_x\xi_c - f_c\alpha_c,$$
отличающимися тем, что полные производные от компонент векторного поля вычисляются с учетом того, что их зависимость от функции $f(t,x,c)$ является не только прямой, но и опосредованной функцией $\mathcal{F}$, которая также является одним из аргументов этих компонент. Однако, ввиду соотношения (22), в котором зависимость $\eta$ от $\mathcal{F}$ отсутствует, мы можем, удалив слагаемые с $\eta_{\mathcal{F}}$, пользоваться обычными формулами
$$\eta^t = \eta_t + \eta_ff_t - f_t\tau_t - f_x\xi_t - f_c\alpha_t,\quad
\eta^x = \eta_x + \eta_ff_x - f_t\tau_x - f_x\xi_x - f_c\alpha_x,$$
$$\eta^c = \eta_c + \eta_ff_c - f_t\tau_c - f_x\xi_c - f_c\alpha_c.\eqno(25)$$

Подставим формулы для преобразований производных $f$ и $\mathcal{F}$ в условие (24), заменим $f_t$ на $-cf_x - \mathcal{F}f_c - \mathcal{F}_cf - \mathcal{F}_ff_cf$ и преобразуем полученное выражение с учетом того, что $\tau$ и $\xi$ зависят только от $(t,x)$:
$$D_F\eta - \eta_f(\mathcal{F}_cf + \mathcal{F}_ff_cf) + \mathcal{F}_c\eta + \mathcal{F}_ff_c\eta + \mathcal{F}_cf D_F\tau  + \mathcal{F}_ff_cfD_F\tau - \mathcal{F}_fff_c\alpha_c +$$
$$+ f(\phi_c + \phi_{\mathcal{F}} \mathcal{F}_c - \mathcal{F}_c\alpha_c) + \mathcal{F}_f ff_c \eta_f + (\phi_{\mathcal{F}} \mathcal{F}_f -\mathcal{F}_f\eta_f)ff_c = 0,$$
для удобства мы обозначили $D_F = \partial_t + c\partial_x + \mathcal{F}\partial_c$.

Так как $\mathcal{F}_t$, $\mathcal{F}_x$, $\mathcal{F}_c$, $\mathcal{F}_f$ -- произвольные величины, то коэффициенты при них должны равняться нулю:
$$D_F\eta + f\phi_c =0, \,\,\, \eta-\eta_ff + f(\phi_{\mathcal{F}} + D_F\tau-\alpha_c)=0, \,\,\, f_c(\eta -\eta_f f) + f_cf (\phi_{\mathcal{F}}+D_F\tau - \alpha_c)= 0.$$
Подстановка в эти соотношения выражений
$$\phi = D_F\alpha-\mathcal{F}D_F\tau, \quad
\eta = f(c\tau_x - \xi_x - \alpha_c),\quad
\alpha = D_F\xi- cD_F\tau$$
из (20) и (22) обращает все эти соотношения в тождества. Таким образом, инвариантность семейства уравнений (5) обеспечивается условиями инвариантности (2), (3) и (6) автоматически и алгебра эквивалентности для семейства (5) имеет вид (20), (22), что совпадает с частью утверждения теоремы 1, касающейся уравнения (5).

\subsubsection{Алгебра эквивалентности для $\mathcal{F}(t,x,c)$}

Условие инвариантности уравнения (4) совпадает с условием инвариантности из предыдущего пункта, за исключением того, что здесь мы считаем $\mathcal{F}_f = 0$. Повторяя для этой версии рассуждения предыдущего пункта, мы получаем тот же результат: соотношения (20), (22) автоматически обеспечивают инвариантность семейства уравнений (4), что совпадает с частью утверждения теоремы 1, касающейся уравнения (4).

Осталось отметить, что для любого диффеоморфизма $\bar t=\varphi(t,x)$, $\bar x=\psi(t,x)$ соответствующие преобразования $c$ и $\mathcal{F}$ следуют из (2), а преобразования $f$ -- из (6), что и дает нам формулы (8). Доказательство теоремы 1 на этом полностью завершено.

Следует обратить внимание на существенный для будущего рассмотрения уже трехмерного уравнения феномен автоматического сохранения инвариантным дифференциальной части уравнения Больцмана при выполнении условий инвариантности (2)--(3) и (6).

\subsubsection{Алгебра симметрий для $\mathcal{F}(t,x,c,f)$}
Компоненты алгебры симметрий уравнения Больцмана должны на поверхности, заданной уравнением (5), удовлетворять уравнению Ли:
$$\eta^t + c\eta^x + \mathcal{F}\eta^c + \mathcal{F}_c\eta + \alpha f_x + f_c(\tau\mathcal{F}_t + \xi\mathcal{F}_x + \alpha\mathcal{F}_c + \eta\mathcal{F}_{f}) + f(\tau \mathcal{F}_{tc} + \xi \mathcal{F}_{xc} + \alpha \mathcal{F}_{cc} + \eta\mathcal{F}_{cf}) +$$
$$+ \mathcal{F}_f f_c\eta + \mathcal{F}_f f\eta^c + ff_c(\tau \mathcal{F}_{tf} + \xi\mathcal{F}_{xf} + \alpha\mathcal{F}_{cf} + \eta\mathcal{F}_{ff}) = 0. \eqno(26)$$
Формулы продолжения на производные $f$ вычисляются в этом случае стандартно -- по формулам (25), которые ввиду независимости $\tau$ и $\xi$ от $c$ приобретают вид
$$\eta^t = \eta_t + \eta_ff_t - f_t\tau_t - f_x\xi_t - f_c\alpha_t,\
\eta^x = \eta_x + \eta_ff_x - f_t\tau_x - f_x\xi_x - f_c\alpha_x,\
\eta^c = \eta_c + \eta_ff_c - f_c\alpha_c.$$
Подставим их в (26) и преобразуем полученное выражение, имея в виду, что $f_{t} = - cf_{x} - \mathcal{F}f_{c} - \mathcal{F}_cf - \mathcal{F}_ff_{c}f$. Получаем уравнение
$$\eta_t+c\eta_x+F\eta_c-f\eta_f(\mathcal{F}_c + \mathcal{F}_ff_{c})+(cf_{x} + \mathcal{F}f_{c} + \mathcal{F}_cf + \mathcal{F}_ff_{c}f)(\tau_t+c\tau_x) -f_x(\xi_t+c\xi_x)-f_c(\alpha_t+c\alpha_x+\mathcal{F}\alpha_c)+$$
$$+ \mathcal{F}_c\eta + \alpha f_x + f_c(\tau\mathcal{F}_t + \xi\mathcal{F}_x + \alpha\mathcal{F}_c + \eta\mathcal{F}_{f}) + f(\tau \mathcal{F}_{tc} + \xi \mathcal{F}_{xc} + \alpha \mathcal{F}_{cc} + \eta\mathcal{F}_{cf}) +$$
$$+ \mathcal{F}_f f_c\eta + \mathcal{F}_f f(\eta_c + \eta_ff_c - f_c\alpha_c) + ff_c(\tau \mathcal{F}_{tf} + \xi\mathcal{F}_{xf} + \alpha\mathcal{F}_{cf} + \eta\mathcal{F}_{ff}) = 0.$$

Приравняем коэффициенты при $f_x$, $f_c$ и свободные члены к нулю. Это дает три уравнения:
$$\alpha=\xi_t+c\xi_x-c(\tau_t+c\tau_x),$$
$$\mathcal{F}_f(\eta+ f(\tau_t+c\tau_x - \alpha_c))+ f(\tau \mathcal{F}_{tf} + \xi\mathcal{F}_{xf} + \alpha\mathcal{F}_{cf} + \eta\mathcal{F}_{ff})=$$
$$=\alpha_t+c\alpha_x+\mathcal{F}\alpha_c- \mathcal{F}(\tau_t+c\tau_x)- (\tau\mathcal{F}_t + \xi\mathcal{F}_x + \alpha\mathcal{F}_c + \eta\mathcal{F}_{f}),$$
$$\eta_t+c\eta_x+F\eta_c+(\eta-f\eta_f)\mathcal{F}_c + \mathcal{F}_cf(\tau_t+c\tau_x)+ f(\tau \mathcal{F}_{tc} + \xi \mathcal{F}_{xc} + \alpha \mathcal{F}_{cc} + \eta\mathcal{F}_{cf})+ \mathcal{F}_f f\eta_c =0.$$
Первое из них выполнено в силу формулы для $\alpha$ в (23), во втором правая часть идентична левой части уравнения (18), и поэтому равна нулю. Наконец, левая часть второго уравнения и третье уравнение получаются дифференцированием по $c$ и по $f$ соответственно все того же уравнения (18) и использованием формулы для $\eta$ из (23).
Таким образом, и для группы симметрий оказывается, что условие инвариантности уравнения (5) является следствием выполнения условий сохранения (2), (3) и (6).

\subsubsection{Алгебра симметрий для $\mathcal{F}(t,x,c)$}
Нахождение алгебры симметрий для уравнения с $\mathcal{F}(t,x,c)$ дословно повторяет рассуждения предыдущего пункта, с той только разницей, что всюду $\mathcal{F}_f=0$. Таким образом, алгебра симметрий для уравнения (5) определяется формулами (23) и классифицирующим уравнением (19). Теорема 2 полностью доказана.

\section{Доказательство теоремы 4.}
\subsection{Схема доказательства}

Доказательство теоремы 4 мы будем проводить по той же схеме, что и доказательство результата, сформулированного нами в виде теоремы 3, в [3].
В силу теоремы 2 нас будет интересовать размерность пространства решений уравнения
$$\tau\mathcal{F}_t + \xi\mathcal{F}_x + (\xi_t + c(\xi_x - \tau_t) - c^2\tau_x)\mathcal{F}_c + f(3c\tau_x - 2\xi_x + \tau_t)\mathcal{F}_f +(3c\tau_x -\xi_x + 2\tau_t)\mathcal{F} =$$
$$= \xi_{tt} + c(2\xi_{tx} - \tau_{tt}) + c^2(\xi_{xx} - 2\tau_{tx}) - c^3\tau_{xx}. \eqno(27)$$
(это уравнение (18), в котором мы подставили выражения $\alpha$ и $\eta$ из (23), оставив только две неизвестные функции -- $\tau(t,x)$ и $\xi(t,x)$).

Прежде всего отметим, что приведение к канонической форме различных представителей одного подсемейства, осуществляемые с помощью указанных в таблицах 2-4 преобразований, проверяется по формулам (8) и не требует ничего, кроме навыков дифференцирования, поэтому мы эти выкладки приводить не будем, однако начнем доказательство, приведя список соответствующих преобразований.

Что же касается основной части доказательства, то сначала мы покажем, что пространство решений уравнения (27) конечномерно и, более того, не может иметь размерность более восьми. Затем мы воспользуемся классификацией вещественных конечномерных алгебр на плоскости из [15], из которой следует, что за исключением алгебры $\mathfrak{so}(3,\mathbb{R})$, все  алгебры размерности больше единицы содержат двумерные подалгебры, причем таких подалгебр, с точностью до диффеоморфизма (напомним, что диффеоморфизмы пространства $t,x$ -- это группа эквивалентности для нашего уравнения) ровно четыре. Это алгебры с базисами $(\partial_t, \partial_x)$, $(\partial_t, t\partial_t+x\partial_x)$, $(\partial_t,x\partial_t)$ и $(\partial_t,t\partial_t)$. Для каждой из этих алгебр мы установим соответствующий им анзац функций $\mathcal{F}$.

Это даст нам классификацию из таблицы 5. Для определения функций $\mathcal{F}$, для которых решение (27) дает более чем двумерные алгебры, следует провести разделение переменных для каждого из полученных анзацев и разобрать все возможные варианты. В результате получается довольно перепутанная система отношений между вариантами, поэтому мы будем обосновывать результат другим путем, более прозрачным для читателя.

Мы просто проверим, что решения (27) с приведенными в формулировке теоремы 4 функциями дают именно те алгебры, которые указаны в теореме, причем сделаем это "послойно". Сначала функции с пяти- и четырехмерной алгеброй, затем -- с трехмерной (проверяя каждый раз, что большая размерность возникает только в тех частных случаях, которые уже рассмотрены), а затем покажем, что для каждого из четырех основных анзацев разделение переменных приводит либо к одной из указанных функций, либо алгебра оказывается двумерной.

\subsection{Основные замены переменных, приводящих функцию $\mathcal{F}$ к канонической форме}

Как уже отмечалось, для приведения представителей семейства к канонической форме используются замены переменных, указанные в таблицах 2-4, и чисто техническую проверку того, что эти замены действительно преобразуют соответствующие функции к каноническому виду, мы опускаем. Для удобства читателя, мы лишь выпишем, как при этих подстановках меняются остальные переменные (фактически подставляя их в (8)). Обратим внимание, что ниже указывается лишь "внешнее" преобразование функции $\mathcal{F}$, а для полного осуществления замены переменных необходимо еще сделать соответствующие замены в аргументах у $\mathcal{F}$, перейдя от $(t,x,c,f)$ к $(\bar t,\bar x,\bar c,\bar f)$ по указанным ниже формулам.

Итак, замены
$$\bar t=\varphi(t),\qquad \bar x=x,\qquad
\bar c=\frac{c}{\varphi'(t)},\qquad
\bar f=f\varphi'(t), \qquad
\bar{\mathcal{F}}=\frac{1}{{\varphi'}^2(t)}\left[\mathcal{F}-c \frac{\varphi''(t)}{\varphi'(t)}\right],$$
$$\bar t=t,\qquad \bar x=\psi(x),\qquad
\bar c=c\psi'(x),\qquad
\bar f=\frac{f}{{\psi'}^3(x)}, \qquad
\bar{\mathcal{F}}=\psi'(x)\left[\mathcal{F}+c^2 \frac{\psi''(x)}{\psi'(x)}\right]$$
позволяют избавить функцию $\mathcal{F}$ от линейных и квадратичных по $c$ членов. Замена
$$\bar t=t,\qquad \bar x=x-\psi(t),\qquad
\bar c=c-\psi'(t),\qquad
\bar f=f, \qquad
\bar{\mathcal{F}}=\mathcal{F}- \psi''(t)$$
устраняет свободный член, а
$$\bar t=x,\qquad \bar x=t,\qquad
\bar c=\frac 1c,\qquad
\bar f=fc^3, \qquad
\bar{\mathcal{F}}=-\frac 1{c^3}\mathcal{F}$$
позволяет квадратичные и кубические члены превратить в линейные.

В более конкретных ситуациях используются следующие версии этих замен:
$$\bar t=e^{\lambda t},\qquad \bar x=x,\qquad
\bar c=\frac{c}{\lambda}e^{-\lambda t},\qquad
\bar f=\lambda fe^{\lambda t}, \qquad
\bar{\mathcal{F}}=\frac{e^{-2\lambda t}}{\lambda^2}[\mathcal{F}-\lambda c];$$
$$\bar t=t,\qquad \bar x=e^{\mu x},\qquad
\bar c=\mu ce^{\mu x},\qquad
\bar f=\frac{f}{\mu ^2}e^{-2\mu x}, \qquad
\bar{\mathcal{F}}=\mu e^{\mu x}[\mathcal{F}+\mu c^2 ];$$
$$\bar t=t^\alpha\ (\alpha\ne 0),\quad \bar x=x,\quad
\bar c=\frac{c}{\alpha}t^{1-\alpha},\quad
\bar f=\alpha f t^{\alpha-1}, \quad
\bar{\mathcal{F}}=\frac{t^{2-2\alpha}}{\alpha^2}\left[\mathcal{F}-c \frac{\alpha-1}{t}\right];$$
$$\bar t=\ln|t|,\qquad \bar x=x,\qquad
\bar c=ct,\qquad
\bar f=\frac ft, \qquad
\bar{\mathcal{F}}=t^2\left[\mathcal{F}+ \frac ct\right];$$
$$\bar t=x+kt,\qquad \bar x=x-kt,\qquad \bar c=\frac{c-k}{c+k},\qquad \bar f=f\frac{(k+c)^3}{4k^2},\qquad \bar{\mathcal{F}}=\frac{2k}{(k+c)^3}\mathcal{F};$$
$$\bar t=-\frac 1t,\qquad \bar x=\frac xt,\qquad \bar c=ct-x,\qquad \bar f=f,\qquad \bar{\mathcal{F}}=t^3\mathcal{F};$$
$$\bar t=e^{\mu (x+kt)},\qquad \bar x=e^{\mu (x-kt)},\qquad \bar c=e^{-2\mu kt}\frac{c-k}{c+k},\hspace{20mm}$$
$$\bar f=f\frac{e^{-\mu x+3k\mu t}(k+c)^3}{4k^2\mu },\qquad
\bar{\mathcal{F}}=\frac{2ke^{-\mu x-3\mu kt}}{\mu (k+c)^3}[\mathcal{F}-\mu (c^2-k^2)];$$
$$\bar t=e^{-k^2Qt}\cos kQ x,\qquad \bar x=e^{-k^2Qt}\sin kQ x,\qquad
\bar c=\frac{k\sin kQ x-c\cos kQx}{k\cos kQx+c\sin kQx},$$
$$\bar f=f\frac{-e^{k^2Qt}(k\cos kQx+c\sin kQx)^3}{k^3Q},\qquad \bar{\mathcal{F}}=\frac{e^{-k^2Q t}\left[\mathcal{F}-cQ(k^2+c^2)\right]}{Q(-k\cos kQx-c\sin kQx)^3};$$
$$\bar t=e^{Qx}\cos kQt,\qquad \bar x=e^{Qx}\sin kQt,\qquad
\bar c=\frac{k\cos kQt+c\sin kQt}{-k\sin kQt+c\cos kQt},$$
$$ \bar f=f\frac{(-k\sin kQt+c\cos kQt)^3}{k^2Q e^{Qx}},\qquad
\bar{\mathcal{F}}=\frac{k[-\mathcal{F}+Q(c^2+k^2)]}{e^{Qx}Q(-k\sin kQt+c\cos kQt)^3};$$
$$\bar t={\rm tg\,} kQt,\qquad \bar x=\frac{e^{-Qx}}{\cos kQt},\qquad \bar c =\frac 1k e^{-Qx}(-c\cos kQt+k\sin kQt),$$
$$\bar f=\frac kQ e^{2Qx}f,\qquad \bar{ \mathcal{F}}=\frac{e^{-Qx}\cos^3kQt}{k^2Q}[-\mathcal{F}+Q(c^2+k^2)];$$
$$\bar t=e^{-2kQt},\qquad \bar x=e^{-Q(x+kt)},\qquad \bar c =\frac 1{2k} e^{Q(-x+kt)}(c+k),$$
$$\bar f=-\frac {2k}Q e^{2Qx}f,\qquad \bar{ \mathcal{F}}=\frac{e^{-Qx+3kQt}}{4k^2Q}[-\mathcal{F}+Q(c^2-k^2)].$$
\subsection{Конечномерность алгебры симметрий}
Начнем с вопроса о конечномерности множества решений уравнения (27).
Предположим, что уравнение (27) с некоторой функцией $\mathcal{F}(t,x,c,f)$ имеет решение $(\tau(t,x),\xi(t,x))$ размерности не менее, чем семь. В этом случае мы можем написать как минимум семь уравнений вида (27) для каждого оператора $\Xi^i=\tau^i\partial_t+\xi^i\partial_x$, $i=1,\dots,7$. Полученную систему уравнений можно рассматривать как линейную однородную алгебраическую систему относительно семи переменных $\mathcal{F}_t$, $\mathcal{F}_x$, $\mathcal{F}_c$, $c\mathcal{F}_c - 2f\mathcal{F}_f - \mathcal{F}$, $c\mathcal{F}_c - f\mathcal{F}_f - 2\mathcal{F}$, $c^2\mathcal{F}_c - 3cf\mathcal{F}_f - 3c\mathcal{F}$ и единицы. Поскольку эта система имеет нетривиальное решение, ее определитель равен нулю:
$$\begin{vmatrix}
\tau^1 & \xi^1 & \xi^1_t & \xi^1_x & -\tau^1_t &- \tau^1_x  & -(\xi^1_{tt} + c(2\xi^1_{tx} - \tau^1_{tt}) + c^2(\xi^1_{xx} - 2\tau^1_{tx}) - c^3\tau^1_{xx}) \\  \vdots & \vdots & \vdots & \vdots & \vdots & \vdots & \vdots  \\
\tau^7 & \xi^7 & \xi^7_t & \xi^7_x & -\tau^7_t & - \tau^7_x & -(\xi^7_{tt} + c(2\xi^7_{tx} - \tau^7_{tt}) + c^2(\xi^7_{xx} - 2\tau^7_{tx}) - c^3\tau^7_{xx})
\end{vmatrix} = 0.
$$

В силу независимости $\tau^i$ и $\xi^i$ от $c$ мы получаем из этого условия четыре разных определителя, каждый из которых равен нулю.
Условие линейной зависимости столбцов каждого из этих определителей переписывается в виде дифференциального уравнения, которому должны удовлетворять все пары $(\tau^i,\xi^i)$ (причем, поскольку нетривиальность решения указанной системы гарантированно находится в последней компоненте, можно утверждать, что последний столбец выражается через остальные). Получим
$$\left\{\begin{array}{rcl}
\xi_{tt}&=&L^1_{11}(\tau,\xi),\\
2\xi_{tx}-\tau_{tt}&=&L^2_{11}(\tau,\xi),\\
\xi_{xx}-2\tau_{tx}&=&L^3_{11}(\tau,\xi),\\
\tau_{xx}&=&L^4_{11}(\tau,\xi),\\
\end{array}\right.\eqno(28)$$
через $L^k_{ij}(\tau,\xi)$ обозначены линейные дифференциальные операторы (с переменными, вообще говоря, коэффициентами), имеющие порядок $i$ производных от функций $\tau$ и порядок $j$ производных от функций $\xi$, индекс $k$ нумерует эти дифференциальные операторы.

Полученная система дифференциальных уравнений имеет только конечномерное пространство решений. В этом легко убедиться: воспользовавшись условиями согласования производных от $\xi$ в первых трех уравнениях, мы получаем два соотношения
$$\tau_{ttt}=L^5_{21}(\tau,\xi), \quad \tau_{ttx}=L^6_{21}(\tau,\xi),$$
добавление к ним соотношений
$$\tau_{txx}=L^7_{21}(\tau,\xi), \quad \tau_{xxx}=L^8_{21}(\tau,\xi),$$
получаемых дифференцированием последнего уравнения в (28) и первых трех уравнений из (28), дает нам линейную систему относительно $\tau$ и $\xi$ в нормальной форме, размерность пространства решений которой не превосходит количества начальных условий для функций $\tau$, $\xi$ и их производных: для $\tau$ до второго, а для $\xi$ -- до первого порядка. Общее количество таких начальных условий равно девяти. В рассматриваемом нами случае они зависимы, поскольку, в силу последнего уравнения (28), начальное значение для $\tau_{xx}$ определяется по значениям младших производных $\tau$ и $\xi$  однозначно. Значит, интересующая нас алгебра симметрий не может иметь размерность больше, чем восемь.

Как известно, все конечномерные алгебры Ли на плоскости были классифицированы С. Ли в [14], при этом алгебра рассматривалась как линейное пространство над полем комплексных чисел. Уточнение этой классификации для случая алгебры над полем вещественных чисел было сделано в [15]. Из представленного там результата очевидно вытекает, что все конечномерные  алгебры Ли плоскости над $\mathbb{R}$ либо одномерны, либо имеют двумерную подалгебру, либо эквивалентны алгебре вращений $\mathfrak{so}(3,\mathbb{R})$.  Поэтому нам необходимо рассмотреть три основных случая: одномерной алгебры, двумерных алгебр и алгебры $\mathfrak{so}(3,\mathbb{R})$.

\subsection{Случай алгебры $\mathfrak{so}(3,\mathbb{R})$}
Нам удобно использовать следующую реализацию алгебры $\mathfrak{so}(3)$: $\Xi_1=-x\partial_t+t\partial_x$, $\Xi_2 = tx\partial_t + (x^2+1)\partial_x$ и $\Xi_3=(t^2+1)\partial_t + tx\partial_x$, эта алгебра имеет то преимущество, что для нее уравнение (27) оказывается однородным.

Найдем, какие функции $\mathcal{F}(t,x,c,f)$ обладают данной алгеброй. Подставим операторы $\Xi_1$, $\Xi_2$ и $\Xi_3$ в $(27)$, получим систему:


$$-x\mathcal{F}_t + t\mathcal{F}_x + (1 + c^2)\mathcal{F}_c - 3cf\mathcal{F}_f - 3c\mathcal{F} = 0,$$
$$tx\mathcal{F}_t + (x^2+1)\mathcal{F}_x + c(x - ct)\mathcal{F}_c - 3f(x - ct)\mathcal{F}_f + 3ct\mathcal{F} = 0,$$
$$(t^2 + 1)\mathcal{F}_t + tx\mathcal{F}_x + (x - ct) \mathcal{F}_c + 3t\mathcal{F} = 0.$$

Линейным комбинированием эта система приводится к виду
$$-x\mathcal{F}_t + t\mathcal{F}_x + (1 + c^2)\mathcal{F}_c - 3cf\mathcal{F}_f - 3c\mathcal{F} = 0,$$
$$(t^2 + x^2 + 1)\mathcal{F}_x + (t + cx)\mathcal{F}_c - 3xf\mathcal{F}_f = 0,$$
$$(t^2 + x^2 + 1)\mathcal{F}_t + c(t^2 + x^2 + 1)\mathcal{F}_x + 3(t+cx)\mathcal{F} = 0.$$

Из последнего уравнения получаем $\mathcal{F}=\Phi(x-ct,c,f)(1+x^2+t^2)^{-3/2}$,
подстановка этого выражения в первое уравнение дает соотношение
$$ca\Phi_a + (1 + c^2)\Phi_c - 3cf\Phi_f - 3c\Phi = 0,$$
где через $a$ обозначен первый аргумент функции $\Phi$. Отсюда мы получаем, что
$\Phi(a,c,f)=\frac 1f\Psi(fa^3, a^2/(1+c^2))$, а значит,
$$\mathcal{F} = \frac 1f (1+x^2+t^2)^{-\frac32}\Psi(f(x-ct)^3, \frac{(x-ct)^2}{1+c^2}).$$
Наконец, подстановка полученного выражения в среднее уравнение приводит уже к дифференциальному уравнению на $\Psi$: $3b\Psi_b + 2g(1+g)\Psi_g = 0$, где $b = f(x-ct)^3$, $g = \frac{(x-ct)^2}{1+c^2}$. Откуда находим $\Psi = G\left(b(\frac{1+g}{g})^{\frac32}\right)$, где $G$ -- произвольная функция, что и дает
$$\mathcal{F} = \frac 1{(1+x^2+t^2)^{\frac32}f} G\left(((x-ct)^2 + 1+c^2)^{\frac32}f\right),$$
или, эквивалентно,
$$\mathcal{F} = \left(\frac {(x-ct)^2 + 1+c^2}{1+x^2+t^2}\right)^{\frac32} G\left(((x-ct)^2 + 1+c^2)^{\frac32}f\right),$$
указанную в формулировке теоремы 4 под номером III.1.

Обратно, вычисление производных функции $\mathcal{F}$ указанного вида
$$\mathcal{F}_t=\frac 32 \left[\frac{-2c(x-ct)}{(x-ct)^2 + 1+c^2}-\frac{2t}{1+x^2+t^2}\right]\mathcal{F}+$$
$$+\frac 32\left(\frac {(x-ct)^2 + 1+c^2}{1+x^2+t^2}\right)^{\frac32} G'\left(((x-ct)^2 + 1+c^2)^{\frac32}f\right)[-2c(x-ct)((x-ct)^2 + 1+c^2)^{\frac12}f],$$
$$\mathcal{F}_x=\frac 32\left[\frac{2(x-ct)}{(x-ct)^2 + 1+c^2}-\frac{2x}{1+x^2+t^2}\right]\mathcal{F}+$$
$$+\frac 32\left(\frac {(x-ct)^2 + 1+c^2}{1+x^2+t^2}\right)^{\frac32} G'\left(((x-ct)^2 + 1+c^2)^{\frac32}f\right)[2(x-ct)((x-ct)^2 + 1+c^2)^{\frac12}f],$$
$$\mathcal{F}_c=\frac 32\frac{-2t(x-ct)+2c}{(x-ct)^2 + 1+c^2}\mathcal{F}+$$
$$+\frac 32\left(\frac {(x-ct)^2 + 1+c^2}{1+x^2+t^2}\right)^{\frac32} G'\left(((x-ct)^2 + 1+c^2)^{\frac32}f\right)[-2t(x-ct)+2c]((x-ct)^2 + 1+c^2)^{\frac12}f,$$
$$\mathcal{F}_f=\left(\frac {(x-ct)^2 + 1+c^2}{1+x^2+t^2}\right)^{\frac32} G'\left(((x-ct)^2 + 1+c^2)^{\frac32}f\right)((x-ct)^2 + 1+c^2)^{\frac32}$$
и подстановка их в (27), после введения обозначения $z=((x-ct)^2 + 1+c^2)^{\frac32}f$ приводит к уравнению
$$\frac 32\left(\frac {(x-ct)^2 + 1+c^2}{1+x^2+t^2}\right)^{\frac32}\left\{\tau \left[G(z)\left(\frac{-2c(x-ct)}{(x-ct)^2 + 1+c^2}-\frac{2t}{1+x^2+t^2}\right)+
 G'(z)z\frac{-2c(x-ct)}{(x-ct)^2 + 1+c^2}\right]\right.+$$
$$\xi\left[G(z)\left(\frac{2(x-ct)}{(x-ct)^2 + 1+c^2}-\frac{2x}{1+x^2+t^2}\right)+
 G'(z)z\frac{2(x-ct)}{(x-ct)^2 + 1+c^2}\right]+$$
$$+(\xi_t + c(\xi_x - \tau_t) - c^2\tau_x)\left[G(z)\frac{-2t(x-ct)+2c}{(x-ct)^2 + 1+c^2}+G'(z)z\frac{-2t(x-ct)+2c}{(x-ct)^2 + 1+c^2}\right]
+\frac 23 (3c\tau_x - 2\xi_x + \tau_t) G'(z)z$$
$$\left.+\frac 23(3c\tau_x -\xi_x + 2\tau_t)G(z)\right\}=\xi_{tt} + c(2\xi_{tx} - \tau_{tt}) + c^2(\xi_{xx} - 2\tau_{tx}) - c^3\tau_{xx}. $$
Левая часть этого уравнения зависит от $c$ иррациональным образом, правая же является полиномом, поэтому и левая, и правая части равны нулю по отдельности.

Для правой части отсюда следует
$$\xi_{tt} = 2\xi_{tx} - \tau_{tt}=\xi_{xx} - 2\tau_{tx}=\tau_{xx}=0,$$
решение полученных уравнений имеет вид $$\tau=At^2+Btx+Ct+Dx+E, \qquad \xi=Atx+Bx^2+Ft+Gx+H.\eqno(29)$$

А для левой после умножения на $(x-ct)^2+1+c^2$ мы получаем, что
$\frac{(t\tau+x\xi)[(x-ct^2)+1+c^2]}{1+x^2+t^2}$ должно делиться нацело.
Поскольку
$\frac{(x-ct^2)+1+c^2}{1+x^2+t^2}=1+c^2-\frac{(t+cx)}{1+x^2+t^2}$, и получившаяся дробь несократима, нам остается только утверждать, что
$t\tau+x\xi$ должно нацело делиться на $1+t^2+x^2$.
Подставляя в это условие формулы (29), получаем, что нацело делится дробь
$$\frac{t\tau+x\xi}{1+t^2+x^2}=(At+Bx+C)+
\frac{(G-C)x^2+(D+F)tx+(E-A)t+(H-B)x-C}{1+t^2+x^2},$$
откуда $G=C=0$, $D=-F$, $E=A$, $H=B$, то есть мы получили ту самую трехмерную алгебру, без дополнительных расширений.

\subsection{Семейства функций $\mathcal {F}$ для уравнений с одно- и двумерными алгебрами симметрий}

Как известно, любую нетривиальную одномерную алгебру в $\mathbb{R}$ диффеоморфизмом можно превратить в $\Xi=\partial_x$. Условие $\Xi \mathcal{F} = 0$ дает нам условие независимости $\mathcal {F}$ от переменной $x$, и очевидно, что для всех таких функций уравнение (5)  обладает указанной алгеброй. Таким образом мы получаем семейство V из утверждения теоремы. Впрочем, в каких-то ситуациях, возможно, удобнее, чтобы этот оператор был дифференцированием по временн\'{о}й переменной, и тогда мы получаем другую форму для представителей соответствующих классов эквивалентности -- не зависящую уже от $t$.

Как уже отмечалось выше, любая конечномерная алгебра симметрий размерности больше двух, кроме $\mathfrak{so}(3)$, содержит двумерную подалгебру. Поэтому ниже мы найдем все уравнения, допускающие двумерную алгебру симметрий, и уже среди них будем отыскивать уравнения с алгебрами б\'{о}льших размерностей.

Любую двумерную алгебру $\Xi = A\Xi_1 + B\Xi_2$ за счет выбора подходящего базиса можно привести к одному из двух канонических случаев -- $[\Xi_1,\Xi_2] = 0$ (коммутативная алгебра) и $[\Xi_1, \Xi_2] = \Xi_1$ (некоммутативная алгебра), мы далее рассмотрим каждый их этих двух случаев.

Итак, пусть сначала алгебра коммутативна. Тогда она либо приводится к виду $\Xi=A\partial_x+B\partial_t$, либо к виду $\Xi=A\partial_x+Bt\partial_x$. Действительно, один из операторов заменой переменных из группы эквивалентности можно привести к виду $\Xi_1=\partial_x$, при этом соответствующая замена $\bar{t} = \varphi(t,x)$, $\bar{x} = \psi(t,x)$ находится как решение системы $\Xi_1\varphi= 0$, $\Xi_1\psi= 1$.
Тогда из условий коммутирования c $\partial_x$ для оператора $\Xi_2=\tau_2\partial_t+\xi_2\partial_x$ получаем $\tau_2 = \tau_2(t)$, $\xi_2 = \xi_2(t)$.

В случае $\tau_2 \ne 0$  мы можем найти замену переменных, которая переведет $(\tau_2,\xi_2)$ в $(1,0)$ и сохранит первое векторное поле: $\bar t = \int 1/\tau_2(t)dt$, $\bar x = x - \int \xi_2(t)/\tau_2(t) dt$.
В случае $\tau_2=0$ обозначение $\bar t=\xi_2(t)$ приводит второй оператор к виду $t\partial_x$.
Первому случаю -- транзитивной коммутативной группы -- соответствует $\mathcal{F} =\Psi(c,f)$, а второму -- нетранзитивной коммутативной группы -- $\mathcal{F}=\Psi(t,f)$.

Если алгебра некоммутативна, то, приведя первый оператор к виду $\partial_x$, мы для второго оператора из условия $[\Xi_1, \Xi_2] = \Xi_1$
получаем $(\tau_2)_x=0$, $(\xi_2)_x=1$, то есть $\tau_2=\tau(t)$, $\xi_2=x+\tilde\xi(t)$.

Здесь тоже есть два варианта. Если $\tau(t)\ne 0$, то можно, сделав замену $\bar t=\varphi(t)$, $\bar x=x-\psi(t)$, где функции $\varphi$ и $\psi$ определяются из уравнений $\tau\dot\varphi=\varphi$, $\tau\dot\psi=\psi-\tilde\xi(t)$, мы сведем второй оператор к виду $\Xi_2=t\partial_t+x\partial_x$.
Если же $\tau(t)=0$, то замена $\bar t=t$, $\bar x=x+\tilde \xi(t)$ превращает второй оператор в $x\partial_x$.

В первом случае (транзитивная группа) подстановка нашей алгебры в (27) дает $\mathcal{F}=\frac 1t\Psi(c,tf)$, во втором (нетранзитивная группа) -- $\mathcal{F}=c\Psi(t,c^2f)$

Четыре полученных семейства и есть семейства IV.1-IV.4, приведенные в формулировке теоремы 4 в таблице 5. Очевидно, что для любой функции $\mathcal{F}$ этих семейств соответствующая группа симметрий будет содержать как минимум отвечающую этому семейству двумерную алгебру. Для некоторых функций из этих семейств алгебра симметрий может оказаться шире (эти функции приведены в таблицах 2-4). В следующих параграфах мы проверим, что этим функциям отвечают именно указанные алгебры, и покажем, что других функций, для которых алгебра симметрий уравнения (5) имеет размерность более, чем два, нет.

\subsection{Группы эквивалентности семейств}

Для дальнейшего анализа нам будет удобно вычислить также алгебры эквивалентности каждого из четырех выделенных нами семейств функций $\mathcal{F}$.

{\bf Лемма 1.} {\em Алгебры эквивалентности для выделенных нами семейств функций $\mathcal{F}$ имеют вид, приведенный в таблице 6:
\begin{figure}[h]
Таблица 6. Алгебры эквивалентности семейств IV.
\begin{center}
\begin{tabular}{|c|l|}
\hline
 &\\
\text{Семейство функций }$\mathcal{F}(t,x,c,f)$ & \text{          }Базис алгебры эквивалентности\\
\hline
\hline
&\\
$\mathcal{F} = \Psi(c,f)$ & $\Xi_1= \partial_t$, $\Xi_2= \partial_x$, $\Xi_3=t\partial_t $,\\
&  $\Xi_4 = x\partial_x$, $\Xi_5 = x\partial_t$, $\Xi_6 = t\partial_x$\\
\hline
&\\
$\mathcal{F} = \Psi(t,f)$ & $\Xi_1= \partial_t$, $\Xi_2= t\partial_t$, $\Xi_3=t^2\partial_t+tx\partial_x$,\\
& $\Xi_4 = x\partial_x$, $\Xi_\infty=h(t)\partial_x,$ \\
\hline
&\\
$\mathcal{F} = \frac 1t\Psi(c,tf)$ & $\Xi_1= t\partial_t$, $\Xi_2= \partial_x$, $\Xi_3=x\partial_x $, $\Xi_4 = t\partial_x $\\
&\\
\hline
$\mathcal{F} = c\Psi(t,c^2f)$& $\Xi_1= \partial_x$, $\Xi_2= x\partial_x$, $\Xi_\infty=h(t)\partial_t$.\\
\hline
\end{tabular}
\end{center}
\end{figure}}

{\bf Доказательство.} Воспользуемся формулой (7) из теоремы 1, дающей нам компоненту
$$\Xi\mathcal{F}=\phi = \xi_{tt} + c(2\xi_{tx} - \tau_{tt}) + c^2(\xi_{xx} - 2\tau_{tx}) - c^3\tau_{xx} + \mathcal{F}(\xi_x - 2\tau_t - 3c\tau_x),$$
найдем продолжения оператора на производные функции $\mathcal{F}$, и заметим, что каждое семейство определяется парой дифференциальных соотношений. Это соотношения $\mathcal{F}_t=\mathcal{F}_x=0$ для первого семейства, $\mathcal{F}_x=\mathcal{F}_c=0$ для второго, $\mathcal{F}_x=t\mathcal{F}_t-f\mathcal{F}_f+\mathcal{F}=0$ для третьего и $\mathcal{F}_x=c\mathcal{F}_c-2f\mathcal{F}_f-\mathcal{F}=0$ для четвертого. Поэтому нам достаточно найти те $\phi$, для которых остаются инвариантными эти пары соотношений.

Продолжение действия оператора на производные имеет вид
$$\Xi\mathcal{F}_t=\xi_{ttt} + c(2\xi_{ttx} - \tau_{ttt}) + c^2(\xi_{txx} - 2\tau_{ttx}) - c^3\tau_{txx} + \mathcal{F}(\xi_{tx} - 2\tau_{tt} - 3c\tau_{tx})+\mathcal{F}_t(\xi_x - 3\tau_t - 3c\tau_x)-$$
$$-\mathcal{F}_x\xi_t-\mathcal{F}_c (\xi_{tt} + c\xi_{tx} - c(\tau_{tt} + c\tau_{tx}))-f\mathcal{F}_f(3c\tau_{tx} - 2\xi_{tx} + \tau_{tt}) ,$$
$$\Xi\mathcal{F}_x=\xi_{ttx} + c(2\xi_{txx} - \tau_{ttx}) + c^2(\xi_{xxx} - 2\tau_{txx}) - c^3\tau_{xxx} + \mathcal{F}(\xi_{xx} - 2\tau_{tx} - 3c\tau_{xx})+\mathcal{F}_x(- 2\tau_t - 3c\tau_x)-$$
$$-\mathcal{F}_t\tau_x-\mathcal{F}_c (\xi_{tx} + c\xi_{xx} - c(\tau_{tx} + c\tau_{xx}))-f\mathcal{F}_f(3c\tau_{xx} - 2\xi_{xx} + \tau_{tx}),$$
$$\Xi\mathcal{F}_c=2\xi_{tx} - \tau_{tt} + 2c(\xi_{xx} - 2\tau_{tx})-3c^2\tau_{xx}-3\mathcal{F}\tau_{x}+
\mathcal{F}_c(- \tau_t - c\tau_x)-3f\mathcal{F}_f \tau_{x},$$
$$\Xi\mathcal{F}_f=\mathcal{F}_f(3\xi_x - 3\tau_t - 6c\tau_x).$$

Условие инвариантности соотношений $\mathcal{F}_t=\mathcal{F}_x=0$ дает равенства
$$\xi_{ttt} + c(2\xi_{ttx} - \tau_{ttt}) + c^2(\xi_{txx} - 2\tau_{ttx}) - c^3\tau_{txx} + \mathcal{F}(\xi_{tx} - 2\tau_{tt} - 3c\tau_{tx})-$$
$$-\mathcal{F}_c (\xi_{tt} + c\xi_{tx} - c(\tau_{tt} + c\tau_{tx}))-f\mathcal{F}_f(3c\tau_{tx} - 2\xi_{tx} + \tau_{tt})=0,$$
$$\xi_{ttx} + c(2\xi_{txx} - \tau_{ttx}) + c^2(\xi_{xxx} - 2\tau_{txx}) - c^3\tau_{xxx} + \mathcal{F}(\xi_{xx} - 2\tau_{tx} - 3c\tau_{xx})-$$
$$-\mathcal{F}_c (\xi_{tx} + c\xi_{xx} - c(\tau_{tx} + c\tau_{xx}))-f\mathcal{F}_f(3c\tau_{xx} - 2\xi_{xx} + \tau_{tx}),$$
из которых, в силу произвольности величин $\mathcal{F}_c$ и $\mathcal{F}_f$, следует, что все вторые производные функций $\tau$ и $\xi$ равны нулю. Что и дает нам алгебру линейных преобразований, указанную в первой строке.

Условие инвариантности соотношений $\mathcal{F}_x=\mathcal{F}_c=0$ имеет вид
$$\xi_{ttx} + c(2\xi_{txx} - \tau_{ttx}) + c^2(\xi_{xxx} - 2\tau_{txx}) - c^3\tau_{xxx} + \mathcal{F}(\xi_{xx} - 2\tau_{tx} - 3c\tau_{xx})-$$
$$-\mathcal{F}_t\tau_x-f\mathcal{F}_f(3c\tau_{xx} - 2\xi_{xx} + \tau_{tx})=0,$$
$$2\xi_{tx} - \tau_{tt} + 2c(\xi_{xx} - 2\tau_{tx})-3c^2\tau_{xx}-3\mathcal{F}\tau_{x}-3f\mathcal{F}_f \tau_{x}=0,$$
откуда следует $\tau_x=\xi_{xx}=0$, $\xi_{ttx} =0$, $2\xi_{tx} - \tau_{tt} =0$, что и приводит к алгебре, указанной во второй строчке. В отличие от первой, она бесконечномерна и позволяет из функции $\mathcal{F}=\Psi(t,x)$ вычитать любую функцию от $t$.

Условие инвариантности соотношений $\mathcal{F}_x=t\mathcal{F}_t-f\mathcal{F}_f+\mathcal{F}=0$ имеет вид
$$\xi_{ttx} + c(2\xi_{txx} - \tau_{ttx}) + c^2(\xi_{xxx} - 2\tau_{txx}) - c^3\tau_{xxx} + \mathcal{F}(\xi_{xx} - 2\tau_{tx} - 3c\tau_{xx})-$$
$$-\mathcal{F}_t\tau_x-\mathcal{F}_c (\xi_{tx} + c\xi_{xx} - c(\tau_{tx} + c\tau_{xx}))-(t\mathcal{F}_t+\mathcal{F})(3c\tau_{xx} - 2\xi_{xx} + \tau_{tx})=0,$$
$$t[\xi_{ttt} + c(2\xi_{ttx} - \tau_{ttt}) + c^2(\xi_{txx} - 2\tau_{ttx}) - c^3\tau_{txx} + \mathcal{F}(\xi_{tx} - 2\tau_{tt} - 3c\tau_{tx})+\mathcal{F}_t(\xi_x - 3\tau_t - 3c\tau_x)-$$
$$-\mathcal{F}_x\xi_t-\mathcal{F}_c (\xi_{tt} + c\xi_{tx} - c(\tau_{tt} + c\tau_{tx}))-f\mathcal{F}_f(3c\tau_{tx} - 2\xi_{tx} + \tau_{tt})]+\tau \mathcal{F}_t-$$
$$-(t\mathcal{F}_t+\mathcal{F})(\xi_x - 2\tau_t - 3c\tau_x)+\xi_{tt} + c(2\xi_{tx} - \tau_{tt}) + c^2(\xi_{xx} - 2\tau_{tx}) - c^3\tau_{xx} + \mathcal{F}(\xi_x - 2\tau_t - 3c\tau_x)=0,$$
откуда, в силу произвольности $\mathcal{F}_c$ немедленно следует, что все вторые производные от $\tau$ и $\xi$ равны нулю, а тогда условия редуцируются к $\mathcal{F}_t\tau_x=\mathcal{F}_t[\tau -t(\xi_x - 2\tau_t - 3c\tau_x)]=0$, откуда $\tau_x=0$, $\tau=t\tau_t$, и линейность обеих функций дает нам алгебру в третьей строчке.


Наконец, для соотношений $\mathcal{F}_x=c\mathcal{F}_c-2f\mathcal{F}_f-\mathcal{F}=0$ условие инвариантности приобретает вид
$$\xi_{ttx} + c(2\xi_{txx} - \tau_{ttx}) + c^2(\xi_{xxx} - 2\tau_{txx}) - c^3\tau_{xxx} + \mathcal{F}(\xi_{xx} - 2\tau_{tx} - 3c\tau_{xx})-$$
$$-\mathcal{F}_t\tau_x-\mathcal{F}_c (\xi_{tx} + c\xi_{xx} - c(\tau_{tx} + c\tau_{xx}))-\frac 12(c\mathcal{F}_c-\mathcal{F})(3c\tau_{xx} - 2\xi_{xx} + \tau_{tx})=0,$$
$$c[2\xi_{tx} - \tau_{tt} + 2c(\xi_{xx} - 2\tau_{tx})-3c^2\tau_{xx}-3\mathcal{F}\tau_{x}+
\mathcal{F}_c(- \tau_t - c\tau_x)-\frac 32(c\mathcal{F}_c-\mathcal{F})\tau_{x}]+\mathcal{F}_c(\xi_t + c(\xi_x - \tau_t) -c^2\tau_x))-$$
$$-(c\mathcal{F}_c-\mathcal{F})[\xi_x - 2\tau_t - 3c\tau_x)-
[\xi_{tt} + c(2\xi_{tx} - \tau_{tt}) + c^2(\xi_{xx} - 2\tau_{tx}) - c^3\tau_{xx} + \mathcal{F}(\xi_x - 2\tau_t - 3c\tau_x)]=0,$$
откуда, в силу произвольности $\mathcal{F}_t$, $\mathcal{F}_c$, получаем $\tau_x=\xi_t=0=0$, $\xi_{xx}= 0$, что и приводит к алгебре в четвертой строке. Лемма доказана.

\subsection{Алгебры симметрий для семейств I и II.*}
В этом параграфе мы покажем, что уравнения (5) с функциями семейств I и II.* имеют в точности пяти- и четырехмерную алгебру симметрий соответственно. Для этого мы каждую из этих функций подставим в классифицирующее уравнение (27) и найдем соответствующую алгебру.
\subsubsection{$\mathcal{F}={P}/{f}$}
Подставляя эту функцию в (27), получаем
$$(\xi_x + \tau_t)\frac{P}{f}=\xi_{tt} + c(2\xi_{tx} - \tau_{tt}) + c^2(\xi_{xx} - 2\tau_{tx}) - c^3\tau_{xx}.$$
Расщепление сначала по $f$, а потом по $c$ дает $\xi_x+\tau_t=0$ и равенство нулю всех вторых производных, откуда мы и получаем пятимерную алгебру, указанную в I.

\subsubsection{$\mathcal{F}=Pf^j$}
Подстановка в (27) дает
$$[(3c\tau_x - 2\xi_x + \tau_t)j +3c\tau_x -\xi_x + 2\tau_t]Pf^j =
\xi_{tt} + c(2\xi_{tx} - \tau_{tt}) + c^2(\xi_{xx} - 2\tau_{tx}) - c^3\tau_{xx},$$
откуда после расщепления
$$3(j+1)\tau_x=0,\qquad -(2j+1)\xi_x + (j+2)\tau_t=0,\qquad
\xi_{tt}=2\xi_{tx} - \tau_{tt}=\xi_{xx} - 2\tau_{tx}=\tau_{xx}=0.$$
При $j\ne -1$ мы имеем $\tau_x=0$, и решение системы дает алгебру, указанную в II.1, при $j=-1$ мы получаем уже рассмотренную выше функцию I.

\subsubsection{$\mathcal{F}=P(c^2\pm k^2)^3f$}
Для этой функции подстановка в (27) дает
$$[(\xi_t + c(\xi_x - \tau_t) - c^2\tau_x)6c + (6c\tau_x - 3\xi_x + 3\tau_t)(c^2\pm k^2)]P(c^2\pm k^2)^2f=\hspace{50mm}$$
$$\hspace{70mm}=\xi_{tt} + c(2\xi_{tx} - \tau_{tt}) + c^2(\xi_{xx} - 2\tau_{tx}) - c^3\tau_{xx},$$
расщеплением по $f$, а затем по $c$ получаем $\xi_x=\tau_t$, $\xi_t\pm k^2\tau_x=0$, $\tau_{tt}=\tau_{xx}= \tau_{tx}=0$, что и дает нам алгебру II.2.

\subsubsection{$\mathcal{F}=Pf^{-2}x^{-3}$}
Аналогично подстановка и расщепление приводят к системе
$\tau_x=\xi_x-\frac \xi x =0$, $\tau_{tt}-2\xi_{tx}=\xi_{tt}=0$, которая дает алгебру II.3.

\subsubsection{$\mathcal{F}=P(t^2\pm k^2)^{-3/2}f^{-1/2}$}
И снова расщепление дает систему
$\tau_x=0$, $\tau_t(t^2\pm k^2)=2t\tau$,
$\xi_{tt} =2\xi_{tx} - \tau_{tt}=\xi_{xx}=0$, которая дает алгебру II.4.

\subsubsection{$\mathcal{F}=P\ln |f|+Q$}
Все так же расщепляя уравнение, получаем
$\tau_x=0$, $\xi_x = 2\tau_t$,
$\xi_{tt}=- 3P\tau_t$, $\tau_{tt}=0$, откуда немедленно приходим к алгебре II.5.


\subsection{Алгебры симметрий для семейств III.*.*}

В этом параграфе мы докажем, что для всех функций семейств III.*.*, за исключением тех, которые попадают в семейства I и II.*, алгебра симметрий уравнения (5) в точности трехмерна. При этом, помимо прямого выхода на случаи II.*.*, мы будем использовать для приведения к этим функциям некоторые простые замены из группы эквивалентности . Это замена $\bar t=x$, $\bar x=t$, для нее $\bar c=1/c$, $\bar {\mathcal{F}}=-c^{-3}\mathcal{F}$, $\bar f=c^3f$, и замена $\bar t=t$, $\bar x=x-kt$, для которой $\bar c=c-k$, $\bar{\mathcal{F}}=\mathcal{F}$, $\bar f=f$.

Случай III.1 с алгеброй $\mathfrak{so}(3)$ был рассмотрен выше, поэтому мы перейдем к следующему случаю.

\subsubsection{$\mathcal{F} = G(f)$}
Подстановка этой функции в (27) приводит к равенству
$$ (3c\tau_x - 2\xi_x + \tau_t)fG'(f) +(3c\tau_x -\xi_x + 2\tau_t)G(f) = \xi_{tt} + c(2\xi_{tx} - \tau_{tt}) + c^2(\xi_{xx} - 2\tau_{tx}) - c^3\tau_{xx}.$$
Если функции $fG'(f)$, $G(f)$ и единица являются линейно независимыми, то из этого равенства следует
$$3c\tau_x - 2\xi_x + \tau_t=3c\tau_x -\xi_x + 2\tau_t= \xi_{tt} + c(2\xi_{tx} - \tau_{tt}) + c^2(\xi_{xx} - 2\tau_{tx}) - c^3\tau_{xx}=0,$$
расщепление полученных равенств по $c$ дает
$\tau_x=\xi_x=\tau_t=0$, $\xi_{tt}=0$, откуда мы и получаем алгебру III.2.1.

Если же эти функции линейно зависимы, то учитывая, что $G(f)$ не является константой по предположению теоремы, мы получаем, что $fG'(f)=\alpha G(f)+\beta$, где $\alpha$ и $\beta$ -- некоторые константы. Если $\alpha\ne 0$, то $\mathcal{F}=G(f)=Pf^{\alpha}-\frac \beta\alpha$, замена
$\bar t=t$, $\bar x=x+\frac {\beta t^2}{2\alpha}$ приводит эту функцию к виду II.1 при $\alpha\ne -1$ или I при $\alpha=-1$. Если же $\alpha=0$, то мы получаем функцию II.5.

\subsubsection{$\mathcal{F} ={G(f)}/{x^3}$}
Подставляя эту функцию в (27), получаем уравнение
$$-3\xi\frac{G(f)}{x^4} + (3c\tau_x - 2\xi_x + \tau_t)\frac{fG'(f)}{x^3} +(3c\tau_x -\xi_x + 2\tau_t)\frac{G(f)}{x^3} = \xi_{tt} + c(2\xi_{tx} - \tau_{tt}) + c^2(\xi_{xx} - 2\tau_{tx}) - c^3\tau_{xx}. $$
И здесь, если $fG'(f)$, $G(f)$ и единица являются линейно независимыми, то мы получаем
$$3c\tau_x - 2\xi_x + \tau_t=0,\quad -3\xi\frac 1x + 3c\tau_x -\xi_x + 2\tau_t=0,\quad \xi_{tt} + c(2\xi_{tx} - \tau_{tt}) + c^2(\xi_{xx} - 2\tau_{tx}) - c^3\tau_{xx}=0,$$
расщепление по $c$ дает нам
$\tau_x=0$, $2\xi_x=\tau_t$, $-3\xi\frac 1x +3\xi_x=0$,
$\xi_{tt}=\xi_{xx}=0$, откуда мы немедленно получаем алгебру III.2.2.

Если же $fG'(f)=\alpha G(f)+\beta$, то мы получаем другое расщепление по $f$ (поскольку $G(f)$ по условиям теоремы непостоянная):
$$-3\xi\frac{1}{x} + 3(\alpha+1) c\tau_x - (2\alpha+1)\xi_x + (\alpha+2)\tau_t=0,$$
$$\frac {\beta}{x^3} (3c\tau_x - 2\xi_x + \tau_t) = \xi_{tt} + c(2\xi_{tx} - \tau_{tt}) + c^2(\xi_{xx} - 2\tau_{tx}) - c^3\tau_{xx},$$
и дальнейшее расщепление по $c$ приводит к системе
$$3(\alpha+1) \tau_x=0,\qquad -3\frac{\xi}{x} - (2\alpha+1)\xi_x + (\alpha+2)\tau_t=0, $$
$$\xi_{tt}=\frac{\beta}{x^3} ( - 2\xi_x + \tau_t),\qquad 2\xi_{tx} - \tau_{tt}=\frac {3\beta}{x^3}\tau_x, \qquad \xi_{xx} - 2\tau_{tx}=\tau_{xx}=0.$$
Если $\alpha\ne -1$, то мы из первого уравнения получаем $\tau_x=0$, из последнего -- $\xi_{xx}=0$, а значит, $\xi=p(t)x+q(t)$. Подстановка этой функции во второе уравнение дает
$(\alpha+2)(\tau_t-2p(t))=0$, $q(t)=0$, из третьего $p''(t)=\beta(\tau_t-2p(t))=0$, из четвертого $2p'(t)=\tau_{tt}$. В итоге $p(t)=at+b$, $\xi=(at+b)x$, $\tau=at^2+kt+l$ и $(\alpha+2)(k-2b)=\beta(k-2b)=0$.
Если $k=2b$, то мы получаем всю ту же трехмерную алгебру III.2.2, случай же $k\ne 2b$ возможен только при $\alpha=-2$, $\beta=0$, что дает нам функцию II.3.

\subsubsection{$\mathcal{F}=G(c)f$}
Подстановка в (27) и расщепление по $f$ приводит нас к условиям
$$(\xi_t + c(\xi_x - \tau_t) - c^2\tau_x)G'(c) + 3(2c\tau_x - \xi_x + \tau_t)G(c)=0,$$
$$\xi_{tt} =2\xi_{tx} - \tau_{tt}=\xi_{xx} - 2\tau_{tx}= \tau_{xx}=0.$$
Из второго мы получаем для $\tau$ и  $\xi$ формулы
$$\tau=At^2+Btx+Ct+Dx+E,\qquad \xi=Atx+Bx^2+Ft+Gx+H,$$
подстановка которых в первое уравнение и расщепление по переменным $t,x$ дает равенства
$$-(Ac+Bc^2)G'(c) + 3(2Bc+A)G(c)=0,\qquad
(A+cB)G'(c) - 3BG(c)=0, $$
$$(F+c(G-C) - c^2D)G'(c) + 3(2Dc+C-G)G(c)=0.$$
Складывая второе равенство, умноженное на $c$, с первым, мы получаем $G(c)(Bc+A)=0$, откуда, в силу $G(c)\ne 0$, $A=B=0$. Из третьего же равенства мы получаем, что либо $F=D=0$, $G=C$, что дает нам трехмерную алгебру, указанную в III.2.3, либо $G(c)=(Pc^2+Qc+R)^3$, $\mathcal{F}=(Pc^2+Qc+R)^3f$.

В последнем случае, если $P\ne 0$, то заменой переменных $\bar t=t$, $\bar x=x+\frac {Q}{2P} t$ мы приводим нашу функцию к виду $\mathcal{F}=P(c^2\pm k^2)^3f$, то есть к функции II.2, если $P=0$, но $Q\ne 0$, то заменой $\bar t=t$, $\bar x=x+\frac {R}{Q} t$ мы приводим нашу функцию к виду $\mathcal{F}=Pc^3f$, также являющуюся представителем семейства II.2. Если $P=Q=0$, но $R\ne 0$, то функция сразу оказывается частным случаем II.1.

\subsubsection{$\mathcal{F}=G(t)f^{-1/2}$}
Здесь подстановка и расщепление по $f$ и по $c$ дает
$$\tau_x=0,\qquad \tau G'(t)+\frac 32\tau_t G(t) =0, \qquad \xi_{tt}=2\xi_{tx} - \tau_{tt}=\xi_{xx}=0,$$
откуда $\xi=atx+bx+kt+l$, $\tau=at^2+pt+q$, и либо $\tau=0$, и мы получаем трехмерную алгебру III.2.4, либо
$G(t)=(at^2+pt+q)^{-3/2}$, $\mathcal{F}=(at^2+pt+q)^{-3/2}f^{-1/2}$, что сдвигами по $t$ приводится либо к $\mathcal{F}=P(t^2\pm k^2)^{-3/2}f^{-1/2}$, либо к $\mathcal{F}=Pt^{-3}f^{-1/2}$, либо к $\mathcal{F}=Pt^{-3/2}f^{-1/2}$, либо к $\mathcal{F}=Pf^{-1/2}$. Первая третья из перечисленных функций относятся к семейству II.4, а вторая и четвертая -- к семейству II.1.

\subsubsection{$\mathcal{F}=(c^2+n)^{3/2}G((c^2+n)^{3/2}f)$}
Случай $n=0$ дает нам функцию из семейства III.2.1, уже рассмотренного ранее, поэтому далее мы считаем $n\ne 0$.
Подстановка в уравнение и обозначение аргумента функции $G$ через $z$ приводит к равенству
$$ (\xi_t + c(\xi_x - \tau_t) - c^2\tau_x)3c(c^2+n)^{1/2}(G(z)+zG'(z)) + (3c\tau_x - 2\xi_x + \tau_t)(c^2+n)^{3/2}zG'(z) +$$
$$+(3c\tau_x -\xi_x + 2\tau_t)(c^2+n)^{3/2}G(z)= \xi_{tt} + c(2\xi_{tx} - \tau_{tt}) + c^2(\xi_{xx} - 2\tau_{tx}) - c^3\tau_{xx}.$$
Поскольку слева стоит иррациональная функция от $c$, а справа -- многочлен, правая и левая часть равенства равны нулю по отдельности:
$$ (\xi_t + c(\xi_x - \tau_t) - c^2\tau_x)3c(G(z)+zG'(z)) + (3c\tau_x - 2\xi_x + \tau_t)(c^2+n)zG'(z) +(3c\tau_x -\xi_x + 2\tau_t)(c^2+n)G(z)= 0,$$
$$\xi_{tt}=2\xi_{tx} - \tau_{tt}=\xi_{xx} - 2\tau_{tx}=\tau_{xx}=0.$$
Если $G(z)+zG'(z)= 0$, то мы получаем функцию I. Если же $G(z)+zG'(z)\ne 0$, то в первом уравнении второе и третье слагаемые содержат множитель $c^2+n$, поэтому первый множитель в первом слагаемом должен делиться на $c^2+n$, что дает $\xi_x=\tau_t$, $\xi_t+n\tau_x=0$, откуда следует равенство нулю всех вторых производных $\xi$ и $\tau$, первое же равенство при этом приобретает вид
$\tau_t[-zG'(z)+G(z)]=0$. Если $G(z)\ne Pz$, то $\tau_t=\xi_x=0$, и мы получаем трехмерную алгебру, указанную в III.3.1, в противном случае мы приходим к функции II.2.

\subsubsection{$\mathcal{F} = c^mG(c^{3-m}f)$}

Подставив $\mathcal{F} = c^mG(c^{3-m}f)$ в (27), обозначив через $z$ аргумент функции $G$, получим после группировки
$$ [(m\xi_t + c((m-1)\xi_x - (m-2)\tau_t) - (m-3)c^2\tau_x]c^{m-1}G(z)+$$
$$+[((3-m)\xi_t + c((1-m)\xi_x - (2-m)\tau_t) +mc^2\tau_x)]c^{m-1}zG'(z)=$$
$$= \xi_{tt} + c(2\xi_{tx} - \tau_{tt}) + c^2(\xi_{xx} - 2\tau_{tx}) - c^3\tau_{xx}.\eqno(30)$$

Если $m\ne 0,1,2,3$ (эти случаи мы рассмотрим отдельно), то либо функции $zG'(z)$, $G(z)$ и единица линейно независимы -- тогда $\tau_x=\xi_t=0$, $(m-1)\xi_x=(m-2)\tau_t$, и мы получаем трехмерную алгебру III.3.2. Либо они линейно зависимы, и тогда $zG'(z)=\alpha G(z)+\beta$, где $\alpha\ne 0$ (если коэффициент при $G(z)$ нулевой, то нулевым оказывается и коэффициент при $zG'(z)$), что дает $G(z)=Pz^j+Q$.

Для этой $G(z)$ наше уравнение приводится к виду
$$ [(m\xi_t + c((m-1)\xi_x - (m-2)\tau_t) - (m-3)c^2\tau_x]c^{m-1}(Pz^j+Q)+$$
$$+[((3-m)\xi_t + c((1-m)\xi_x - (2-m)\tau_t) +mc^2\tau_x)]c^{m-1}jPz^j=\xi_{tt} + c(2\xi_{tx} - \tau_{tt}) + c^2(\xi_{xx} - 2\tau_{tx}) - c^3\tau_{xx},$$
приравнивание коэффициентов при $z^j$ и при членах, не содержащих $z$, дает
$$ [m+j(3-m)]\xi_t + c(1-j)[(m-1)\xi_x - (m-2)\tau_t] - [(m-3)-jm]c^2\tau_x=0,$$
$$ Q[(m\xi_t + c((m-1)\xi_x - (m-2)\tau_t) - (m-3)c^2\tau_x]c^{m-1}=\xi_{tt} + c(2\xi_{tx} - \tau_{tt}) + c^2(\xi_{xx} - 2\tau_{tx}) - c^3\tau_{xx},$$
и из первого уравнения при $j\ne 1,\frac{m-3}{m}, \frac{m}{m-3}$ снова получаем ту же трехмерную алгебру. Если $j=1$, то из первого уравнения $\xi_t=\tau_x=0$, второе уравнение приобретает вид
$$ Q c^m((m-1)\xi_x - (m-2)\tau_t)= -c\tau_{tt} + c^2\xi_{xx},$$
и поскольку $m\ne 1,2$, мы снова получаем ту же трехмерную алгебру. Исключение составляет случай $Q=0$, что дает нам функцию из семейства II.2.

При $j=\frac m{m-3}$ мы из первого уравнения получаем
$$ -c\frac 3{m-3}[(m-1)\xi_x - (m-2)\tau_t] - \frac{9-6m}{m-3}c^2\tau_x=0,$$
откуда при $m\ne \frac 32$ немедленно следует $(m-1)\xi_x=(m-2)\tau_t$, $\tau_x=0$,
$$ Q m\xi_t c^{m-1}=\xi_{tt} + c(2\xi_{tx} - \tau_{tt}) + c^2\xi_{xx},$$
и в силу $m\ne 1,2,3$ получаем при $Q\ne 0$, что $\xi_t=0$, приходя ко все той же трехмерной алгебре. Если же $Q=0$, то мы получаем функцию II.1.

Аналогично при $j=\frac{m-3}m$ мы из первого уравнения при $m\ne \frac 32$ получаем $\xi_t=0$, $(m-1)\xi_x=(m-2)\tau_t$, второе уравнение приобретает вид
$$ -Q(m-3)\tau_x c^{m+1}= -c\tau_{tt} + c^2(\xi_{xx} - 2\tau_{tx}) - c^3\tau_{xx},$$
и опять же, поскольку $m\ne 0,1,2$, мы получаем при $Q\ne 0$, что $\tau_x=0$, приходя ко все той же трехмерной алгебре. Если же $Q=0$, мы получаем функцию $\mathcal{F}=Pc^3(fc^3)^j$, которая заменой $\bar t=x$, $\bar x=t$ приводится к той же функции II.1.

Наконец, случай $m=3/2$, $j=-1$ дает нам уравнения $\xi_x+\tau_t=0$,
$$ Q[(\frac 32 \xi_t + c\frac 12(\xi_x +\tau_t) +\frac 32 c^2\tau_x]c^{1/2}=\xi_{tt} + c(2\xi_{tx} - \tau_{tt}) + c^2(\xi_{xx} - 2\tau_{tx}) - c^3\tau_{xx},$$
откуда получаем при $Q\ne 0$ $\xi_t=\tau_x=0$ и трехмерную алгебру, а при $Q=0$ -- функцию I.

У нас остались не рассмотренными случаи $m=0,1,2,3$. Поскольку замена $\bar t=x$, $\bar x=t$ переводит случай $m=0$ в случай $m=3$, а случай $m=1$ в случай $m=2$, достаточно рассмотреть только случаи $m=1$ и $m=3$.

Если $m=3$, то уравнение (30) приобретает вид
$$ [3\xi_t + c(2\xi_x - \tau_t)]c^2 G(z)
+[c(-2\xi_x+\tau_t)+3c^2\tau_x]c^2 zG'(z)=\xi_{tt} + c(2\xi_{tx} - \tau_{tt}) + c^2(\xi_{xx} - 2\tau_{tx}) - c^3\tau_{xx},$$
и ввиду того, что $G(z)$ в нашей теореме предполагается непостоянной, мы немедленно получаем, что равенство нулю коэффициента при $c^4$ означает $\tau_x=0$, а тогда
$$(2\xi_x - \tau_t)c^3 (G(z)-zG'(z))=0, \qquad \xi_{tt}=0,\qquad 2\xi_{tx} - \tau_{tt}=0,\qquad 3\xi_t G(z)= \xi_{xx}.$$
Опять же, поскольку $G(z)$ непостоянная, из последнего соотношения мы получаем $\xi_t=\xi_{xx}=0$, а из предпоследнего -- $\tau_{tt}=0$. И у нас остается либо четырехмерная алгебра, если $G(z)=zG'(z)$, то есть $G(z)=Pz$ -- это функция из семейства  II.2, либо $2\xi_x=\tau_t$, и мы возвращаемся ко все той же трехмерной алгебре.

В случае $m=1$ уравнение (30) приобретает вид
$$ (\xi_t + c\tau_t +2c^2\tau_x)G(z)+(2\xi_t - c\tau_t +c^2\tau_x)zG'(z)=\xi_{tt} + c(2\xi_{tx} - \tau_{tt}) + c^2(\xi_{xx} - 2\tau_{tx}) - c^3\tau_{xx},$$
его расщепление по степеням $c$ дает
$$\tau_{xx}=0, \ \tau_x [2G(z)+zG'(z)]=\xi_{xx} - 2\tau_{tx}, \ \tau_t [G(z)-zG'(z)]=2\xi_{tx} - \tau_{tt}, \
\xi_t [G(z)+2zG'(z)]= \xi_{tt},$$
откуда, если $G(z)\ne Pz^{-2}+Q, Pz+Q, Pz^{-1/2}+Q$, мы получаем $\xi_t=\xi_x=\tau_x=\tau_{tt}=0$, что опять дает нам трехмерную алгебру. Если $G(z)= Pz+Q$, то мы получаем функцию семейства II.2, если $G(z)= Pz^{-1/2}+Q$, то функцию из семейства II.4, а если $G(z)= Pz^{-2}+Q$, то мы функцию $\mathcal{F}=\frac P{c^3f^2}+Qc$, которая заменой $\bar t=x$, $\bar x=t$ превращается $\mathcal{F}=\frac P{f^2}+Qc^2$, то есть в функцию из семейства II.3.

\subsubsection{$\mathcal{F} = e^{mc}G(e^{-mc}f)$}
Поскольку при $m=0$ мы попадаем в уже рассмотренное семейство III.2.1, нас будет интересовать случай $m\ne 0$.
Подстановка в (27) и обозначение аргумента функции $G$ через $z$ приводит к равенству
$$e^{mc}[(m\xi_t -\xi_x + 2\tau_t + c(m\xi_x - m\tau_t+3\tau_x) - c^2m\tau_x)G(z)-$$
$$-(m\xi_t+2\xi_x - \tau_t + c(m\xi_x - m\tau_t-3\tau_x) - c^2m\tau_x)zG'(z)]=$$
$$= \xi_{tt} + c(2\xi_{tx} - \tau_{tt}) + c^2(\xi_{xx} - 2\tau_{tx}) - c^3\tau_{xx}.$$
Иррациональность экспоненты влечет равенство нулю по отдельности левой и правой частей этого равенства, расщепление их по $c$ дает систему
$$m\tau_x(-G(z)  + zG'(z))=0,\quad (m\xi_x - m\tau_t+3\tau_x)G(z)- (m\xi_x - m\tau_t-3\tau_x)zG'(z)=0,$$
$$(m\xi_t -\xi_x + 2\tau_t)G(z) -(m\xi_t+2\xi_x - \tau_t)zG'(z)=0,$$
$$\xi_{tt}=2\xi_{tx} - \tau_{tt}=\xi_{xx} - 2\tau_{tx}=\tau_{xx}=0.$$
Если $G(z)$ и $zG'(z)$ линейно независимы, то мы получаем
$\tau_x=\xi_x-\tau_t=m\xi_t+\tau_t=0$, и приходим к  трехмерной алгебре, указанной в III.3.3.

В противном случае $G(z)=Pz^j$, и система приобретает вид
$$\tau_x(j-1)=0,\quad m(1-j)\xi_x - m(1-j)\tau_t+3(1+j)\tau_x=0,\quad m(1-j)\xi_t -(2j+1)\xi_x + (2+j)\tau_t=0,$$
$$\xi_{tt}=2\xi_{tx} - \tau_{tt}=\xi_{xx} - 2\tau_{tx}=\tau_{xx}=0.$$
При $j=1$ мы получаем линейную функцию из семейства II.1, а при $j\ne 1$ из первого уравнения следует $\tau_x=0$, из второго -- $\xi_x=\tau_t$, а из третьего -- $m\xi_t+\tau_t=0$, и мы снова приходим к трехмерной алгебре III.3.3.

\subsubsection{$\mathcal{F}=(t^2\pm k^2)^{-3/2}e^{\lambda h(t)}G(fe^{2\lambda h(t)})$, где  $h(t)=\int\frac{dt}{t^2\pm k^2}$}
Подставим в (27), и обозначим аргумент $G$ через $z$, получим уравнение
$$[\tau((\lambda-3t)G(z)+2\lambda zG'(z)) + (3c\tau_x - 2\xi_x + \tau_t)(t^2\pm k^2)zG'(z) +\eqno(31)$$
$$+(3c\tau_x -\xi_x + 2\tau_t)(t^2\pm k^2)G(z)](t^2\pm k^2)^{-5/2}e^{\lambda h(t)}=
\xi_{tt} + c(2\xi_{tx} - \tau_{tt}) + c^2(\xi_{xx} - 2\tau_{tx}) - c^3\tau_{xx}.$$
Расщепление по $c$ приводит к системе
$$\xi_{tt}\!\!=\!\![\tau((\lambda-3t)G(z)+2\lambda zG'(z)) + (\tau_t- 2\xi_x)(t^2\pm k^2)zG'(z) +(2\tau_t-\xi_x)(t^2\pm k^2)G(z)](t^2\pm k^2)^{-5/2}e^{\lambda h(t)},$$
$$2\xi_{tx} - \tau_{tt}=3\tau_x [zG'(z) +G(z)](t^2\pm k^2)^{-3/2}e^{\lambda h(t)},\qquad \xi_{xx} - 2\tau_{tx}=\tau_{xx}=0.$$
Если $zG'(z)$, $G(z)$ и единица линейно независимы, то из этой системы получаем $\tau_x=0$, $\xi_{xx}=0$, $\xi_{tt}=0$, $2\xi_{tx} - \tau_{tt}=0$, откуда $\xi=ax+btx+pt+q$, $\tau=bt^2+rt+s$ и
$$2\tau\lambda+(\tau_t- 2\xi_x)(t^2\pm k^2)=0,\qquad
\tau (\lambda-3t)+(2\tau_t-\xi_x)(t^2\pm k^2)=0.$$
Подставляя в эти уравнения $\tau$ и $\xi$, и расщепляя их по $t$, получаем
$$2\lambda b+(r-2a)=0, \qquad \lambda r=0,\qquad 2\lambda s\pm k^2(r-2a)=0,$$
$$\lambda b-3r+(2r-a)=0,\qquad \lambda r-3s \pm 3k^2 b=0,\qquad \lambda s\pm k^2(2r-a)=0.$$
Из первого и четвертого следует $r=0$ и $a=\lambda b$, из пятого $s=\pm k^2b$, второе, третье и шестое тогда выполняются автоматически. Таким образом,
$\tau=b(t^2\pm k^2)$, $\xi=b(t+\lambda)x+pt+q$, и мы получаем трехмерную алгебру III.4.1.

Если же $zG'(z)=\alpha G(z)+\beta$, то мы, после выполнения этой подстановки в (31), получаем другое расщепление:
$$\tau(\lambda-3t+2\lambda\alpha ) + [(\alpha+2)\tau_t- (2\alpha+1)\xi_x](t^2\pm k^2) =0,$$
$$\xi_{tt}=[2\tau\lambda \beta + \beta(\tau_t- 2\xi_x)(t^2\pm k^2)](t^2\pm k^2)^{-5/2}e^{\lambda h(t)},$$
$$3\tau_x (\alpha+1)=0,\qquad 2\xi_{tx} - \tau_{tt}=3\tau_x \beta (t^2\pm k^2)^{-3/2}e^{\lambda h(t)},\qquad \xi_{xx} - 2\tau_{tx}=\tau_{xx}=0,$$
и при $\alpha\ne -1$ мы получаем из третьего, четвертого и пятого уравнений $\tau_x=0$, $\xi_{xx}=0$, $2\xi_{tx} - \tau_{tt}=0$, откуда $\tau=\tau(t)$, $\xi=\frac 12\tau'(t)x+ax+b(t)$. Подставляя эти функции в первое и второе уравнения и расщепляя второе по $x$, получаем
$$\tau(\lambda-3t+2\lambda\alpha ) + [\frac 32\tau'(t)-(2\alpha+1)a](t^2\pm k^2) =0,$$
$$\tau'''(t)=0,\qquad b''(t)=[2\tau\lambda \beta -2a \beta(t^2\pm k^2)](t^2\pm k^2)^{-5/2}e^{\lambda h(t)}.$$
Из второго соотношения $\tau=pt^2+qt+r$. Подставляя это в первое соотношение, получаем
$$\lambda(1+2\alpha)p-\frac 32 q-(2\alpha+1)a=0,\qquad
\lambda(1+2\alpha)q-3r\pm 3p k^2 =0,\qquad \lambda(1+2\alpha)r\pm k^2[\frac 32q-(2\alpha+1)a]=0.$$
Выразим из первого равенства $q$, из второго $r$ и подставим в третье:
$$q=\frac 23(1+2\alpha)(\lambda p-a),\qquad
r=\frac 13\lambda(1+2\alpha)q\pm p k^2 ,\qquad 2\left[\frac 19\lambda^2(1+2\alpha)^2 \pm k^2\right](1+2\alpha)(\lambda p-a)=0.$$
Из последнего равенства следует, что возможны три варианта.

Либо $a=\lambda p$, тогда $q=0$, $r=\pm pk^2$, $\tau=p(t^2\pm k^2)$, $b''(t)=0$, и мы получаем все ту же трехмерную алгебру III.4.1.

Либо $\alpha=-\frac 12$, тогда $\mathcal{F}=(t^2\pm k^2)^{-3/2}(Pf^{-1/2}+Q)$, и подходящей заменой $\bar t=t$, $\bar x=x-\phi(t)$ эта функция приводится к II.4.

Наконец, в последнем случае $\frac 19\lambda^2(1+2\alpha)^2 \pm k^2=0$ (возможном только если $1+2\alpha\ne 0$ и $\pm k^2=-k^2$, ) мы для $\lambda=\pm\frac{3k}{1+2\alpha}$ получаем при $\alpha\ne 0$
$$h(t)=\frac 1{2k}\ln\frac {t-k}{t+k}, \qquad  e^{\lambda h(t)}=\left(\frac {t-k}{t+k}\right)^{\pm \frac 3{2(1+2\alpha)}},$$
$$\mathcal{F}=(t^2-k^2)^{-\frac 32}e^{\lambda (2\alpha+1)h(t)}Pf^\alpha+Q(t^2-k^2)^{-\frac 32}e^{\lambda h(t)}=$$
$$=(t^2-k^2)^{-\frac 32}\left(\frac {t-k}{t+k}\right)^{\pm 3/2}Pf^\alpha+Q(t^2-k^2)^{-\frac 32}\left(\frac {t-k}{t+k}\right)^{\pm \frac 3{2(1+2\alpha)}} = $$
$$=(t\pm k)^{-3}Pf^\alpha +Q \frac{(t\mp k)^{-3\alpha/(2\alpha+1)}}{(t\pm k)^{3(\alpha+1)/(2\alpha+1)}},$$
эта функция сдвигом переменной $t$ и подходящей заменой $\bar t=t$, $\bar x=x-\phi(t)$ приводится к II.1. Если же $\alpha=0$, то мы получаем уже функцию
$$\mathcal{F}=\beta (t\pm k)^{-3} \ln |f| \pm (t\pm k)^{-3}\frac 32\ln\frac {t-k}{t+k},$$
которая теми же преобразованиями приводится к II.5.

\subsubsection{$\mathcal{F}=t^{k-2}G(t^{2k-1}f)$}
Подстановка этой функции в (27) приводит, после расщепления по $c$, к системе (как и ранее, аргумент функции $G$ обозначен через $z$):
$$\tau_{xx}=\xi_{xx} - 2\tau_{tx}=0,\qquad
2\xi_{tx} - \tau_{tt}=3\tau_x t^{k-2}[zG'(z) +G(z)],$$
$$\xi_{tt}=\tau[(k-2)t^{k-3}G(z)+(2k-1)t^{k-3}zG'(z)]+(- 2\xi_x + \tau_t)t^{k-2}zG'(z) +(-\xi_x + 2\tau_t)t^{k-2}G(z).$$
Все так же, предполагая, что $zG'(z)$, $G(z)$ и единица линейно независимы, получаем, что
$$\tau_x=0,\qquad \xi_{xx}=0,\qquad \xi_{tt}=0\qquad, 2\xi_{tx} - \tau_{tt}=0,$$
$$(k-2)\frac \tau t+(-\xi_x + 2\tau_t)=0,\qquad (2k-1)\frac \tau t+(- 2\xi_x + \tau_t)=0.$$
Из последних двух равенств следует $\tau_t=\frac\tau t$, $\xi_x=k\frac\tau t$, а значит $\tau=at$, $\xi=ax+pt+q$ -- мы получили трехмерную алгебру III.4.2.

Если же $zG'(z)=\alpha G(z)+\beta$, то мы получаем
$$\mathcal{F}=t^{(2\alpha+1)k-2-\alpha}f^\alpha-\frac \beta\alpha t^{k-2}\qquad \mbox{ при }\alpha\ne 0,$$
$$\mathcal{F}=t^{k-2}\beta \ln |f|+(2k-1)t^{k-2}\ln t\qquad \mbox{ при }\alpha=0.$$
И та, и другая функция подходящей заменой вида $\bar t=t$, $\bar x=x-\phi(t)$ сводится к
$$\mathcal{F}=t^{(2\alpha+1)k-2-\alpha}Pf^\alpha \mbox{ \ и \ }\mathcal{F}=t^{k-2}\beta \ln |f|$$
соответственно.

Для первой из них уравнение (27) приобретает вид (через $m$ обозначено $k(2\alpha+1)-\alpha-2$ -- показатель степени у $t$)
$$\tau mt^{m-1}f^\alpha  + \alpha (3c\tau_x - 2\xi_x + \tau_t)t^mf^\alpha +(3c\tau_x -\xi_x + 2\tau_t)t^mf^\alpha= \xi_{tt} + c(2\xi_{tx} - \tau_{tt}) + c^2(\xi_{xx} - 2\tau_{tx}) - c^3\tau_{xx},$$
и его расщепление по $f$ и по $c$ дает систему
$$(\alpha+1)\tau_x=0,\quad m\frac \tau t - (2\alpha+1)\xi_x + (\alpha+2)\tau_t=0,\quad \xi_{tt}=2\xi_{tx} - \tau_{tt}=\xi_{xx} - 2\tau_{tx}=\tau_{xx}=0.$$
Из этих уравнений, кроме второго, следует при $\alpha\ne -1$, что $\tau=pt^2+qt+r$, $\xi=ptx+ax+bt+l$, подстановка этих функций во второе уравнение и расщепление его по $t$ дают три равенства:
$$mr=p(m+3)=(kq - a)(2\alpha+1)=0,$$
из которых следует, что при $\alpha\ne -\frac 12$ и $m\ne 0,-3$ мы получаем $p=r=0$, $\tau=qt$, $\xi=kqx+bt+l$, то есть ту же самую трехмерную алгебру III.4.2. Случаи $m=0$ и $m=-3$ приводят нас к функциям из семейства II.1, а случай $\alpha=-\frac 12$ -- к функции из семейства II.4.

Случай $\alpha=-1$ чуть сложнее: из
$\xi_{tt}=2\xi_{tx} - \tau_{tt}=\xi_{xx}-2\tau_{tx}=\tau_{xx}=0$
мы получаем, что $\tau=At^2+Btx+Ct+Dx+E$, $\xi=Atx+Bx^2+Ft+Gx+H$, и подстановка этих функций в оставшееся уравнение
$m\frac \tau t +\xi_x +\tau_t=0$ дает
$mD=mE=0$, $(m+3)A=(m+3)B=0$, $G=-C(m+1)$, откуда при $m\ne 0,-3$ мы снова получаем всю ту же трехмерную алгебру, а при $m=0, -3$ -- функцию из семейства I.

Что же касается второй -- логарифмической -- функции, то для нее уравнение (27) приобретает вид
$$(k-2)\tau t^{k-3}\beta \ln |f| + (3c\tau_x - 2\xi_x + \tau_t)t^{k-2}\beta +(3c\tau_x -\xi_x + 2\tau_t)t^{k-2}\beta \ln |f| =\hspace{50mm} $$
$$\hspace{60mm}=\xi_{tt} + c(2\xi_{tx} - \tau_{tt}) + c^2(\xi_{xx} - 2\tau_{tx}) - c^3\tau_{xx},$$
и расщепление по $f$ и $c$ дает систему
$$\tau_x=0, \qquad (k-2)\frac \tau t -\xi_x + 2\tau_t=0,\qquad  \xi_{tt}=(- 2\xi_x + \tau_t)t^{k-2}\beta, \qquad 2\xi_{tx} - \tau_{tt}= \xi_{xx}=0,$$
из которой мы по той же схеме получаем сначала $\tau=pt^2+qt+r$, $\xi=ptx+ax+b(t)$, затем подставновкой во второе приходим к
$(k-2)r=0$ $p(k+1)=0$, $a=kq$, откуда при $k\ne 2,-1$ мы получаем в очередной раз трехмерную алгебру, а случаи $k=2,-1$ приводят нас к семейству II.5.

\subsubsection{$\mathcal{F}=e^{kt}G(e^{2kt}f)$}
Здесь подстановкой в (27) получаем уравнение
$$e^{kt}[(3c\tau_x - 2\xi_x + \tau_t+2k\tau)zG'(z) +(3c\tau_x -\xi_x + 2\tau_t+k\tau)G(z)]
= \xi_{tt} + c(2\xi_{tx} - \tau_{tt}) + c^2(\xi_{xx} - 2\tau_{tx}) - c^3\tau_{xx},$$
которое в случае линейной независимости функций $zG'(z)$, $G(z)$ и единицы приводит к системе
$$\tau_x=0,\qquad  - 2\xi_x + \tau_t+2k\tau=0,\qquad -\xi_x + 2\tau_t+k\tau=0,\qquad
\xi_{tt}=2\xi_{tx} - \tau_{tt}=\xi_{xx} - 2\tau_{tx}=\tau_{xx}=0,$$
из которой $\tau=a={\rm const}$, $\xi=kax+bt+l$, и мы получаем трехмерную алгебру III.4.3.

Если же $zG'(z)=\alpha G(z)+\beta$, то, как и в предыдущем случае, мы подходящей заменой приводим соответствующие $\mathcal{F}$ к виду
$$\mathcal{F}=e^{k(2j+1)t}f^j\qquad \mbox{или}\qquad \mathcal{F}=e^{kt}\ln|f|.$$
Для первой функции подстановка в (27) и расщепление по $f$ и по $c$ дает систему
$$\tau_x=0,\qquad k(2j+1)\tau - (2j+1)\xi_x + (j+2)\tau_t=0,\qquad
\xi_{tt}=2\xi_{tx} - \tau_{tt}=\xi_{xx}=0,$$
из которой $\tau=pt^2+qt+r$, $\xi=ptx+ax+bt+l$, и подстановка этих функций во второе уравнение дает
$k(2j+1)(pt^2+qt+r) - (2j+1)(pt+a) + (j+2)(2pt+q)=0$, откуда при $k\ne 0$, $j\ne -\frac12$ получаем $p=q=0$, $a=kr$, то есть ту же трехмерную алгебру III.4.3. Если $k=0$ или $j=-\frac 12$, то мы получаем функцию семейства II.1.

Для второй функции те же рассуждения дают сначала систему
$$\tau_x=0,\qquad  -\xi_x + 2\tau_t+k\tau=0,\qquad
\xi_{tt}=(- 2\xi_x + \tau_t)e^{kt},\qquad
2\xi_{tx} - \tau_{tt}=\xi_{xx}=0.$$
Из первого, третьего и последнего уравнений мы получаем
$\tau=pt^2+qt+r$, $\xi=ptx+ax+b(t)$, подстановка во второе уравнение приводит к условию
$-(pt+a) + 2(2pt+1)+k(pt^2+qt+r)=0$, из которого при $k\ne 0$ опять $p=q=0$, $a=kr$, $b''(t)=-2ae^{kt}$, и мы снова получаем трехмерную алгебру, а при $k=0$ -- функцию из семейства II.5.

\subsubsection{$\mathcal{F}=\left(1\pm \frac{(c\mp t)^2}{2x\mp t^2}\right)^{3/2}G\left((2x\mp t^2\pm(c\mp t)^2)^{3/2}f\right)$}
Здесь подстановка функции в (27) приводит к уравнению
$$\left(1\pm \frac{(c\mp t)^2}{2x\mp t^2}\right)^{1/2}\left\{\tau\left [3\frac{(c\mp t)(ct-2x)}{(2x\mp t^2)^2}G(z)-\frac{3c}{2x\mp t^2\pm(c\mp t)^2}
\left(1\pm \frac{(c\mp t)^2}{2x\mp t^2}\right)zG'(z)\right]\right. +$$
$$+ \xi\left [\mp 3\frac{(c\mp t)^2}{(2x\mp t^2)^2}G(z)
+\frac 3{2x\mp t^2\pm(c\mp t)^2}\left(1\pm \frac{(c\mp t)^2}{2x\mp t^2}\right)zG'(z)\right] +$$
$$+(\xi_t + c(\xi_x - \tau_t) - c^2\tau_x)\left[\pm 3\frac{c\mp t}{2x\mp t^2}G(z)
\pm \frac{3(c\mp t)}{2x\mp t^2\pm(c\mp t)^2}\left(1\pm \frac{(c\mp t)^2}{2x\mp t^2}\right)zG'(z)\right] +$$
$$\left .+(3c\tau_x - 2\xi_x + \tau_t)\left(1\pm \frac{(c\mp t)^2}{2x\mp t^2}\right)zG'(z) +(3c\tau_x -\xi_x + 2\tau_t)\left(1\pm \frac{(c\mp t)^2}{2x\mp t^2}\right)G(z)\right\} =$$
$$= \xi_{tt} + c(2\xi_{tx} - \tau_{tt}) + c^2(\xi_{xx} - 2\tau_{tx}) - c^3\tau_{xx}.$$
Поскольку иррациональный (по $c$) множитель не может быть представлен в виде отношения рациональных функций, левая и правая части этого равенства равны нулю по отдельности, что дает
$$\xi_{tt}=2\xi_{tx} - \tau_{tt}=\xi_{xx} - 2\tau_{tx}=\tau_{xx}=0,$$

$$\left[\frac{3\tau(c\mp t)(ct-2x)\mp3\xi(c\mp t)^2}{(2x\mp t^2)^2}
\pm 3(\xi_t + c(\xi_x - \tau_t) - c^2\tau_x)\frac{c\mp t}{2x\mp t^2}+\right.\hspace{35mm}$$
$$\hspace{65mm}\left.+(3c\tau_x -\xi_x + 2\tau_t)\left(1\pm \frac{(c\mp t)^2}{2x\mp t^2}\right)\right]G(z)+$$
$$+\left [\frac{-3c\tau+3\xi
\pm 3(\xi_t + c(\xi_x - \tau_t) - c^2\tau_x)(c\mp t)}{2x\mp t^2}+(3c\tau_x - 2\xi_x + \tau_t)\left(1\pm \frac{(c\mp t)^2}{2x\mp t^2}\right)\right]zG'(z) =0.$$
Во втором равенстве слева стоит многочлен от $c$, который равен нулю. Для того, чтобы получить из него нужную нам систему уравнений, удобно ввести обозначение $v=\frac{c\mp t}{2x\mp t^2}$, тогда второе равенство примет вид
$$\left[3\tau(v^2t-v)\mp 3\xi v^2
\pm 3(\xi_t + (v(2x\mp t^2)\pm t)(\xi_x - \tau_t) - (v(2x\mp t^2)\pm t)^2\tau_x)v+\right.$$
$$\left.+(3(v(2x\mp t^2)\pm t)\tau_x -\xi_x + 2\tau_t)\left(1\pm v^2(2x\mp t^2)\right)\right]G(z)+$$
$$+\left [-3v\tau- \frac{\pm 3t\tau- 3\xi}{2x\mp t^2}
\pm 3(\xi_t + (v(2x\mp t^2)\pm t)(\xi_x - \tau_t) - (v(2x\mp t^2)\pm t)^2\tau_x)v+\right.$$
$$\left.+(3(v(2x\mp t^2)\pm t)\tau_x - 2\xi_x + \tau_t)\left(1\pm v^2(2x\mp t^2)\right)\right]zG'(z) =0,$$
и уже легко расщепляется по $v$:
$$\left[\pm 3\tau t-3\xi +(2x\mp t^2)(2\xi_x - \tau_t \mp 3 t\tau_x)\right]G(z)+(2x\mp t^2)\left [\xi_x - 2\tau_t \mp 3t\tau_x\right]zG'(z) =0,$$
$$\left[-3\tau
\pm 3(\xi_t \pm t(\xi_x - \tau_t) - t^2\tau_x)+3(2x\mp t^2)\tau_x \right](G(z)+zG'(z)) =0,$$
$$\left[\mp 3t\tau_x +\xi_x - 2\tau_t\right]G(z)+\left [
\frac{\pm 3t\tau-3\xi}{2x\mp t^2}\mp 3t\tau_x + 2\xi_x - \tau_t\right]zG'(z) =0.$$
Первое и третье уравнения образуют линейную однородную систему относительно $G(z), zG'(z)$ с определителем
$$\left[\pm 3\tau t-3\xi + (2x\mp t^2)(2\xi_x - \tau_t \mp 3 t\tau_x)\right]^2
-(2x\mp t^2)^2\left [\xi_x - 2\tau_t \mp 3t\tau_x\right]^2.$$
Поскольку $G(z)$ предполагается ненулевой, определитель равен нулю, а значит, либо
$$\pm 3\tau t-3\xi + (2x\mp t^2)(2\xi_x - \tau_t \mp 3 t\tau_x)=
(2x\mp t^2)\left [\xi_x - 2\tau_t \mp 3t\tau_x\right],$$
либо
$$\pm 3\tau t-3\xi + (2x\mp t^2)(2\xi_x - \tau_t \mp 3 t\tau_x)=
-(2x\mp t^2)\left [\xi_x - 2\tau_t \mp 3t\tau_x\right].$$
Подставляя в эти уравнения следующие из
$\xi_{tt}=2\xi_{tx} - \tau_{tt}=\xi_{xx} - 2\tau_{tx}=\tau_{xx}=0$ формулы для $\tau$, $\xi$
$$\tau=At^2+Btx+Ct+Dx+E,\qquad \xi=Atx+Bx^2+Ft+Gx+H,$$
получаем в первом случае
$G=2C$, $D=\mp A$, $F=\pm E$, $B=H=0$,
а во втором $G=2C$, $F=\pm E$, $A=D=B=H=0$, что в итоге дает трехмерную алгебру $\tau=A(t^2\mp x)+Ct+E$, $\xi=Atx\pm Et+2Cx$, которая обращает все остальные уравнения в тождества и которая и приведена в III.5.


\subsection{Алгебры симметрий уравнений с функцией $\mathcal{F}=\mathcal{F}(c,f)$}
В этом параграфе мы покажем, что для всех функций $\mathcal{F}=\mathcal{F}(c,f)$, за исключением функций из семейств I-III, алгебра симметрий уравнения (5) в точности двумерна.
\subsubsection{Анзацы и классифицирующие уравнения} Подставим $\mathcal{F}(c,f)$ в уравнение $(27)$:
$$(\xi_t + c(\xi_x - \tau_t) - c^2\tau_x)\mathcal{F}_c + f(3c\tau_x - 2\xi_x + \tau_t)\mathcal{F}_f + (3c\tau_x -\xi_x + 2\tau_t)\mathcal{F} =$$
$$= \xi_{tt} + c(2\xi_{tx} - \tau_{tt}) + c^2(\xi_{xx} - 2\tau_{tx}) - c^3\tau_{xx}. \eqno(32)$$

Поскольку для $\mathcal{F}(c,f)$ уравнение (32) выполняется для любых постоянных $\tau$, $\xi$, алгебра симметрий содержит двумерную коммутативную алгебру, и для завершения рассуждений с алгебрами, содержащими двумерную коммутативную подалгебру, нам необходимо показать, что в случаях, отличных от рассмотренных в предыдущих параграфах, алгебра симметрий имеет размерность, в точности равную двум.

При этом мы будем пользоваться результатами леммы 1, в соответствии с которой для рассматриваемого семейства функций группа эквивалентности -- это группа линейных преобразований переменных $(t,x)$. Поэтому функции, переводящиеся друг в друга такими преобразованиями, принадлежат одному классу эквивалентности, из которого мы будем обычно указывать только одного представителя.

Итак, предположим, что уравнение (32) выполнено для некоторых непостоянных $\tau(t,x)$ , $\xi(t,x)$. Подставим их в уравнение, которое теперь будем рассматривать как уравнение относительно $\mathcal{F}(c,f)$.

1) Предположим вначале, что $\xi_t + c(\xi_x - \tau_t) - c^2\tau_x\equiv 0$. Тогда $\tau_x=\xi_t=0$, а $\tau_t=\xi_x={\rm const}$. Эта константа -- ненулевая в силу предположения о том, что хотя бы одна из функций $\tau$, $\xi$ не является постоянной. Но тогда оказывается, что $2c\tau_x -\xi_x + \tau_t = 0$, а выражение $3c\tau_x - 2\xi_x + \tau_t$ равно минус этой константе, в итоге наше уравнение превращается в $-f\mathcal{F}_f+\mathcal{F}=0$, а значит, $\mathcal{F}$ -- линейная по $f$ функция.

2) Пусть теперь $\xi_t + c(\xi_x - \tau_t) - c^2\tau_x\not\equiv 0$. Тогда, зафиксировав временно те значения $(t,x)$, при которых указанное выражение не нуль (соответствующие числовые значения функций мы будем обозначать звездочкой), получаем, что решение уравнения (32) имеет вид
$$\mathcal{F} = e^{I_2}\left(G(fe^{I_1}) - I_3\right),\quad I_1 = \int \frac{(3c\tau^*_x - 2\xi^*_x + \tau^*_t)dc}{c^2\tau^*_x - c(\xi^*_x - \tau^*_t) - \xi^*_t}, \,\,\, I_2 = \int \frac{(3c\tau^*_x - \xi^*_x + 2\tau^*_t)dc}{c^2\tau^*_x - c(\xi^*_x - \tau^*_t) - \xi^*_t}, $$
$$I_3 = \int \frac{\xi^*_{tt} + c(2\xi^*_{tx} - \tau^*_{tt}) + c^2(\xi^*_{xx} - 2\tau^*_{tx}) - c^3\tau^*_{xx}}{c^2\tau^*_x - c(\xi^*_x - \tau^*_t) - \xi^*_t} e^{-I_2}dc.$$
Таким образом, и в том, и в другом случае мы получили следующий предварительный анзац для $\mathcal{F}$:
$$\mathcal{F} = u(c)G(fv(c)) + w(c). \eqno(33)$$
В нем, в силу предположения теоремы ($\mathcal{F}$ зависит от $f$) функции $u(c)$ и $v(c)$ предполагаются ненулевыми, а $G(\cdot)$ -- непостоянной.

Подстановка предварительного анзаца (33) в (32), дает, после введения обозначения $z = fv(c)$, равенство
$$\left(3c\tau_x - 2\xi_x + \tau_t + (\xi_t + c(\xi_x - \tau_t) - c^2\tau_x)\frac{v_c}{v}\right)uzG_z + $$
$$+\left(3c\tau_x -\xi_x + 2\tau_t + (\xi_t + c(\xi_x - \tau_t) - c^2\tau_x)\frac{u_c}{u} \right)uG =$$
$$ = \xi_{tt} + c(2\xi_{tx} - \tau_{tt}) + c^2(\xi_{xx} - 2\tau_{tx}) - c^3\tau_{xx} - (\xi_t + c(\xi_x - \tau_t) - c^2\tau_x)w_c - (3c\tau_x -\xi_x + 2\tau_t)w.$$

I. Если $3c\tau_x - 2\xi_x + \tau_t + (\xi_t + c(\xi_x - \tau_t) - c^2\tau_x)(\ln v)_c = 0$, то, поскольку в силу условий теоремы $G$ -- не константа, коэффициент при $G$ -- тоже нулевой, а значит нулевой и свободный член, в (33) функция $G(\cdot)$ произвольная, а $u$, $v$ и $w$ удовлетворяют системе уравнений
$$3c\tau_x -\xi_x + 2\tau_t + (\xi_t + c(\xi_x - \tau_t) - c^2\tau_x)(\ln |u|)_c = 0,\eqno(34)$$
$$3c\tau_x - 2\xi_x + \tau_t + (\xi_t + c(\xi_x - \tau_t) - c^2\tau_x)(\ln |v|)_c = 0,\eqno(35)$$
$$\xi_{tt} + c(2\xi_{tx} - \tau_{tt}) + c^2(\xi_{xx} - 2\tau_{tx}) - c^3\tau_{xx} = (\xi_t + c(\xi_x - \tau_t) - c^2\tau_x)w_c + (3c\tau_x -\xi_x + 2\tau_t)w.\eqno(36)$$

II. Если $3c\tau_x - 2\xi_x + \tau_t - (\xi_t + c(\xi_x - \tau_t) - c^2\tau_x)(\ln v)_c \ne 0$, а коэффициент при $G$ нулевой, то $G = G_1\ln |z|$ и $\mathcal{F} = u(c)\ln |f| + w(c)$, где
$$3c\tau_x -\xi_x + 2\tau_t + (\xi_t + c(\xi_x - \tau_t) - c^2\tau_x)(\ln |u|)_c = 0,\eqno(37)$$
$$(\xi_t + c(\xi_x - \tau_t) - c^2\tau_x)w_c + (3c\tau_x -\xi_x + 2\tau_t)w + (3c\tau_x - 2\xi_x + \tau_t)u = \xi_{tt} + c(2\xi_{tx} - \tau_{tt}) + c^2(\xi_{xx} - 2\tau_{tx}) - c^3\tau_{xx}.\eqno(38) $$

III. Коэффициенты при $G$ и $G_z$ в уравнении ненулевые, тогда $G = G_1z^{j} + G_2$ и $\mathcal{F} = u(c)f^j + w(c)$
$$3c\tau_x -\xi_x + 2\tau_t + j(3c\tau_x - 2\xi_x + \tau_t) + (\xi_t + c(\xi_x - \tau_t) - c^2\tau_x)(\ln |u|)_c = 0,\eqno(39)$$
$$(\xi_t + c(\xi_x - \tau_t) - c^2\tau_x)w_c + (3c\tau_x -\xi_x + 2\tau_t)w = \xi_{tt} + c(2\xi_{tx} - \tau_{tt}) + c^2(\xi_{xx} - 2\tau_{tx}) - c^3\tau_{xx}.\eqno(40)$$
Таким образом, для первого семейства функций $\mathcal{F}=\mathcal{F}(c,f)$ мы получили три анзаца, с которыми и будем далее работать. Первый из них мы, для удобства, будем называть {\em общим}, второй -- {\em логарифмическим}, а третий -- {\em степенным}.

\subsubsection{Первичная классификация} Можно заметить, что уравнение
(39) -- общее для всех трех случаев (уравнения (34) и (37) для случаев I и II получаются из него полаганием $j=0$), поэтому естественно за основание классификации взять именно его.

{\bf Лемма 2.} {\em Уравнение (39) как уравнение относительно $\tau$, $\xi$ имеет решения, отличные от констант, только для функций $u(c)$, приведенных в таблице 7:
\begin{center}
{\rm Таблица 7. Первичная классификация: отличные от констант решения $(\tau,\xi)$ уравнения (39).}

\begin{tabular}{|c||c|c|c|}
\hline
№ & $j$ & $u$& $\tau$, $\xi$ \\
\hline
\hline
\hline
1&  $j =1$ & $u(c)$ произвольная & $\tau=at+b$, $\xi=ax+l$,  \\
\hline
2 & любое &$u=u_1(c-k)^{j+2}(c+k)^{2j+1}$ &$\tau=\tau(x+kt)$, $\xi=k\tau(x+kt)+l$\\
\hline
2.1 &$j=1$ & $u=u_1(c^2+n)^3$ & $\tau_t=\xi_x$, $n\tau_x+\xi_t=0$\\
\hline
3 & любое & $u = u_1c^m$ & $\tau = (2j+1-m)\lambda t + b$, \\
& & $m \ne 3(j+1), j+2, 2j+1, 0$ &$\xi = (2+j-m)\lambda x + l$ \\
\hline
3.1 & $j \ne -1, -1/2$ & $u=u_1c^{3(j+1)}$ & $\xi=\xi(x)$, $(2j+1)\tau = (2+j)\xi_xt + b(x)$ \\
\hline
3.2 & $j \ne -2, -1/2, 1$ & $u = u_1c^{j+2}$& $\tau = \tau(t)$, $\xi = l$\\
\hline
3.3 & $j \ne -2, -1/2, 1$ & $u = u_1c^{2j+1}$& $\tau = b$, $\xi = \xi(x)$\\
\hline
3.4 & $j = -1/2$ & $u=u_1c^{3/2}$& $\tau = \tau(t,x)$, $\xi = {\rm const}$\\
\hline
3.5 & $j = 1$ & $u = u_1c^3$& $\tau = \tau(t)$, $\xi = \xi(x)$\\
\hline
4 & любое & $u = u_1e^{mc}$ & $\tau = at + b$,  \\
& & &$\xi = a(x + \frac{j-1}{m}t) + l$ \\
\hline
5 & $j \ne -1, -1/2$ & $u = u_1$ & $\tau=\tau(t)$, $(2j+1)\xi = (2+j)\tau_t x + l(t)$ \\
\hline
5.1 &  $j = -1$ & $u = u_1$ & $\tau=\tau(t,x)$, $\xi=\xi(t,x)$,  \\
& & &$\tau_t + \xi_x = 0$ \\
\hline
5.2 &  $j = -1/2$ & $u = u_1$ & $\tau={\rm const}$, $\xi=\xi(t,x)$,  \\
\hline
\end{tabular}
\end{center}}

{\bf Доказательство.} Рассмотрим (39) как уравнение на функцию $u$. Если в этом уравнении коэффициент при $(\ln u)_c$ равен нулю тождественно, то равен нулю и свободный член, и мы получаем случай 1.

Если же этот коэффициент не равен нулю, то, так как $u$ не зависит от $t$, $x$, то можно, временно зафиксировав соответствующие значения $t,x$ (и, соответственно, $\tau$, $\xi$ и их производных), получить, что уравнение относительно $u$ имеет вид
$$(\ln u)_c = \frac{3k(j+1)c - p(2j+1) + q(j+2)}{kc^2 - (p - q)c - r},\eqno(41)$$
где $k$, $p$, $q$, $r$ -- некоторые константы. Это дает нам анзац для $u$, который мы далее будем подставлять обратно в (39), разрешая его уже относительно $\tau$ и $\xi$. При этом имеет место несколько вариантов:

A) Если $k \ne 0$, и дробь несократима (случай сократимой дроби, как и случай $k=0$ будет рассматриваться в пункте б) ), то заменой переменных $\bar{t} = t$, $\bar{x} = x - (p-q)t/(2k)$ мы можем упростить формулу для $u$, приведя ее к виду:
$$(\ln u)_c = \frac{3(j+1)c + m}{c^2 + n}.\eqno(42)$$

Б) Если в (41) $k = 0$, но $p-q\ne 0$, то заменой переменных $\bar{t} = t$, $\bar{x} = (p-q)x + rt$ мы можем упростить (41), приведя его к виду
$$(\ln u)_c = {m}/{c} \eqno(43)$$
к этому же виду приводится (41) и в случае $k\ne 0$, но дробь сократима. Мы здесь будем предполагать, что $m\ne 0$.

В) Если $k=p-q=0$, то дробь в (41) сводится к константе $$(\ln u)_c = m. \eqno(44)$$

Рассмотрим каждый из этих случаев.

A) Случай {\bf функции (42)}. Считая, что (42) -- несократимая дробь, подставим (42) в (39), получим:
$$[3c\tau_x -\xi_x + 2\tau_t + j(3c\tau_x - 2\xi_x + \tau_t)](c^2 + n) + (\xi_t + c(\xi_x - \tau_t) - c^2\tau_x)(3(j+1)c + q) = 0.$$
Приведем подобные и приравняем к нулю коэффициенты при различных степенях $c$. Это дает три равенства:
$$-(2j+1)\tau_t + (j+2)\xi_x - q\tau_x =0,$$ $$3n(j+1)\tau_x + 3(j+1)\xi_t + q\xi_x - q\tau_t =0,\eqno(45)$$ $$(2+j)n\tau_t - (2j+1)n\xi_x + q\xi_t = 0.$$

Мы получили систему вида
$$a_{1j}\tau_t+a_{2j}\tau_x+a_{3j}\xi_t+a_{4j}\xi_x=0, \qquad j=1,2,3,$$
которой соответствует матрица
$$\left(\begin{array}{cccc}
-(2j+1)&-q&0&j+2\\
-q&3n(j+1)&3(j+1)&q\\
n(j+2)&0&q&-n(2j+1)\\
\end{array}\right)\eqno(46)$$
Основные варианты в дальнейшем рассмотрении связаны с рангом этой матрицы.

Нулевого ранга у системы быть не может: величины $2j+1$ и $j+2$ одновременно в нуль не обращаются.
Ранг равен единице, если $j+1=q=0$, но это случай сократимой дроби в (42), мы его не рассматриваем.

Ранг равен двум, если равны нулю все миноры третьего порядка. Первый из них совпадает с четвертым, а второй равен третьему, умноженному на $n$. Поэтому у нас остается два уравнения:
$$q[9n(j+1)^2+q^2]=0,\quad (j-1)[9n(j+1)^2+q^2]=0.$$
Если равен нулю второй множитель, то либо $j+1=q=0$, либо  $n=-\frac{q^2}{9(j+1)^2}$, в обоих случаях дробь в (42) сократима, что снова противоречит нашему предположению. Значит, ранг системы равен двум при $q=j-1=0$. Во всех остальных случаях ранг матрицы равен трем.

В случае $q=j-1=0$ получаем $u=A(c^2+n)^3$, и мы имеем два уравнения: $\tau_t=\xi_x$, $n\tau_x+\xi_t=0$, то есть случай 2.1 из формулировки леммы.

Во всех остальных случаях, когда ранг матрицы (46) равен трем, какие-то три из производных $\tau_t$, $\tau_x$, $\xi_t$, $\xi_x$ выражаются через четвертую. Одно из этих соотношений содержит одну и ту же функцию, и поэтому эта функция -- функция от линейной комбинации $t,x$. Но тогда, в силу двух других соотношений, и другая функция имеет такой же вид. Подстановка $\tau=\tau(\alpha t+\beta x)$, $\xi=\xi(\alpha t+\beta x)$ в (45) дает три линейных уравнения относительно $\tau'$ и $\xi'$:
$$-[(2j+1)\alpha + q\beta]\tau'+ (j+2)\beta\xi'  = 0,$$
$$[3n(j+1)\beta - q\alpha]\tau'+ [3(j+1)\alpha + q\beta]\xi'  = 0,$$
$$(2+j)n\alpha\tau' + [q\alpha- (2j+1)n\beta]\xi' = 0.$$

Ранг этой системы относительно $\tau'$ и $\xi'$  равен нулю, если все коэффициенты равны нулю. В этом случае если $j\ne -2$, то мы получаем $\alpha=\beta=0$, то есть $\tau$  и $\xi$ являются константами, алгебра двумерна. Если же $j=-2$, то из уравнений $-3\alpha + q\beta= 0$, $3n\beta + q\alpha=0$
мы либо снова получаем $\alpha=\beta=0$, либо $n=-q^2/9$, что соответствует случаю, когда дробь в (42) сократима. Таким образом, в этом случае более, чем двумерных алгебр мы не получаем.

Если ранг этой системы равен двум, то мы получаем $\tau'=\xi'=0$, что также приводит нас к двумерной алгебре.

Остается случай, когда ранг матрицы (46) равен единице. Для этого необходимо, чтобы были равны нулю все миноры второго порядка. То есть
$$[(2j+1)\alpha + q\beta][3(j+1)\alpha + q\beta]+ (j+2)\beta[3n(j+1)\beta - q\alpha]= 0,$$
$$[(2j+1)\alpha + q\beta][q\alpha- (2j+1)n\beta] = 0,$$
$$[3n(j+1)\beta - q\alpha][q\alpha- (2j+1)n\beta]-[3(j+1)\alpha + q\beta](j+2)n\alpha= 0.$$
Упростим уравнения
$$3(j+1)(2j+1)\alpha^2 + q(4j+2)\alpha\beta + [q^2 + 3n(j+1)(j+2)]\beta^2 = 0,$$
$$q(2j+1)\alpha^2 + [q^2 - 3n(j+1)(j-1)] \alpha\beta - nq(2j+1)\beta^2  = 0,$$
$$[q^2 + 3n(j+1)(j+2)]\alpha^2 - qn(4j+2)\alpha\beta + 3n^2(2j+1)(j+1)\beta^2 = 0.$$
Комбинируя первое и третье уравнения, получаем
$$[q^2 + 9n(j+1)^2](\alpha^2+ n\beta^2)=0,$$
откуда (поскольку равенство нулю первого множителя приводит к случаю сократимой дроби в (42)), $\alpha^2+ n\beta^2=0$.
Очевидно, что ненулевое решение здесь возможно только при $n\leq 0$, $n=-k^2$ и тогда $\alpha=k\beta$. Подставляя это в остальные уравнения, получаем, что нетривиальное решение (то есть не сводящееся к $\alpha=\beta=0$) возможно только если
$$q^2+2q(2j+1)k+3k^2(j+1)(j-1)= 0
\quad\text{или}\quad (q+3k(j+1))(q+k(j-1))=0.$$
Снова исключая равенство нулю первого сомножителя, приводящий к сократимой дроби в (42), приходим к окончательному условию: $q=-k(j-1)$. При этом $u=(c-k)^{j+2}(c+k)^{2j+1}$, $\tau=\tau(x+kt)$, $\xi=k\tau(x+kt)+C$. Это дает нам случай 2 из формулировки леммы.
\\

Б) Случай {\bf функции (43)}. Предполагая, что $m\ne 0$, подставим (43) в (39), получим:
$$[3c\tau_x -\xi_x + 2\tau_t + j(3c\tau_x - 2\xi_x + \tau_t)]c + m(\xi_t + c(\xi_x - \tau_t) - c^2\tau_x) = 0.$$
Расщепление этого равенства по $c$ дает
$$(3+3j-m)\tau_x = 0, \qquad (2+j-m)\tau_t = (2j+1-m)\xi_x, \qquad m\xi_t = 0.$$
Так как $m\ne 0$, то $\xi_t = 0$, то есть $\xi=\xi(x)$.

Рассмотрим несколько вариантов:

Если $m\ne 3(1+j)$, то $\tau_x = 0$, $\tau=\tau(t)$, и из уравнения $(2+j-m)\tau_t = (2j+1-m)\xi_x$ следует, что и слева, и справа стоит константа. Обозначив ее через $\lambda$, получаем $(2+j-m)\tau = \lambda t + b$, $(2j+1-m)\xi = \lambda x + l$. Это случай 3 (он имеет место, если $m\ne j+2, 2j+1$), либо случаи 3.2 (при $m=j+2$) или 3.3 (при $m=2j+1$). Одновременное выполнение равенств $m=j+2=2j+1$ возможно только при $j=1$, и это дает нам случай 3.5.\\

Если же $m=3(1+j)$, мы получаем уравнение $(2j+1)\tau_t = (2+j)\xi'(x)$, откуда $(2j+1)\tau = (2+j)\xi'(x)t + b(x)$ и это случай 3.1 (при $j\ne -1/2$) и 3.4 (при $j=-1/2$).
\\
\\

В) Случай {\bf функции (44)}. Подставим (44) в (39), получим:
$$3c\tau_x -\xi_x + 2\tau_t + j(3c\tau_x - 2\xi_x + \tau_t) + m(\xi_t + c(\xi_x - \tau_t) - c^2\tau_x) = 0,$$
откуда, расщепляя по $c$, можно вывести
$$m\tau_x = 0, \qquad 3(j+1)\tau_x + m\xi_x - m\tau_t = 0, \qquad (2+j)\tau_t - (2j+1)\xi_x + m\xi_t = 0.$$

Если $m\ne 0$, то $\tau_x = 0$, $\xi_x = \tau_t$ и $(1-j)\tau_t + m\xi_t = 0$. Из второго и третьего уравнений следует $\tau_{tt}=\xi_{tx}=\frac{j-1}m\tau_{tx}=0$, а значит $\tau = at + b$. В итоге получаем $\tau = at + b$, $\xi = a(x + \frac{j-1}{m}t) + l$. Это случай 4.\\

Если же $m=0$, то $(j+1)\tau_x = 0$, $(2+j)\tau_t - (2j+1)\xi_x = 0$. Отсюда для $j\ne-1$ получаем $\tau_x = 0$ (то есть $\tau=\tau(t)$) и $(2j+1)\xi_x = (2+j)\tau_t$. Это случаи 5 (если $j\ne -1/2$) и 5.2 (если $j=-1/2$); случай же $j=-1$ дает нам уравнение $\tau_t + \xi_x = 0$ для $\tau(t,x)$, $\xi(t,x)$ и $u = {\rm const}$, это случай 5.1.

Лемма доказана.\\

Возвращаясь к задаче классификации, заметим, что замена переменных $\bar{t}=x$, $\bar{x}=t$ устанавливает эквивалентность случаев 3.1 и 5, 3.4  и 5.2, 3.2 и 2.3.   Поэтому случаи 3.1, 3.3 и 3.4 мы рассматривать далее не будем.

Напомним, что для анзацев I и II уравнение (39) рассматривается только при $j=0$, поэтому для них мы будем рассматривать только случаи 2, 3, 3.2, 4, 5.

Начнем подстановку указанных алгебр в оставшиеся уравнения, соответствующие полученным нами трем анзацам.

\subsubsection{Анзац I. $\mathcal{F} = u(c)G(fv(c)) + w(c)$, $G(\cdot)$ -- произвольная функция}

Напомним оставшиеся уравнения (35) и (36):
$$3c\tau_x - 2\xi_x + \tau_t + (\xi_t + c(\xi_x - \tau_t) - c^2\tau_x)(\ln |v|)_c = 0,$$
$$\xi_{tt} + c(2\xi_{tx} - \tau_{tt}) + c^2(\xi_{xx} - 2\tau_{tx}) - c^3\tau_{xx} = (\xi_t + c(\xi_x - \tau_t) - c^2\tau_x)w_c + (3c\tau_x -\xi_x + 2\tau_t)w.$$

Начнем со {\bf случая 2}. Подстановка $\tau$, $\xi$ при $j=0$ дает
$$\tau'[3c - k + (k^2 - c^2)(\ln |v|)_c] = 0, \quad \tau''(k + c)^2(k - c) = \tau'[(k^2 - c^2)w_c + (3c + k)w].$$
Если $\tau'=0$, то алгебра двумерна. Если же $\tau'\ne 0$, то $v = A(c - k)(c + k)^2$,
из второго уравнения следует $\tau'' = \lambda\tau'$ и $w= (c^2-k^2)[Q(c-k) - \lambda]$. Мы получили функцию
$$\mathcal{F} = (c-k)^2(c+k)G((c+k)^2(c-k)f)-\lambda (c^2-k^2).\eqno(47)$$
Замена переменных $\bar t=x+kt$, $\bar x=x-kt$ при $k\ne 0$ переводит эту функцию в функцию $\bar{\mathcal{F}} = \bar c^2 \left[2kG(4k^2\bar c f)+\lambda\right]-\lambda \bar c$, то есть в функцию вида III.3.2.

Если же $k=0$, то замена $\bar t=x$, $\bar x=t$ приводит нас к функции $\mathcal{F}=-G(f)+\lambda c$, принадлежащей семейству III.4.2 при $\lambda\ne 0$ или III.2.1 при $\lambda=0$. Таким образом, в этом случае алгебр размерности более двух для уравнений, не входящих в семейства III.3.2, III.4.2 и III.2.1, нет. \\

В {\bf случае 3} либо $\lambda = 0$ и алгебра двумерна, либо $\lambda \ne 0$, и тогда подстановка $\tau$ и $\xi$ дает уравнения
$$m - 3 + c (\ln |v|)_c = 0, \quad cw_c - mw = 0,$$
то есть $v = v_1c^{3-m}$, $w = w_1 c^m$ при $(2-m)\tau = \lambda t + b$, $(1-m)\xi = \lambda x + l$ (напомним что в случае 3 при $j=0$ действует ограничение $m\ne 0,1,2,3$). Мы получили функцию $\mathcal{F}=c^m[G(c^{3-m}f)+w_1]$ из семейства III.3.2.\\

В {\bf случае 3.2} подстановка $\xi=l$, $\tau=\tau(t)$ в уравнения для $v$ и $w$ дает соотношения
$$ \tau_t[1 - c (\ln |v|)_c] = 0, \quad \tau_t(-cw_c+2w) = - c\tau_{tt}.$$
Алгебра размерности больше двух получается только если $\tau_t\ne 0$,
и тогда $\tau_{tt}=\lambda\tau_t$, $v=v_1c$, $w=w_1c^2-\lambda c$, и мы получаем $\mathcal{F}=c^2[G(cf)+w_1]-\lambda c$ из семейства III.3.2.\\

В {\bf случае 4} при $a=0$ мы получаем двумерную алгебру. Если же $a \ne 0$, то
$$1 + \frac{1}{m}(\ln |v|)_c = 0, \qquad w - \frac{1}{m}w_c = 0.$$
Отсюда $v = v_1 e^{-mc}$, $w = w_1 e^{mc}$ (напомним, что $m\ne 0$), $\tau = at + b$, $\xi = a(x - \frac{1}{m}t) + l$. Это функция $\mathcal{F} = e^{mc}G(e^{-mc}f)$ из семейства III.3.3.\\

В {\bf случае 5} подстановка $\tau$, $\xi$ при $j=0$ в уравнения позволяет каждое уравнение рассматривать как многочлен по $x$, тождественно равный нулю. Это дает четыре уравнения:
$$ - 3\tau' + (l'+ c\tau')(\ln |v|)_c = 0,\qquad 2\tau''(\ln |v|)_c=0,$$
$$l''+ 3c\tau''= (l'+ c\tau')w_c,\qquad 2\tau'''=2\tau''w_c.$$
Прежде всего, заметим, что $\tau''=0$: иначе  из первого уравнения следует, что $(\ln |v|)_c = 0$, а тогда из первого $\tau'= 0$, что противоречит предположению.

Но раз $\tau''=0$, то $\tau=at+b$, и у нас остается два уравнения
$$ - 3a + (l'+ ac)(\ln |v|)_c = 0,\qquad l''= (l'+ ac)w_c.$$
Если $l'$ и $a$ одновременно обращаются в нуль, то алгебра двумерна. Если же нет -- из первого уравнения мы получаем, что либо $a=0$, $v=v_1$, и тогда из второго $l''=\lambda l'$, $w=\lambda c$ (что дает нам функцию $\mathcal{F} = G(f)+\lambda c$, принадлежащую семейству III.4.2 при $\lambda\ne 0$ или III.2.1 при $\lambda=0$). Либо $l'={\rm const}$, $l=l_1t+l_2$, $v=(ac+l_1)^3$, $w=w_1$. Полученная функция $\mathcal{F}=G((c+k)^3f)$ заменой $\bar t=t$, $\bar x=x+kt$ преобразуется к функции III.3.2 (при $m=0$).
Таким образом, мы перебрали все случаи из леммы 1, подходящие для этого анзаца, и убедились, что алгебра размерности больше двух встречается только у функций семейства III.*.*.

\subsubsection{Анзац II. $\mathcal{F} = u(c)\ln f + w(c)$}

Напомним оставшееся уравнение (38):
$$(\xi_t + c(\xi_x - \tau_t) - c^2\tau_x)w_c + (3c\tau_x -\xi_x + 2\tau_t)w + (3c\tau_x - 2\xi_x + \tau_t)u = \xi_{tt} + c(2\xi_{tx} - \tau_{tt}) + c^2(\xi_{xx} - 2\tau_{tx}) - c^3\tau_{xx}.\eqno(38)$$

{\bf Случай 2} при подстановке в (38) $\tau=\tau(x+kt)$, $\xi=k\tau(x+kt)+l$ дает нам равенство
$$\tau'[(k^2- c^2)w_c + (3c+k)w +(3c-k)u_1(c-k)^2(c+k)] = \tau''(k^2+ck^2-c^2k-c^3).$$
Если $\tau'=0$, то алгебра получается двумерной. Если же $\tau'\ne 0$, то в уравнении можно совершить разделение переменных: $\tau'' = \lambda\tau'$,
$$(k^2- c^2)w_c + (3c+k)w +(3c-k)u_1(c-k)^2(c+k) = \lambda (k^2+ck^2-c^2k-c^3).$$

Если $\lambda \ne 0$, то $\tau = ae^{\lambda(x+kt)} + b$, $\xi = ake^{\lambda(x+kt)} + l$, $w= (c^2-k^2)(c-k)\{w_1 + u_1\ln[(c+k)^2(c-k)]\} - \lambda (c^2-k^2)$.

Если же $\lambda = 0$, то $\tau = a(x+kt) + b$, $\xi = ak(x+kt) + l$ и $w=(c^2-k^2)(c-k)\{w_1 + u_1\ln[(c+k)^2(c-k)]\} $.

Мы получили функцию
$\mathcal{F}= (c^2-k^2)(c-k)\{u_1\ln |f(c+k)^2(c-k)|+w_1\} - \lambda (c^2-k^2)$, то есть
частный случай функции (47) из предыдущего анзаца (для $G(z)=u_1\ln z+w_1$), которая, как уже говорилось, сводится к III.3.2.
\\

В {\bf случае 3} $\lambda=0$ приводит к двумерной алгебре. Если же $\lambda\ne 0$, то уравнение (38) приобретает вид
$$cw_c - mw + (m-3)u_1c^m = 0.$$
Отсюда $w = w_1c^m - u_1(m-3)c^m\ln c$ при $\tau = (m-1)\lambda t + b$, $\xi = (m-2)\lambda x + l$. И здесь мы получаем лишь частный случай (опять же для $G(z)=P\ln z+Q$), функции III.3.2, но только в другой форме: $\mathcal{F}= c^m\{u_1\ln |fc^{3-m}|+w_1\}$.\\

В {\bf случае 3.2} мы подстановкой в уравнение (38) $u = u_1c^2$, $\tau = \tau(t)$, $\xi = l$ приводим его к виду
$$\tau_t(-cw_c + 2 w + u_1c^2)= - c\tau_{tt}.$$
Алгебра размерности более двух возможна только если $\tau_t\ne 0$, а тогда $\tau_{tt}=\lambda\tau_t$, $w=w_1c^2+u_1c^2\ln c -\lambda c$, что снова приводит нас к функции $\mathcal{F}= c^2\{u_1\ln |cf|+w_1\}-\lambda c$ из семейства III.3.2.\\

{\bf Случай 4} дает нам при $a=0$ двумерную алгебру, а при $a\ne 0$ (и $m\ne 0$) -- уравнение $w_c - mw + mu_1e^{mc} = 0$, из которого получаем $w = w_1e^{mc} - mu_1ce^{mc}$, $\mathcal{F}= e^{mc}\{u_1\ln |fe^{-mc}|+w_1\}$, то есть частный случай функции III.3.3. \\

В {\bf случае 5} подстановка $u=u_1$, $\tau=\tau(t)$, $\xi = 2\tau' x + l(t)$ приводит уравнение к виду
$$(2\tau'' x + l' + c\tau')w_c - 3\tau' u_1 = 2\tau''' x + l'' + 3c\tau'',$$
и оно расщепляется по степеням $x$:
$$\tau'' w_c = \tau''', \qquad (l' + c\tau')w_c - 3\tau' u_1 = l'' + 3c\tau''.$$
Из этих уравнений следует, что $\tau''=0$: иначе из первого $w_c={\rm const}$, а тогда из второго получаем $\tau''=0$, то есть противоречие предположению. Таким образом, $\tau=at+b$, $\xi=2ax+l(t)$ и у нас остается одно уравнение
$(l' + ac)w_c - 3a u_1 = l'',$
из которого следует, что $l''=\alpha l'+\beta$, где $\alpha,\beta$ -- некоторые константы.

Подставляя эту формулу в наше уравнение, мы получаем
$(l' + ac)w_c - 3a u_1 = \alpha l'+\beta,$
откуда либо $l'$ -- константа, либо уравнение распадается на два:
$$w_c = \alpha ,\qquad acw_c - 3a u_1 = \beta,$$
в последнем случае мы получаем $a\alpha=0$, и либо $\alpha=0$, $w=w_1$, $\beta=-3au_1$, $l=-\frac 32au_1t^2+l_1t+l_2$ (здесь мы получаем функцию $\mathcal{F}= u_1\ln |f|+w_1$ из семейства II.5), либо $a=\beta=0$, $w=\alpha c+w_1$, $l=l_1e^{\alpha t}+l_2$, $\tau=b$, $\xi=l_1e^{\alpha t}+l_2$ -- что приводит нас, после сдвига $c$ заменой $\bar t=t$, $\bar x=x+\frac {w_1}\alpha t$, к $\mathcal{F}= u_1\ln |f|+\alpha c$ из семейства III.4.2 при $\alpha\ne 0$ или III.2.1 при $\alpha=0$.

Если же $l'$ -- константа, то есть $l=l_1t+l_2$, то уравнение для $w(c)$ приобретает вид
$(l_1+ ac)w_c - 3a u_1 =0$, и имеет решением при $a\ne 0$ функцию $w=3u_1\ln(ac+l_1)+w_1$ (что дает в итоге $\mathcal{F}= u_1\ln |fc^3|+w_1$, то есть частный случай функции III.3.2), а при $a=0$ -- константу $w\equiv w_1$ (что приводит тоже к частному случаю, но уже функции III.2.1).

Таким образом, и для этого анзаца мы алгебры размерности больше двух получаем только для функций из семейства III.*.*.

\subsubsection{Анзац III. $\mathcal{F} = u(c)f^j + w(c)$}

Здесь, поскольку $j$ может принимать, вообще говоря, любые значения, нам придется рассмотреть уже десять случаев из леммы 1 -- 1, 2, 2.1, 3, 3.2, 3.5, 4, 5, 5.1, 5.2.

Для классификации у нас остается одно уравнение (40):
$$(\xi_t + c(\xi_x - \tau_t) - c^2\tau_x)w_c + (3c\tau_x -\xi_x + 2\tau_t)w = \xi_{tt} + c(2\xi_{tx} - \tau_{tt}) + c^2(\xi_{xx} - 2\tau_{tx}) - c^3\tau_{xx}.\eqno(40)$$

Начнем со {\bf случая 1}: $j=1$, $u(c)$ -- произвольная, $\tau=at+b$, $\xi=ax+l$. Подставляя их в уравнение (40) для $w(c)$, получаем $aw = 0$, откуда либо $a=0$, и мы получаем двумерную алгебру, либо $w(c)=0$ и мы приходим к семейству III.2.3.\\

Перейдем к {\bf случаю 2}. Подстановка в уравнение (40) $\tau=\tau(x+kt)$, $\xi=k\tau(x+kt)+l$ приводит его к виду
$$\tau'[(k^2 - c^2)w_c + (3c + k)w] = \tau''(k + c)(k^2 - c^2).$$
Предположение $\tau'=0$ приводит к двумерной алгебре. Если же $\tau$ не константа, то, поскольку зависимость от $x+kt$ содержится только в функции $\tau$, получаем $\tau'' = \lambda\tau'$. Если $\lambda \ne 0$, то $\tau = ae^{\lambda(x+kt)} + b$, $\xi = ake^{\lambda(x+kt)} + l$, $w= w_1(c^2-k^2)(c-k) - \lambda(c^2-k^2)$, и мы приходим к частному случаю функции (47) из первого анзаца, приводящейся к III.3.2.  Если же $\lambda = 0$, то $\tau = a(x+kt) + b$, $\xi = ak(x+kt) + l$ и $w=w_1(c^2-k^2)(c-k)$, это также дает частный случай функции (47). \\

В {\bf случае 2.1} подстановка $\tau_t=\xi_x$ и $n\tau_x+\xi_t=0$ (то есть $\tau_{tt} = -n\tau_{xx}$) в уравнение при $j=1$, приводит это уравнение к виду.
$$(-n - c^2)\tau_x w_c + (3c\tau_x + \tau_t)w = (\tau_{tx} + c\tau_{xx})(-n - c^2). \eqno(48)$$

Если $\tau_x = 0$, то мы получаем, что $\tau_t w = 0$, и тогда либо $\tau_t = 0$ и алгебра двумерна, либо $\tau_t$ не равно нулю, и тогда $w \equiv 0$ и мы получаем функцию II.2 (при $n\ne 0$) или линейный случай функции II.1 (при $n=0$).

Если же $\tau_x \not\equiv 0$, то поделим (48) на $\tau_x$:
$$(-n - c^2)w_c + (3c + \frac{\tau_t}{\tau_x})w = (\frac{\tau_{tx}}{\tau_x} + c\frac{\tau_{xx}}{\tau_x})(-n - c^2)$$
и продифференцируем полученное уравнение по $t$ и по $x$:
$$\left(\frac{\tau_t}{\tau_x}\right)_t w = -\left(\left(\frac{\tau_{tx}}{\tau_x}\right)_t + c\left(\frac{\tau_{xx}}{\tau_x}\right)_t\right)(n + c^2),$$
$$\left(\frac{\tau_t}{\tau_x}\right)_x w = -\left(\left(\frac{\tau_{tx}}{\tau_x}\right)_x + c\left(\frac{\tau_{xx}}{\tau_x}\right)_x\right)(n + c^2).$$
Если $\left(\frac{\tau_t}{\tau_x}\right)_t =\left(\frac{\tau_t}{\tau_x}\right)_x \equiv 0$ получаем, что $\tau_t = k\tau_x$ (то есть $\tau = \tau(x+kt)$), что возможно только если $n=-k^2$, и тогда $\xi=k\tau +l$. При этом в (48) переменные разделяются: $\tau_{xx} = a\tau_x$,
$$-(k^2-c^2)w_c + (3c + k)w = -a(k + c)(k^2-c^2),$$
откуда
$$w = w_1(c-k)^2(c+k)- a(c^2-k^2),$$
и таким образом мы опять получаем частный случай функции (47).

Если же $\left(\frac{\tau_t}{\tau_x}\right)_t$ и $\left(\frac{\tau_t}{\tau_x}\right)_x$ не равны тождественно нулю, то, выбрав и зафиксировав соответствующее значение $(t,x)$, мы из по крайней мере одного из полученных уравнений получим, что $w = (n + c^2)(Qc+R)$. При $n=\pm k^2\ne 0$ это дает нам, после замены $\bar t=kRt+kQx$, $\bar x=Rx\mp k^2Qt$  функцию $\mathcal{F}= -\frac {u_1}{k^3}(c^2\pm 1)^3f-\frac 1{k^2}(c^2\pm 1)\}$ из семейства II.2, а при $n=0$ мы получаем функцию $\mathcal{F}= u_1c^6f+Qc^3+Rc^2$ из семейства II.1.\\

Рассмотрим теперь {\bf случай 3}: $u = u_1c^m$, $m \ne 3(j+1), j+2, 2j+1, 0$, $\tau = (2j+1-m)\lambda t + b$, $\xi = (2+j-m)\lambda x + l$. Подстановка в уравнение (40) приводит его к виду
$\lambda [c(1 - j)w_c + (3j-m)w ] = 0$.
При $j\ne 1$ для $w = w_1c^{\frac {m-3j}{1 - j}}$ переобозначением $k=\frac {m-3j}{1 - j}$ мы приводим функцию к виду $\mathcal{F}=u_1c^{k+j(3-k)}f^j+w_1c^k$, то есть к частному случаю функции III.3.2 (с той же алгеброй). При $j= 1$, поскольку $m\ne 3$, получаем $w = 0$, то есть снова частный случай функции III.3.2. \\

Перейдем к {\bf случаю 3.2}:
$j \ne -2, -1/2, 1$, $u = u_1c^{j+2}$ и $\tau = \tau(t)$, $\xi =l$. Подстановка приводит уравнение к виду $\tau_t(-cw_c + 2w) = - c\tau_{tt}$.
Алгебра имеет размерность больше двух, если $\tau'\ne 0$, а в этом случае
$\tau''=\lambda\tau'$, $-cw_c +2w = - \lambda c$, $w=w_1c^2-\lambda c$. Это дает нам $\mathcal{F}= c^2[u_1(cf)^j +w_1]-\lambda c$, то есть частный случай функции III.3.2.\\

Рассмотрим далее {\bf случай 3.5}: $j = 1$, $u = u_1c^3$, $\tau = \tau(t)$, $\xi = \xi(x)$.
Подстановка приводит уравнение (40) к виду $c(\xi_x - \tau_t)w_c + (-\xi_x + 2\tau_t)w = - c\tau_{tt} + c^2\xi_{xx}$.
Из этого уравнения следует, что $\tau_{tt}=\alpha \tau_t+\beta$, а $\xi_{xx}=\gamma \xi_{x}+\delta$.
Подставляя эти соотношения, мы получаем
$$c(\xi_x - \tau_t)w_c + (-\xi_x + 2\tau_t)w = - c(\alpha\tau_t+\beta) + c^2(\gamma \xi_x+\delta).$$
Далее для разделения переменных нужно разобрать четыре варианта: когда обе производные $\tau_t$ и $\xi_x$ -- не константы, когда одна из них константа, и когда обе они постоянные.

В первом варианте мы получаем три уравнения
$$-cw_c + 2w = - c\alpha, \qquad cw_c-w =c^2\gamma,\qquad - c\beta+ c^2\delta=0, $$
из которых следует, что $\beta=\delta=0$, а $w=\gamma c^2-\alpha c$, и мы получаем функцию $\mathcal{F}= u_1 c^3f+\gamma c^2-\alpha c$ из семейства II.2.

Если $\tau_t={\rm const}$, $\tau=at+b$, но $\xi_x\ne {\rm const}$ уравнений будет уже два:
$$c w_c - w = c^2\gamma, \qquad -acw_c+ 2aw = - c(\alpha a+\beta) + c^2\delta.$$
Из первого следует, что $w=w_1c+\gamma c^2$, и мы снова получаем функцию семейства II.2. Аналогично рассматривается симметричный вариант $\xi_x={\rm const}$, $\tau_t\ne{\rm const}$, там получается $w=w_1c^2-\alpha c$.

Наконец, когда обе производных постоянны, и $\tau=at+b$, $\xi=\lambda x+l$, мы из уравнения $c(\lambda- a)w_c + (-\lambda + 2a)w = 0$ получаем $w=w_1c^{\frac{\lambda-2a}{\lambda-a}}$. Обозначая показатель степени через $k$, мы приходим к функции $\mathcal{F}=Pc^3f+Qc^k=c^k(Pc^{3-k}f+Q)$ -- частный случай функции III.3.2.\\

Продолжим. {\bf Случай 4}: $u = u_1e^{mc}$, $m\ne 0$,  $\tau = at + b$, $\xi = a(x + \frac{j-1}{m}t) + l$. Здесь мы приходим к уравнению $a[(j-1)w_c + mw] = 0$. При $j\ne 1$ получаем $w = w_1e^{-\frac{mc}{j-1}}$. Обозначая $k=-\frac{m}{j-1}$, мы приводим нашу функцию к виду $e^{kc}[u_1(e^{-kc}f)^j+w_1]$, то есть к частному случаю функции III.3.3. При $j=1$ мы получаем $w=0$, и трехмерная алгебра $\tau = at + b$, $\xi = ax + l$ приводит нас к частному случаю функции III.2.3.
\\

{\bf Случай 5.} $j \ne -1, -1/2$, $u=u_1$, $\tau=\tau(t)$, $\xi = \frac{2+j}{2j+1}\tau' x + l(t)$
Подстановка в уравнение и его расщепление по переменной $x$ дают два равенства:
$$(2+j)\tau'' w_c = (2+j)\tau''',\qquad [(2j+1)l' + c(1-j)\tau']w_c + 3j\tau' w = (2j+1)l'' + 3c\tau''.\eqno (49)$$

Если $j\ne -2$, то из первого уравнения получаем, что $\tau'''=\alpha \tau''$, и либо $\tau''\ne 0$, тогда $w=w_1+\alpha c$, второе уравнение (49) расщепляется на два
соотношения
$$(2j+1)\tau'\alpha  = 3\tau'', \qquad (2j+1)l'\alpha + 3j\tau' w_1 = (2j+1)l'',$$
из которых первое возможно только при $j=1$ (поскольку $\alpha\ne 0$ в силу предположения $\tau''\ne 0$). Мы получаем функцию $\mathcal{F}= u_1 f+w_1+\alpha c$ из семейства  II.1. Либо $\tau''=0$, $\tau=at+b $, и тогда уже во втором уравнении (49)
$$[(2j+1)l' + c(1-j)a]w_c + 3jaw = (2j+1)l''$$
мы можем использовать анзац $l''=\alpha l'+\beta$, преобразовав его к виду
$$[(2j+1)l' + c(1-j)a]w_c + 3ajw = (2j+1)(\alpha l'+\beta)$$
и перейдя ко второму слою вариантов.

Здесь если $l'$ не константа, и тогда уравнение расщепляется на
$$w_c  = \alpha ,\qquad c(1-j)aw_c + 3ajw = (2j+1)\beta,$$
из первого получаем $w=\alpha c+w_1$, подстановка во второе и расщепление уже по $c$ дает два равенства: $(1+2j)a\alpha=0$, $3ajw_1 = (2j+1)\beta$. Независимо от ограничений, налагаемых этими равенствами, мы все равно имеем дело с функцией вида $\mathcal{F}= u_1 f^j +\alpha c+w_1$, являющуюся частным случаем функции из семейства III.4.2.

Если же $l'$ -- константа, $l=l_1t+l_2$, то функция $w(c)$ находится из уравнения $[(2j+1)l_1 + ac(1-j)]w_c + 3aj w=0$ и имеет вид $w=w_1[(2j+1)l_1 + ac(1-j)]^{\frac {3j}{j-1}}$ при $j\ne 1$ и $w=w_1e^{-\frac{a}{l_1}c}$ при $j=1$. Соответственно мы получаем $\mathcal{F}$ вида III.3.2 и III.3.3.

Осталось рассмотреть случай $j=-2$. Тогда от (49) остается только одно уравнение
$$(c\tau'-l')w_c -2\tau' w = - l'' + c\tau'',$$
которое мы будем решать уже относительно $w(c)$. Если $\tau'(t)$ не является тождественным нулем, то, фиксируя соответствующее значение $t$, мы получаем, что $w$ является решением уравнения
$(ac+b)w_c-2aw=\alpha c+\beta$ и имеет вид $w=w_1(ac+b)^2-\frac {\alpha}ac-\frac{\beta-\alpha b}{2a}$, другими словами, представляет собой квадратный трехчлен. При $a\ne 0$ за счет сдвига можно избавиться от линейного члена, и мы пришли к функции $\mathcal{F}= u_1 f^{-2}+Qc^2+R$ II.3, аналогичная ситуация складывается при $a=0$, но $\alpha\ne 0$ если же $a=\alpha=0$, то мы получаем функцию $\mathcal{F}= u_1 f^{-2}+Qc+R$ семейства II.1.

Если же $\tau'=0$, то при $l'=0$ мы получаем двумерную алгебру, а при $l'\ne 0$ -- $w_c=k$, $l''=kl'$, и мы снова приходим к функции из семейства II.1.\\

Довольно больших рассуждений требует {\bf случай 5.1}: $j = -1$, тогда $u = u_1$ и $\tau_t + \xi_x = 0$ для $\tau(t,x)$, $\xi(t,x)$. Здесь уравнение для $w(c)$ приобретает вид
$$(\xi_t + 2c\xi_x - c^2\tau_x)w_c + (3c\tau_x -3\xi_x)w = \xi_{tt} + 3c\xi_{tx} + 3c^2\xi_{xx} - c^3\tau_{xx}. \eqno(50)$$
Это уравнение оказывается вырожденным (коэффициент при $w_c$ равен нулю) только если алгебра двумерна. Поэтому случай алгебры большей размерности может быть в одном из двух вариантов: либо $\xi_x=\tau_x=0$ (а тогда и $\tau_t=0$), но $\xi_t\ne 0$, либо $\xi_x$ и $\tau_x$ не являются нулевыми одновременно. В первом случае уравнение (50) приобретает вид $\xi_t w_c = \xi_{tt}$, переменные в нем разделяются, и мы получаем линейную функцию $w=w_1c+w_2$ (что дает нам $\mathcal{F}= u_1 f^{-1}+Qc+R$, то есть частный случай функции III.4.2) и трехмерную алгебру.

Если же по крайней мере одна из $\xi_x$ и $\tau_x$ не является тождественным нулем, то, зафиксировав соответствующие значения $t,x$, мы можем получить, что $w(c)$ является решением уравнения вида
$$(ac^2+2bc+d)w_c - 3(ac+b)w = \alpha c^3+\beta c^2+\gamma c+\delta. \eqno(51) $$

{\bf Лемма 3.} {\em Решения уравнения (51) имеют следующий вид, в зависимости от коэффициентов $a,b,d$:

1. если $b^2-ad\ne 0$, то
$$w(c)=w_1(ac^2+2bc+d)^{3/2}+pc^3+qc^2+rc+s,$$
где $p,q,r,s$ однозначно выражаются через $\alpha$, $\beta$, $\gamma$, $\delta$, а $w_1$ -- произвольная константа;

2.если $b^2-ad=0$, $a\ne 0$ (то есть $ac^2+2bc+d=a(c+h)^2$), то
$$w(c)=w_1(c+h)^3+p(c+h)^3\ln(c+h)+qc^2+rc+s,$$
где $p,q,r,s$ однозначно выражаются через $\alpha$, $\beta$, $\gamma$, $\delta$, а $w_1$ -- произвольная константа;

3. если $b^2-ad=0$ и $a=0$ (а тогда и $b=0$), но $d\ne 0$, то
$$w=w_1+\frac 1d \left[\frac 14\alpha c^4+\frac 13\beta c^3+\frac 12 \gamma c^2+\delta c\right],$$
где $w_1$ -- произвольная константа.}

{\bf Доказательство леммы 3} фактически сводится к решению неоднородного уравнения, поскольку решение однородного элементарно. Решение неоднородного уравнения осуществляется методом неопределенных коэффициентов: представив искомое решение в виде $w(c)=pc^3+qc^2+rc+s$, мы после подстановки и приравнивания коэффициентов при одинаковых степенях $c$ получаем систему для определения $p,q,r,s$
$$3bp-aq=\alpha,\qquad 3pd+bq-2ar=\beta,\qquad 2dq-br-3as=\gamma, \qquad dr-3bs=\delta.$$
Ее определитель равен $9(b^2-ad)^2$, и поэтому в случае 1 мы и получаем утверждение леммы.

Если $b^2-ad=0$, то квадратный трехчлен является полным квадратом, наше уравнение становится уравнением Эйлера с резонансным слагаемым в правой части:
$$a(c+h)^2w_c - 3a(c+h)w = \alpha (c+h)^3+\tilde\beta c^2+\tilde\gamma c+\tilde\delta,$$
решение которого имеет как раз вид, указанный в пункте 2.

Оставшийся же случай $a=b=0$ решается прямым интегрированием. Лемма доказана.
\\

Таким образом, нам остается подставить в уравнение указанные варианты решений, и осуществить расщепление по $c$. Для упрощения выкладок можно заметить, что заменой переменных $\bar t=t$, $\bar x=x+\lambda t$ с подходящим $\lambda$ мы можем избавиться от коэффициента $b$ (если $a\ne 0$) или от коэффициента $d$ (если $a=0$, $b\ne 0$), и использовать далее следующие анзацы:

1. $w=w_1(c^2+n)^{3/2}+pc^3+qc^2+rc+s$ ($n\ne 0$),

2. $w=w_1 c^{3/2}+pc^3+qc^2+rc+s$,

3. $w=w_1 c^3+pc^3\ln c+qc^2+rc+s$,

4. $w=w_1+pc^4+qc^3+rc^2+sc$.

Для первого анзаца получаем из (50)
$$(\xi_t + 2c\xi_x - c^2\tau_x)[3c w_1(c^2+n)^{1/2}+3pc^2+2qc+r] + (3c\tau_x -3\xi_x)[w_1(c^2+n)^{3/2}+pc^3+qc^2+rc+s] = $$
$$=\xi_{tt} + 3c\xi_{tx} + 3c^2\xi_{xx} - c^3\tau_{xx}.$$
Отделение иррациональности приводит к разделению на два уравнения:
$$3w_1[c(\xi_t+n\tau_x)  + c^2\xi_x -n \xi_x] =0,$$
$$(\xi_t + 2c\xi_x - c^2\tau_x)[3pc^2+2qc+r] + (3c\tau_x -3\xi_x)[pc^3+qc^2+rc+s] = \xi_{tt} + 3c\xi_{tx} + 3c^2\xi_{xx} - c^3\tau_{xx}.$$
Из первого следует, что при $w_1\ne 0$ получаем $\xi_x=0$ (а значит, и $\tau_t=0$), $\xi_t+n\tau_x=0$, что может быть только если $\tau=ax+b$, $\xi=-ant+l$. Но тогда второе уравнение приобретает вид
$$a[-(c^2+n)(3pc^2+2qc+r) + 3c(pc^3+qc^2+rc+s)] = 0,$$
и трехмерная алгебра (при $a\ne 0$) получается только если выполнены условия
$$q=0,\,\,\, r=\frac 32 pn,\,\,\, s=\frac 23 qn,\,\,\, r=0,$$
то есть только если $p=q=r=s=0$. В итоге мы получаем функцию
$$\mathcal{F}= u_1 f^{-1}+w_1(c^2+n)^{3/2}= (c^2+n)^{3/2}([u_1 [f(c^2+n)^{3/2}]^{-1}+w_1)$$ -- частный случай функции из семейства III.3.1.
Если же $w_1=0$, то остается многочлен третьей степени.

Для второго анзаца аналогичное выделение из уравнения
$$(\xi_t + 2c\xi_x - c^2\tau_x)[\frac 32 w_1 c^{1/2}+3pc^2+2qc+r] + (3c\tau_x -3\xi_x)[w_1 c^{3/2}+pc^3+qc^2+rc+s] =$$
$$=\xi_{tt} + 3c\xi_{tx} + 3c^2\xi_{xx} - c^3\tau_{xx}$$
выделение иррациональных слагаемых дает равенства
$$w_1[\frac 32\xi_t + \frac 32c^2\tau_x ]=0,$$
$$(\xi_t + 2c\xi_x - c^2\tau_x)[3pc^2+2qc+r] + (3c\tau_x -3\xi_x)[pc^3+qc^2+rc+s] = \xi_{tt} + 3c\xi_{tx} + 3c^2\xi_{xx} - c^3\tau_{xx},$$
и если $w_1\ne 0$, то $\xi_t=\tau_x=0$, а поскольку $\tau_t+\xi_x=0$, получаем $\tau=at+b$, $\xi=-ax+l$, и второе уравнение приводится к виду
$$-a(2c[3pc^2+2qc+r] -3[pc^3+qc^2+rc+s])=0,$$
откуда либо $a=0$ (то есть алгебра двумерна), либо $p=q=r=s=0$. Мы получаем $\mathcal{F}= c^{3/2}[u_1 (fc^{3/2})^{-1}+w_1]$, то есть частный случай функции III.3.2 (для $m=3/2$, $G(z)=Pz^{-1}+Q$). Если же $w_1=0$, то у нас опять остается лишь многочлен третьей степени.

Для третьего анзаца из
$$(\xi_t + 2c\xi_x - c^2\tau_x)[3pc^2\ln c+pc^2+3w_1c^2+2qc+r] + (3c\tau_x -3\xi_x)[pc^3\ln c+w_1c^3+qc^2+rc+s] = $$
$$=\xi_{tt} + 3c\xi_{tx} + 3c^2\xi_{xx} - c^3\tau_{xx},$$
получаем
$$p[3\xi_t + 3c\xi_x] = 0,$$
$$(\xi_t + 2c\xi_x - c^2\tau_x)[pc^2+3w_1c^2+2qc+r] + (3c\tau_x -3\xi_x)[w_1c^3+qc^2+rc+s] = \xi_{tt} + 3c\xi_{tx} + 3c^2\xi_{xx} - c^3\tau_{xx},$$
и из первого соотношения при $p\ne 0$ следует $\xi_t=\xi_x(=-\tau_t)=0$, $\xi={\rm const}$, $\tau=\tau(x)$. Подставляя во второе уравнение, получаем
$$\tau'[-pc^4+qc^3+2rc^2+3sc] = - c^3\tau'',$$
откуда $p\tau'=r\tau'=s\tau'=0$, $q\tau'=-\tau''$. Поскольку $p\ne 0$, получаем $\tau'=0$, то есть двумерную алгебру.
Если же $p=0$, то у нас и в этот раз остается многочлен третьей степени.

Наконец, для четвертого анзаца
$$(\xi_t + 2c\xi_x - c^2\tau_x)[4pc^3+3qc^2+2rc+s] + (3c\tau_x -3\xi_x)[pc^4+qc^3+rc^2+sc+w_1] =$$
$$=\xi_{tt} + 3c\xi_{tx} + 3c^2\xi_{xx} - c^3\tau_{xx}$$
при $p\ne 0$ мы сразу же получаем $\tau_x=0$, а затем $\xi_x(=-\tau_t)=0$, а потом и $\xi_t=0$, то есть тут возможна только двумерная алгебра. Если же $p=0$, то и в этом случае у нас остается многочлен третьей степени.

Таким образом, нам осталось рассмотреть вариант, когда $w(c)=pc^3+qc^2+rc+s$. В этом случае уравнение (50) приобретает вид
$$(\xi_t + 2c\xi_x - c^2\tau_x)[3pc^2+2qc+r] + (3c\tau_x -3\xi_x)[pc^3+qc^2+rc+s] = \xi_{tt} + 3c\xi_{tx} + 3c^2\xi_{xx} - c^3\tau_{xx},$$
и его расщепление по $c$ дает систему
$$3p\xi_x+q\tau_x=-\tau_{xx},\quad 3p\xi_t+q\xi_x+2r\tau_x=3\xi_{xx},\quad 2q\xi_t-r\xi_x+3s\tau_x=3\xi_{tx},\quad r\xi_t-3s\xi_x=\xi_{tt}.$$
Имея в виду, что $\tau_t=-\xi_x$, ее можно дополнить уравнениями для $\tau_{tt}$ и $\tau_{tx}$:
$$-3p\xi_t-q\xi_x-2r\tau_x=3\tau_{tx},\qquad -2q\xi_t+r\xi_x-3s\tau_x=3\tau_{tt}.$$
Мы получили полную линейную систему относительно функций $\tau$, $\xi$. Условия совместности этой системы имеют вид
$$(6pr-2q^2)\xi_t+(6qs-2r^2)\tau_x=0,$$
$$(2q^2-6rp)\xi_x+(qr-9ps)\tau_x=0,$$
$$(9sp-rq)\xi_t+(2r^2-6qs)\xi_x=0.$$
Это линейная однородная система вида
$$\left (\begin{array}{ccc}
0&\gamma&\-\beta\\
-\gamma&0&\alpha\\
\beta&-\alpha&0
\end{array}\right)
\left (\begin{array}{c}
\xi_x\\
\xi_t\\
\tau_x
\end{array}\right)=0,$$
где $\alpha=qr-9ps$, $\beta=2r^2-6qs$, $\gamma=6pr-2q^2$. Эта система вырожденная, ее ранг равен двум, если $\alpha^2+\beta^2+\gamma^2\ne 0$, и решение пропорционально вектору $(\alpha,\beta,\gamma)$.

Если $qr-9ps=2r^2-6qs=6pr-2q^2=0$, то при $p=0$ мы немедленно получаем $q=r=0$, $w={\rm const}$ и мы приходим к функции $\mathcal{F}= u_1 f^{-1}+w_1$ семейства I. Если же $p\ne 0$, то, выражая $r$ и $s$ через $p$ и $q$, мы найдем, что $w=p(c+\frac{q}{3p})^3$. Пользуясь заменой $\bar t=t$, $\bar x=x+\frac{q}{3p}t$, мы можем привести функцию $\mathcal{F}$ к виду $\mathcal{F}=Pf^{-1}+Qc^3$, а затем заменой $\bar t=x$, $\bar x=t$ превратить ее в $\mathcal{F}= -Pf^{-1}-Q$, то есть в функцию того же семейства I.

Если же $qr-9ps$, $2r^2-6qs$, $6pr-2q^2$ не все равны нулю, то, положив $\xi_x=-\tau_t=h(t,x)(qr-9ps)$, $\xi_t=h(t,x)(2r^2-6qs)$, $\tau_x=h(t,x)(6pr-2q^2)$, где $h(t,x)$ -- неизвестная функция, запишем условия согласования этих производных
$$h_t(t,x)(qr-9ps)=h_x(t,x)(2r^2-6qs),\quad h_t(t,x)(6pr-2q^2)=-h_x(t,x)(qr-9ps).$$
Поскольку из трех фигурирующих здесь коэффициентов по крайней мере один ненулевой, это означает, что $h(t,x)$, а значит, и $\tau(t,x)$, и $\xi(t,x)$ являются функциями от некоторой комбинации $\alpha t+\beta x$, причем в силу $\xi_x=-\tau_t$ их можно выразить через одну и ту же функцию: $\tau=\beta \sigma(\alpha t+\beta x)+A$, $\xi=-\alpha  \sigma(\alpha t+\beta x)+B$.
Остается подставить эти выражения в исходную систему
$$(-3p\alpha\beta+q\beta^2)\sigma'=-\beta^3 \sigma'',$$
$$(-3p\alpha^2-q\alpha\beta+2r\beta^2) \sigma'=-3\alpha\beta^2\sigma'',$$
$$(-2q\alpha^2+r\alpha\beta+3s\beta^2)\sigma'=-3\alpha^2\beta\sigma'',$$
$$(-r\alpha^2+3s\alpha\beta)\sigma' =-\alpha^3\sigma''.$$
и увидеть, что если $\beta=0$, то мы (считая, без ограничения общности, $\alpha=1$) получаем, что либо $\sigma'=0$ и группа двумерна, либо $p=q=0$ (и тогда, в силу предположения $r\ne 0$), и мы получаем функцию $\mathcal{F}=\frac Pf+rc+s$, которая заменой $\bar t=t$, $\bar x=x+\frac sr t$ приводится к частному случаю функции
III.4.2. А если $\beta\ne 0$, то мы опять же линейной заменой $\bar t=t$, $\bar x=x+\frac \alpha\beta t$ приводим нашу ситуацию к случаю $\alpha=0$, дающему соответственно $r=s=0$ (и $q\ne 0$), что дает нам функцию $\mathcal{F}=\frac Pf+pc^3+qc^2$, приводящейся к предыдущему случаю заменой $\bar t=x$, $\bar x=t$.\\

Наконец, рассмотрим последний случай -- {\bf случай 5.2}: $j=-1/2$, $u=u_1$, $\xi = \xi(t,x)$, $\tau = {\rm const}$. Подставляя в уравнение, получаем
$$(\xi_t + c\xi_x)w_c -\xi_x w = \xi_{tt} + 2c\xi_{tx} + c^2\xi_{xx}.$$

При $\xi_x = \xi_t = 0$ алгебра двумерна. При $\xi_x = 0$, $\xi_t \ne 0$ получим $w=w_1 +\lambda c$ и функцию семейства II.4 при $\lambda\ne 0$ и функцию семейства II.1 при $\lambda=0$.

При $\xi_x \not\equiv 0$ мы, фиксируя временно соответствующее значение $(t,x)$ получаем, что $w$ удовлетворяет уравнению вида $(ac+b)w_c -a w = pc^2+qc+r,$
решение которого имеет вид
$$w=w_1(ac+b)+\alpha c^2+\frac qa (ac+b)\ln(ac+b)+\beta,$$
где $w_1$ -- произвольная константа, а $\alpha$ и $\beta$ выражаются через $p,q,r$ однозначно. Подставляя этот анзац обратно в уравнение, получаем
$$(\xi_t + c\xi_x)(w_1a+2\alpha c+q\ln(ac+b) + q) -\xi_x (w_1(ac+b)+\alpha c^2+\frac qa (ac+b)\ln(ac+b)+\beta) = \xi_{tt} + 2c\xi_{tx} + c^2\xi_{xx}.$$
Выделение логарифмических членов и расщепление оставшегося уравнения по степеням $c$ дает уравнения
$$q[a\xi_t-b\xi_x] =0,\qquad
\alpha\xi_x=\xi_{xx},\qquad
q\xi_x +2\alpha \xi_t = 2\xi_{tx},\qquad
\xi_t(w_1a+q)-\xi_x(bw_1b+\beta)=\xi_{tt}.$$
Из первого уравнения следует, что либо $q=0$, и тогда $w$ -- квадратный трехчлен, и наша функция сводится к функции семейства II.4, если этот квадратный трехчлен не константа, и к функции семейства II.1 в противном случае. Либо $q\ne 0$, тогда $\xi(t,x)=\xi(ax+bt)$, и остальные три уравнения приобретают вид
$$\alpha a\xi'=a^2\xi'',\quad (aq +2b\alpha) \xi' = 2ab\xi'', \quad [b(w_1a+q)-a(bw_1b+\beta)]\xi'=b^2\xi''.$$
Вычитая из второго уравнения первое, умноженное на $b/a$, получаем $aq\xi'=0$, что, в силу предположений $q\ne 0$, $a\ne 0$ влечет $\xi'=0$, то есть алгебра оказывается двумерной.\\

Таким образом, для функций вида $\mathcal{F}=\mathcal{F}(c,f)$ классификация полностью обоснована: уравнения со всеми функции такого вида, кроме функций семейств I-III, имеют только двумерную алгебру симметрий $\partial_t$, $\partial_x$.

\subsection{Алгебры симметрий уравнений с функцией вида $\mathcal{F} = T(c,tf)/t$}
Здесь, как и в предыдущем параграфе, мы доказываем двумерность алгебры симметрий уравнения (5), за исключением уравнений с функциями $\mathcal{F}$ из семейств I-III.
\subsubsection{Анзацы и классифицирующие уравнения}
Для этого семейства группа эквивалентности содержит, помимо сдвигов по переменной $x$ и растяжений, являющихся симметриями, еще два типа преобразований -- растяжение только по $x$ и замены $\bar t=t$, $\bar x=x-kt$, которые не изменяют $\mathcal{F}$ и  $f$, но позволяют осуществлять сдвиги и растяжения переменной $c$, чем мы будем в дальнейшем пользоваться.

Этот случай нам придется разбирать по той же схеме, что и случай функции $\mathcal{F}=\mathcal{F}(c,f)$. Для удобства дальнейших вычислений обозначим второй аргумент функции $T$ через $g$: $tf = g$. Тогда $\mathcal{F} = T(c,g)/t$, и уравнение (27) превращается в уравнение относительно $T$:
$$[\tau + t(3c\tau_x - 2\xi_x + \tau_t)]gT_g + t(\xi_t + c(\xi_x - \tau_t) - c^2\tau_x)T_c - [\tau - t(3c\tau_x -\xi_x + 2\tau_t)]T =$$
$$= t^2[\xi_{tt} + c(2\xi_{tx} - \tau_{tt}) + c^2(\xi_{xx} - 2\tau_{tx}) - c^3\tau_{xx}].\eqno(52)$$

Алгебра симметрий для $\mathcal{F} = T(c,g)/t$ содержит двумерную некоммутативную алгебру $\partial_x$, $t\partial_t+x\partial_x$, и для завершения рассуждений с алгебрами, нам необходимо показать, что в случаях, отличных от рассмотренных в предыдущем параграфе, алгебра симметрий имеет размерность, в точности равную двум.

Для этого предположим, что уравнение (52) выполнено для некоторых $(\tau(t,x),\xi(t,x))$, отличных от представленных в этой алгебре. Подставим их в уравнение, которое теперь будем рассматривать как уравнение относительно $T$. Здесь возможны несколько вариантов.

1) Если $\xi_t + c(\xi_x - \tau_t) - c^2\tau_x\equiv 0$, то $\tau_x=\xi_t=0$, а $\tau_t=\xi_x=a={\rm const}$. Но тогда $\tau = at+b$, $\xi=ax+l$ и уравнение упрощается до $b(gT_g - T) = 0$. Если $b = 0$, то мы получаем двумерную алгебру для произвольной $\mathcal{F} = T(c,tf)/t$ -- это общий случай семейства IV.2. Если же $b\ne 0$, то $T(c,g)=G(c)g$, а значит, $\mathcal{F} = M(c)f$, это семейство III.2.3.

2) Пусть теперь $\xi_t + c(\xi_x - \tau_t) - c^2\tau_x\not\equiv 0$. Тогда мы можем, зафиксировав соответствующее значение $(t,x)=(t^*,x^*)$ и, пометив звездочкой соответствующие значения функций, фигурирующих в уравнении, получить уравнение в переменных $(c,g)$ относительно функции $T(c,g)$. Решение этого уравнения относительно $T$ имеет следующий вид:
$$T = e^{I_2}\left(G(ge^{-I_1}) - I_3\right), \ I_1 =\!\!\int \frac{\tau^* + t^*(3c\tau^*_x - 2\xi^*_x + \tau^*_t)dc}{t^*(c^2\tau^*_x - c(\xi^*_x - \tau^*_t) - \xi^*_t)}, \,\,\, I_2 =\!\!\int \frac{\tau^* - t^*(3c\tau^*_x - \xi^*_x + 2\tau^*_t)dc}{t^*(c^2\tau^*_x - c(\xi^*_x - \tau^*_t) - \xi^*_t)}, $$
$$I_3 = \int \frac{t^*(\xi^*_{tt} + c(2\xi^*_{tx} - \tau^*_{tt}) + c^2(\xi^*_{xx} - 2\tau^*_{tx}) - c^3\tau^*_{xx})}{c^2\tau^*_x - c(\xi^*_x - \tau^*_t) - \xi^*_t} e^{-I_2}dc,$$
где $G(\cdot)$ -- произвольная функция. Таким образом, в этом случае мы получили следующий анзац для $T$:
$$T = u(c)G\left(gv(c)\right) + w(c). \eqno(53)$$

Подставим (53) в (52), обозначив $z = gv(c)$:
$$[\frac{\tau}{t} + (3c\tau_x - 2\xi_x + \tau_t) + \frac{v_c}v(\xi_t + c(\xi_x - \tau_t) - c^2\tau_x)]uzG_z -$$
$$- [\frac{\tau}t- (3c\tau_x -\xi_x + 2\tau_t) - \frac{u_c}u(\xi_t + c(\xi_x - \tau_t) - c^2\tau_x)]uG =$$
$$= t[\xi_{tt} + c(2\xi_{tx} - \tau_{tt}) + c^2(\xi_{xx} - 2\tau_{tx}) - c^3\tau_{xx}] - w_c(\xi_t + c(\xi_x - \tau_t) - c^2\tau_x) + w(\frac{\tau}t - (3c\tau_x -\xi_x + 2\tau_t)).$$

Здесь возможны три варианта:

I. Если коэффициент при $zG'(z)$ нулевой: $\tau/t + 3c\tau_x - 2\xi_x + \tau_t + (\xi_t + c(\xi_x - \tau_t) - c^2\tau_x)(\ln v)_c = 0$, то либо коэффициент при $G$ ненулевой и $G = const$, а значит $F$ не зависит от $f$ (что противоречит предположениям теоремы). Либо коэффициент при $G$ также нулевой, а значит нулевой и свободный член, в (53) функция $G(\cdot)$ произвольная, а $u$, $v$ и $w$ удовлетворяют системе уравнений
$$\frac{\tau}{t} - (3c\tau_x -\xi_x + 2\tau_t ) - (\xi_t + c(\xi_x - \tau_t) - c^2\tau_x)(\ln |u|)_c = 0,\eqno(54)$$
$$\frac{\tau}{t} + 3c\tau_x - 2\xi_x + \tau_t + (\xi_t + c(\xi_x - \tau_t) - c^2\tau_x)(\ln |v|)_c = 0,\eqno(55)$$
$$t[\xi_{tt} + c(2\xi_{tx} - \tau_{tt}) + c^2(\xi_{xx} - 2\tau_{tx}) - c^3\tau_{xx}] = (\xi_t + c(\xi_x - \tau_t) - c^2\tau_x)w_c - (\frac{\tau}{t} -(3c\tau_x -\xi_x + 2\tau_t))w.\eqno(56)$$

II. Если коэффициент при $zG'(z)$ ненулевой: $\tau/t + 3c\tau_x - 2\xi_x + \tau_t - (\xi_t + c(\xi_x - \tau_t) - c^2\tau_x)(\ln v)_c \ne 0$, а коэффициент при $G$ нулевой, то $G = G_1\ln |z|$ и $T = u(c)\ln |g| + w(c)$, где
$$\frac{\tau}{t} - (3c\tau_x -\xi_x + 2\tau_t ) - (\xi_t + c(\xi_x - \tau_t) - c^2\tau_x)(\ln |u|)_c = 0,\eqno(57)$$
$$[\frac{\tau}{t} + (3c\tau_x - 2\xi_x + \tau_t)]u + (\xi_t + c(\xi_x - \tau_t) - c^2\tau_x)w_c - (\frac{\tau}t - (3c\tau_x -\xi_x + 2\tau_t))w =$$
$$= t(\xi_{tt} + c(2\xi_{tx} - \tau_{tt}) + c^2(\xi_{xx} - 2\tau_{tx}) - c^3\tau_{xx}).\eqno(58)$$

III. Коэффициенты при $G$ и $G_z$ в уравнении ненулевые, тогда $G = G_1z^{j} + G_2$ и функция $T$ имеет вид $T = u(c)g^j + w(c)$
$$j(\frac{\tau}{t} + (3c\tau_x - 2\xi_x + \tau_t)) - \frac{\tau}t + (3c\tau_x -\xi_x + 2\tau_t) + (\xi_t + c(\xi_x - \tau_t) - c^2\tau_x)(\ln |u|)_c = 0,\eqno(59)$$
$$w_c(\xi_t + c(\xi_x - \tau_t) - c^2\tau_x) -  w(\frac{\tau}t - (3c\tau_x -\xi_x + 2\tau_t)) = t(\xi_{tt} + c(2\xi_{tx} - \tau_{tt}) + c^2(\xi_{xx} - 2\tau_{tx}) - c^3\tau_{xx}).\eqno(60)$$

Как видно, и в этом семействе выделяются {\em общий, логарифмический и степенной} анзацы.

\subsubsection{Первичная классификация}

Как и в предыдущем случае, уравнение относительно $u(c)$
$$j(\frac{\tau}{t} + (3c\tau_x - 2\xi_x + \tau_t)) - \frac{\tau}t + (3c\tau_x -\xi_x + 2\tau_t) + (\xi_t + c(\xi_x - \tau_t) - c^2\tau_x)(\ln u)_c = 0 \eqno(59)$$
-- общее для всех трех анзацев (для анзацев I и II $j=0$), поэтому естественно за основание классификации взять именно его.

{\bf Лемма 4}. {\em Уравнение (59) как уравнение относительно $\tau$, $\xi$ имеет решения, отличные от $\tau=at$, $\xi=ax+l$, только для функций $u(c)$, приведенных в таблице 8.
\begin{figure}[h]
    \centering
  Таблица 8. Первичная классификация: отличные от констант решения $\tau=at$, $\xi=ax+l$ уравнения (59).
\begin{tabular}{|c||c|c|c|}
\hline
№ & $j$ & $u$& $\tau$, $\xi$ \\
\hline
\hline
1 &$j=1$ & $u$ произвольная & $\tau=at+b$, $\xi=ax+l$\\
\hline
2 &$j\ne 1$ & $u=u_1(c^2+n)^{3(j+1)/2}$ & $\tau=2atx+bt$, \\
&&&$\xi=a(x^2-nt^2)+bx+l$\\
\hline
2.1&$j=1$ & $u=u_1(c^2+n)^3$ & $\tau_t=\xi_x$, $\xi_t+n\tau_x=0$\\
\hline
3 & $j\ne 1$ & $u = u_1c^m$ & $\tau = at + bt^{\frac{j-1}{m-(j+2)}}$, $\xi = ax + l$\\
& & $m \ne 3(j+1), j+2, 2j+1, 0$ & \\
\hline
3.1& $j \ne 1, -2, -1/2$ & $u = u_1c^{2j+1}$& $\tau = at$, $\xi = \xi(x)$\\
\hline
3.2& $j \ne -1, -2, -1/2$ & $u=u_1c^{3(j+1)}$ & $\xi=\xi(x)$, $\tau = \xi_xt + b(x)t^{\frac{j-1}{2j+1}}$ \\
\hline
3.3& $j = 1$ & $u = u_1c^3$& $\tau = \tau(t)$, $\xi = \xi(x)$\\
\hline
3.4& $j = -1/2$ & $u=u_1c^{3/2}$& $\tau = t\xi'(x)$, $\xi = \xi(x)$\\
\hline
3.5& $j = -2$ & $u = u_1c^{-3}$& $\tau = ta(x)$, $\xi = \xi(x)$\\
\hline
4& $j\ne 1$ & $u = u_1e^{mc}$ & $\tau = at + mb$,  \\
&$m\ne 0$ & &$\xi = ax - (j-1)b\ln t  + l$ \\
\hline
5 & $j \ne -1, -1/2$ & $u = u_1$ & $\tau=\tau(t)$,\\
&&&$\xi = (\frac{j-1}{2j+1}\frac \tau t + \frac{j+2}{2j+1}\tau_t)x + l(t)$\\
\hline
5.1 &  $j = -1$ & $u = u_1$ & $\tau=\tau(t,x)$, $\xi=\xi(t,x)$,  \\
& & &$t(\xi_x + \tau_t) = 2\tau$ \\
\hline
5.2 &  $j = -1/2$ & $u = u_1$ & $\tau=at$, $\xi=\xi(t,x)$.\\
\hline
\end{tabular}
\end{figure}}

{\bf Доказательство.}
Будем рассматривать (59) как уравнение на функцию $u$. Если это уравнение выполнено тождественно (то есть и коэффициент при $(\ln u)_c$, и свободный член равны нулю), то $\tau=at+b$, $\xi=ax+l$ и $(j-1)b=0$. Если $b=0$, то алгебра двумерна, а если $j=1$ -- это случай 1.

Если же уравнение (59) невырождено, то, так как $u$ не зависит от $t$, $x$, то можно зафиксировать произвольным образом значения $t$ и $x$, при которых коэффициент при $(\ln u)_c$ ненулевой, а также соответствующие значения $\tau$, $\xi$ и их производных, представив уравнение (59) в виде
$$(\ln u)_c = \frac{(j-1)p + 3k(j+1)c - r(2j+1) + q(j +2)}{kc^2 - (r - q)c - s},\eqno(61)$$
где $k$, $p$, $q$, $r$, $s$ -- некоторые константы, причем выбор $(t,x)$ можно осуществить так, что знаменатель будет не равен нулю тождественно (противный случай рассмотрен в предыдущем пункте). Здесь, как и в предыдущем случае, возможны три варианта:\\

А) Если $k \ne 0$, и дробь в (61) несократима (случай сократимой дроби, как и случай $k=0$ будет рассматриваться в пункте Б) ) и не равна нулю (этот случай будет рассматриваться в пункте В)), то заменой переменных $\bar{t} = t$, $\bar{x} = x - (r-q)t/(2k)$ мы можем упростить формулу для $u$:
$$(\ln u)_c = \frac{3(j+1)c + m}{c^2 + n}.\eqno(62)$$

Б) Если в (61) $k = 0$, но $r-q\ne 0$, то заменой переменных $\bar{t} = t$, $\bar{x} = (r-q)x + nt$ мы можем упростить (61), приведя его к виду
$$(\ln u)_c = {m}/{c},\eqno(63)$$
к этому же виду приводится (61) и в случае $k\ne 0$, но дробь сократима. Мы здесь будем предполагать, что $m\ne 0$.

B) Если $k=r-q=0$, то дробь в (61) сводится к константе $$(\ln u)_c = m.\eqno(64)$$

А) Случай {\bf функции (62)}. Предполагая, что (62) -- несократимая дробь, подставим (62) в (59), получим:
$$(c^2 + n)[j(\frac{\tau}{t} + (3c\tau_x - 2\xi_x + \tau_t)) - \frac{\tau}t + (3c\tau_x -\xi_x + 2\tau_t)] + (\xi_t + c(\xi_x - \tau_t) - c^2\tau_x)(3(j+1)c + m) = 0.$$

Приравняем к нулю коэффициенты при различных степенях $c$:
$$(j-1)\frac{\tau}{t} - (2j+1)\tau_t - m\tau_x + (j+2)\xi_x = 0,$$
$$- m\tau_t + 3n(j+1)\tau_x + 3(j+1)\xi_t + m\xi_x = 0,$$
$$n(j-1)\frac{\tau}{t} + n(j+2)\tau_t + m\xi_t - n(2j+1)\xi_x = 0.$$

Мы получили систему вида
$$a_{0i}\frac\tau t + a_{1i}\tau_t+a_{2i}\tau_x+a_{3i}\xi_t+a_{4i}\xi_x=0, \qquad i=1,2,3,$$
которой соответствует матрица
$$\left(\begin{array}{ccccc}
(j-1) &-(2j+1)&-m&0&j+2\\
0 &-m&3n(j+1)&3(j+1)&m\\
n(j-1)&n(j+2)&0&m&-n(2j+1)\\
\end{array}\right)$$
Вычтем из третьей строки первую, умноженную на $n$:
$$\left(\begin{array}{ccccc}
(j-1) &-(2j+1)&-m&0&j+2\\
0 &-m&3n(j+1)&3(j+1)&m\\
0&3n(j+1)&nm&m&-3n(j+1)\\
\end{array}\right)$$
Основные варианты связаны с рангом этой матрицы. Этот ранг не может быть нулевым (так как равенства $j+2=0$, $j-1=0$ одновременно не могут быть выполнены). Ранг, равный 1, возможен только при $m = 0$, $j=-1$ (что противоречит предположению о том, что дробь в (62) ненулевая).

Ранг равный двум будет в случае равенства нулю всех миноров третьего порядка.
При $j\ne 1,-1$ для ранга два необходимо $n = -m^2/9(j+1)^2$, но в этом случае дробь снова оказывается сократимой. При $j=-1$ ранг, равный двум, невозможен (если $m=0$, то ранг равен единице, если нет -- то трем). При $j=1$ ранг, равный двум, возможен либо при $m=0$, что дает $u(c)=(c^2+n)^{3}$, $\xi_x=\tau_t$, $\xi_t+n\tau_x=0$ -- это случай 2.1,  либо при $m=36n^2$, что вновь противоречит несократимости дроби в (62).

Если ранг равен трем, то имеем три независимых уравнения. Здесь есть два варианта.
Если $m = 0$, то $n\tau_x+\xi_t= 0$, $\tau_t=\xi_x$ и
$(j-1)\frac{\tau}{t} - (j-1)\tau_t= 0$, где $j\ne 1$ (иначе ранг системы будет равен двум).

Если же $m\ne 0$, то комбинированием уравнений наша матрица приводится к виду
$$\left(\begin{array}{ccccc}
(j-1) &-(2j+1)&-m&0&j+2\\
0 &-1&0&0&1\\
0&0&n&1&0\\
\end{array}\right)$$
(при этом исключается случай $m^2 + 9n(j+1)^2=0$, поскольку он соответствует сократимой или нулевой дроби в (62)). В результате мы получаем систему уравнений
$$\tau_t = \xi_x,\qquad n\tau_x+\xi_t= 0,\qquad t[(j-1)\tau_t + m\tau_x] = (j-1)\tau.$$
Нетрудно видеть, что эти формулы при $m=0$ включают и предыдущий случай, при этом ранг системы равен трем тогда и только тогда, когда $m$ и $j-1$ не обращаются одновременно в нуль.

Решение третьего уравнения этой системы имеет вид $\tau = th((j-1)x-qt)$, где  $h(\cdot)$ -- некоторая функция. Подставляя полученную $\tau$ в условие совместности первого и второго $\tau_{tt}+n\tau_{xx}=0$, получаем
$$tm^2h''(s)-2mh'(s)+tn(j-1)^2h''(s)=0,$$
где $s=(j-1)x-mt$.

Если $j\ne 1$, аргумент функции $h$ является независимым с $t$, и поэтому можно осуществить расщепление по $t$, что дает $mh'=(m^2+n(j-1)^2)h''=0$.

Если $m\ne 0$, то $h'=0$, и мы получаем двумерную алгебру $\tau=at$, $\xi=ax+l$. Если же $m=0$, то при $n=0$ мы получаем сократимую дробь в (62), что противоречит предположению, а при $n\ne 0$ оказывается $h''=0$, и мы приходим к трехмерной алгебре $\tau=2atx+bt$, $\xi=a(x^2-t^2)+bx+l$ для $u=(c^2+n)^{3(j+1)/2}$, что и дает нам случай 2.

Наконец, при $j=1$ мы получаем, в силу $m\ne 0$, что $\tau_x=\xi_t=0$, $\tau=at+b$, $\xi=ax+l$ -- это частный вариант случая 1.
\\

Б) Случай {\bf функции (63)}. Подставим (63) в (59), предполагая $m\ne 0$:
$$(j-1)c\frac{\tau}{t} + 3(j+1)c^2\tau_x - (2j+1)c\xi_x + (j+2)c\tau_t + m(\xi_t + c(\xi_x - \tau_t) - c^2\tau_x) = 0.$$

Получаем систему уравнений:
$$[3(j+1) - m]\tau_x = 0,\qquad (j-1)\frac{\tau}{t} +[m - (2j+1)]\xi_x + [(j+2)-m]\tau_t = 0,\qquad m\xi_t = 0.$$
Так как $m\ne 0$, то из последнего уравнения $\xi= \xi(x)$. Дальнейшее ветвление вариантов определяется соотношениями между $m$ и $j$, точнее, равно или не равно $m$ одному из чисел $3(j+1)$, $j+2$ и $2j+1$. Попарное совпадение этих чисел происходит при $j=1,-2,-1/2$, поэтому случай этих $j$ мы рассмотрим отдельно.

a) Пусть $j\ne 1,-2,-1/2$.

Если $m\ne 3(j+1)$ то из первого уравнения $\tau = \tau(t)$. А тогда из второго уравнения получаем следуюшие варианты.

Либо $m\ne 2j + 1$ и $m\ne j+2$, тогда переменные разделяются, $\xi_x={\rm const}$, а значит, $\xi = ax+l$, и $\tau = at+bt^{\frac{j-1}{m - (j+2)}}$. Это  случай 3.

Либо $m=2j+1$, тогда получаем $\xi=\xi(x)$, $\tau=at$. Это случай 3.1.

Если $m=j+2$, то $\xi=\xi(x)$, $\tau=t\xi_x$, а поскольку $\tau$ зависит только от $t$, то мы получаем $\xi=ax+b$, $\tau=at$, то есть только двумерную алгебру.

Если $m = 3(j+1)$, то первое уравнение выполняется автоматически, а второе приобретает вид $(j-1)\tau + [j+2]t\xi_x - [2j+1]t\tau_t = 0$. Тогда, поскольку $j \ne -2, -1/2$, мы получаем алгебру $\tau = p(x)t^{\frac{j-1}{2j+1}} + \xi_xt$, $\xi = \xi(x)$ -- это случай 3.2.

б) Пусть $j=1$, тогда по-прежнему $\xi_t=0$, а остальные уравнения приобретают вид
$$[6 - m]\tau_x = 0,\qquad [m - 3]\xi_x + [3-m]\tau_t = 0,$$
что при $m\ne 3,6$ дает трехмерную алгебру $\tau=at+b$, $\xi=ax+l$. Это частный вариант случая 3. При $m=3$ алгебра имеет вид $\tau=\tau(t)$, $\xi=\xi(x)$ -- случай 3.3, при $m=6$ -- $\xi=\xi(x)$, $\tau=t\xi'(x)+p(x)$, это частный вариант случая 3.2.

в) При $j = -1/2$ все также $\xi_t = 0$, а остальные уравнения приобретают вид
$$[\frac 32 - m]\tau_x = 0,\qquad -\frac 32 \frac{\tau}{t} + m\xi_x + [\frac 32 -m]\tau_t = 0.$$
При $m\ne \frac 32 $ алгебра имеет вид $\xi=ax+l$, $\tau=at+bt^{3/(3-2m)}$ -- это частная версия случая 3, при $m=\frac 32$ -- $\xi=\xi(x)$, $\tau=t\xi'(x)+p(x)$ -- это случай 3.4.

г) При $j = -2$ имеем $\xi_t = 0$
$$-[3+m]\tau_x = 0,\qquad -3\frac{\tau}{t} +[m+3]\xi_x -m\tau_t = 0.$$
При $m\ne -3$ мы получаем $\tau=at+bt^{-3/m}$, $\xi=ax+l$ -- это частный вариант случая 3, а при $m = -3$ алгебра имеет вид $\tau = a(x)t$, $\xi= \xi(x)$, это случай 3.5.
\\

В) Случай {\bf функции (64)}. Подставим в (59) функцию (64), получим систему:
$$m\tau_x = 0,\qquad 3(j+1)\tau_x + m\xi_x - m\tau_t = 0,\qquad
(j-1)\frac{\tau}{t} - (2j+1)\xi_x + (j+2)\tau_t + m\xi_t = 0.$$

a) Если $m\ne 0$, то $\tau_x = 0$ и $\xi_x = \tau_t$, $(j-1)\frac{\tau}{t} - (j-1)\tau_t + m\xi_t = 0$. Из первого соотношения $\xi = \tau_tx + b(t)$, тогда из второго уравнения $\tau = at + p$ и $b_t = -(j-1)\frac{p}{mt}$. Таким образом, в этом случае алгебра имеет вид $\tau = at + p$, $\xi = ax - (j-1)\frac{p}{m}\ln t + l$, это случай 4.

б) Если $m=0$, то $(j+1)\tau_x = 0$, $(j-1)\frac{\tau}{t} - (2j+1)\xi_x + (j+2)\tau_t = 0$.

При $j\ne -1,-1/2$, получаем $\tau_x = 0$, $(j-1)\frac{\tau}{t} - (2j+1)\xi_x + (j+2)\tau_t = 0$, откуда находим $\tau = \tau(t)$, $\xi = (\frac{j-1}{2j+1}\frac \tau t + \frac{j+2}{2j+1}\tau_t)x + b(t)$. Это случай 5.

При $j=-1$ получаем, что первое уравнение выполняется автоматически, а второе приобретает вид $t(\xi_x + \tau_t) = 2\tau$ для $\tau(t,x)$, $\xi(t,x)$, это случай 5.1.

Наконец, при $j = -1/2$ получаем $\tau_x=0$, $\tau = t\tau_t$, откуда $\tau = at$, $\xi = \xi(t,x)$. Это случай 5.2.

Лемма доказана.

\subsubsection{Анзац I. $T = u(c)G\left(gv(c)\right) + w(c)$, $G(\cdot)$ -- произвольная функция}
Функции $v(c)$ и $w(c)$ удовлетворяют уравнениям (55) и (56):
$$\frac{\tau}{t} + 3c\tau_x - 2\xi_x + \tau_t + (\xi_t + c(\xi_x - \tau_t) - c^2\tau_x)(\ln |v|)_c = 0,\eqno(55)$$
$$t(\xi_{tt} + c(2\xi_{tx} - \tau_{tt}) + c^2(\xi_{xx} - 2\tau_{tx}) - c^3\tau_{xx}) = (\xi_t + c(\xi_x - \tau_t) - c^2\tau_x)w_c - (\frac{\tau}{t} -(3c\tau_x -\xi_x + 2\tau_t))w.\eqno(56)$$
Подставим в эти уравнения полученные нами в лемме 4 выражения для $\tau$ и $\xi$ в каждом из случаев, отвечающих $j=0$. Это случаи 2, 3, 3.1, 3.2, 4, 5.

Начнем со {\bf случая 2}. Подставим $\tau=2atx+bt$, $\xi=a(x^2-nt^2)+bx+l$ в уравнения для $v(c)$ и $w(c)$. Получим
$$a[3c-(n+c^2)(\ln |v|)_c] = 0,\qquad -2a(n + c^2) = -2a(n + c^2)w_c +6ca w.$$
Случай $a=0$ приводит нас к двумерной алгебре, если же $a\ne 0$, то $v=v_1(c^2+n)^{3/2}$, $w=w_1(c^2+n)^{3/2}+\frac {c(c^2+n)}n$ при $n\ne 0$, и $w=w_1c^3-\frac 12 c$ при $n=0$, что дает нам
$$ T = (c^2+n)^{\frac 32}G\left(g(c^2+n)^{\frac 32}\right) + \frac{c(c^2+n)}{n}, \quad n\ne 0,$$
$$T = c^3G\left(gc^3\right) - \frac{c}{2}, \quad n= 0,$$
то есть функции $\mathcal{F}=\frac 1t T(c,tf)$ вида III.5 и III.2.2 соответственно\\

{\bf Случай 3}. Здесь подстановка $\tau = at+bt^{\frac{1}{2 - m}}$, $\xi = ax+l$ дает
$$b[(3 - m) - c(\ln |v|)_c] = 0, \quad b[c(m - 1) + (2-m)(- cw_c + m w)] = 0.$$
Если $b=0$, то получаем двумерную алгебру. Если $b\ne 0$, то $v = v_1c^{3-m}$, $w = w_1c^m + \frac 1{m-2}c$, и мы получаем функцию $\mathcal{F} = \frac {c^m}tG(tc^{3-m}f) + \frac c{t(m-2)}$, из семейства III.3.2. \\

{\bf Случай 3.1}. $\tau = at$, $\xi = \xi(x)$. Подстановка в уравнения для $v$ и $w$ дает
$$(\xi_x - a)(c(\ln |v|)_c -2) = 0, \qquad tc^2\xi_{xx} = (\xi_x - a)(cw_c - w).$$

Из второго уравнения следует, что $\xi(x)$ -- линейная функция. Если $\xi_x=a$, то мы получаем двумерную алгебру. Если же $\xi_x\ne a$, то $v=v_1c^2$, $w=w_1c$, что приводит нас к $T = cG(gc^2)$ (как видно, $w(c)$ поглощается главной частью), и мы получаем функцию $\mathcal{F}=\frac{c}{t}G(c^2tf)$ из семейства III.3.2.\\

{\bf Случай 3.2}. Подстановка (напомним, при $j=0$) $\tau = b(x)t^{-1} + \xi_xt$, $\xi = \xi(x)$
в уравнение (55) для $v$ и расщепление его по переменной $t$ дает три равенства:
$$b(x)(\ln |v|)_c = 0, \qquad b'(x)[3 - c(\ln |v|)_c] = 0,\qquad \xi_{xx} [3 - c(\ln |v|)_c] = 0.$$
Если $b(x)\not\equiv 0$, то $v={\rm const}$, тогда $b'=\xi''=0$, и подстановка этих функций в уравнение (56) для $w$ дает $w=w_1c^3+c$. Мы получили функцию $\mathcal{F}=\frac {c^3}tG(tf)+\frac ct$, из семейства III.3.2 (при $m=3$).

Если $b(x)=0$ и $\xi''=0$, то алгебра двумерна. Если же $b(x)=0$, а $\xi''(x)\ne 0$, то $v=v_1c^3$, а подстановка в уравнение для $w$ приводит к равенству
$$-c^2\xi''- c^3t\xi'' = - c^2\xi''w_c +3c\xi''w,$$
откуда $\xi'''(x)=0$, $cw_c =3w+c$,
$w=w_1c^3-\frac c2$, и мы получаем функцию $\mathcal{F}=\frac {c^3}tG(c^3tf)-\frac c{2t}$ из семейства III.2.2.

{\bf Случай 4}. Подстановка $\tau = at + bm$, $\xi = ax + b\ln t + l$ в уравнения (55)-(56) для $v$ и $w$ дает
$$b[m + (\ln |v|)_c] = 0,\qquad b[1 + w_c - mw] = 0.$$

Если $b=0$, то мы получаем двумерную алгебру. Если же $b\ne 0$, то $v = v_1e^{-mc}$, $w = w_1e^{mc} + \frac 1{m}$, $T = e^{mc}G(ge^{-mc}) +  \frac 1{m}$, $\mathcal{F}=\frac 1t e^{mc}G(te^{-mc}f) +  \frac 1{mt}$ из семейства III.2.3.

{\bf Случай 5}. И тут, подставив $\tau = \tau(t)$, $\xi = (2\tau_t - \frac \tau t)x + l(t)$ в уравнения (55)-(56), получим
после расщепления по $x$
$$(2\tau_{tt} - \frac {\tau_t}t + \frac \tau{t^2})(\ln |v|)_c = 0,\quad
(l_{t} + c(\tau_t - \frac \tau t))(\ln |v|)_c = 3(\tau_t - \frac{\tau}{t}),$$
$$2t\tau_{ttt} - \tau_{tt} + 2\frac {\tau_t}{t} - \frac {2\tau}{t^2} = (2\tau_{tt} - \frac {\tau_t}t + \frac \tau{t^2})w_c,$$
$$tl_{tt} + c(3t\tau_{tt} - 2\tau_t + 2\frac \tau{t}) = (l_{t} + c(\tau_t - \frac \tau t))w_c.$$

Из первого уравнения есть два варианта: равен нулю либо первый, либо второй множитель.

а) $(\ln v)_c=0$, $v={\rm const}$, тогда из второго уравнения $\tau=at$, $\xi=ax+l(t)$. Третье уравнение превращается в тождество, а четвертое приобретает вид $tl_{tt}=l_tw_c$. Если $l_t=0$, то мы получаем двумерную алгебру. Если же $l_t\ne0$, то $tl_{tt}=\lambda l_t$, $l(t)=kt^{\lambda+1}+h$, $w=\lambda c$, и мы приходим к функции $\mathcal{F}=\frac 1tG(tf)+\lambda\frac ct$ из семейства III.4.2.

В исключительном случае $\lambda=-1$ будет $l(t)=k\ln t+h$, а функция
$\mathcal{F}=\frac 1tG(tf)-\frac ct$ принадлежит уже семейству III.4.3.

б) $(\ln v)_c\ne 0$, тогда $2\tau_{tt} - \frac {\tau_t}t + \frac \tau{t^2}=0$, $\tau=at+bt^{1/2}$. При этом второе уравнение приобретает вид
$(bc - 2t^{1/2}l_{t})(\ln |v|)_c = 3b$, откуда, в силу предположения $(\ln v)_c\ne 0$, мы получаем, что $2t^{1/2}l_{t}=\lambda$, $l(t)=\lambda t^{1/2}+h$. Если $\lambda=b=0$, то мы получаем двумерную алгебру, иначе $v=(bc-\lambda)^3$. Третье уравнение оказывается опять тождеством, а четвертое приобретает вид $\frac 14(bc-\lambda)=\frac 12(\lambda-bc)w_c$, откуда снова либо $\lambda=b=0$ и мы получаем двумерную алгебру, либо $w=-\frac c2+w_1$. Итак, алгебра размерности больше двух получается только для $T = G\left(g(bc - \lambda)^3\right) - c/2$, $\mathcal{F}=\frac 1tG\left(t(bc - \lambda)^3f\right) - \frac {c}{2t}$.

Если $b\ne 0$, то сдвигом $c$ (то есть заменой $\bar t=t$, $\bar x=x-\frac \lambda bt$) эта функция приводится к $\mathcal{F}=\frac 1tG\left(tc^3f\right) - \frac {c}{2t}$ из семейства III.3.2 (при $m=0$).

Если же $b=0$ (но $\lambda\ne 0$), то функция $\mathcal{F}=\frac 1tG(tf)- \frac {c}{2t}$ попадает в семейство III.4.3.

\subsubsection{Анзац II. $T = u(c)\ln |g| + w(c)$}
Здесь, как и в предыдущем параграфе, нам следует рассматривать только те случаи из леммы 4, в которых допускается $j=0$, то есть 2, 3, 3.1, 3.2, 4, 5. Соответствующее уравнение (58) для $w(c)$ имеет вид
$$[\frac{\tau}{t} + (3c\tau_x - 2\xi_x + \tau_t)]u + (\xi_t + c(\xi_x - \tau_t) - c^2\tau_x)w_c - (\frac{\tau}t - (3c\tau_x -\xi_x + 2\tau_t))w =$$
$$= t(\xi_{tt} + c(2\xi_{tx} - \tau_{tt}) + c^2(\xi_{xx} - 2\tau_{tx}) - c^3\tau_{xx}).\eqno(58)$$

{\bf Случай 2}. Подставляя $\tau=2atx+bt$, $\xi=a(x^2-nt^2)+bx+l$, $u=u_1(c^2+n)^{3/2}$ в уравнение (58) для $w$, получаем
$$2at[3cu_1(c^2+n)^{3/2} -(n + c^2)w_c + 3c w +(n + c^2)]=0.$$

Если $a=0$, то группа двумерна, если же $a\ne 0$, то $w=w_1(c^2+n)^{3/2}+\frac{c(c^2+n)}n +\frac 32 (c^2+n)^{3/2}\ln|c^2+n|$ при $n\ne 0$ и
$w=w_1c^3-\frac c2 +3c^3\ln |c|$ при $n=0$. Это дает нам функции
$\mathcal{F}=\frac{(c^2+n)^{3/2}}t[\ln|(c^2+n)^{3/2}tf|+w_1]+\frac {c(c^2+n)}{nt}$ и
$\mathcal{F}=\frac{c^3}t[\ln|c^3tf|+w_1]-\frac c{2t}$, то есть частные случаи функций III.5 и III.2.2. \\

{\bf Случай 3}. Подстановка (напомним: при $j=0$) $\tau = at+bt^{\frac{1}{2 - m}}$, $\xi = ax+l$ в уравнение (58) для $w$ приводит его к виду
$$b\left[(3-m)u_1c^m - c w_c + mw + \frac{c(m-1)}{2 - m} \right]= 0.$$
Если $b=0$, то алгебра двумерна, если же $b\ne 0$, то $w = w_1c^m + \frac c{m-2} +(3-m) u_1c^m\ln |c|$. Мы получаем $T = u_1c^m \ln (gc^{3-m}) + w_1c^m + \frac c{m-2}$ и функцию $\mathcal{F}=\frac 1t c^m \left[u_1\ln (f t c^{3-m}) + w_1\right] + \frac c{t(m-2)}$, являющуюся частным случаем III.3.2.\\

{\bf Случай 3.1}. Подстановка $\tau = at$, $\xi = \xi(x)$ в уравнение (58) для $w$ приводит его к виду
$$2[a - \xi_x]u + c(\xi_x - a)w_c - (\xi_x - a)w = tc^2\xi_{xx}.$$
Так как зависимость от $t$ есть только в правой части, то $\xi = bx+l$. Если $b=a$, то алгебра двумерна, если же $b\ne a$, тогда $-2u_1c + cw_c - w = 0$, откуда получаем $w=w_1c+2u_1c\ln |c|$,  $T = c\left[u_1\ln(|g|c^2) + w_1\right]$, $\mathcal{F}=\frac ct\left[u_1\ln|c^2tf| + w_1\right]$, то есть частный случай функции III.3.2.\\

{\bf Случай 3.2}. И здесь, подставляя $u = u_1c^3$, $\tau = b(x)t^{-1} + \xi_xt$, $\xi = \xi(x)$ в уравнение (58), мы получим слева и справа многочлены по $t$, и расщепление уравнения по этой переменной дает пять равенств:
$$b(x)[cw_c - 3w + 2c] = 0, \qquad b'(x)[3u_1c^3 - cw_c + 3w - 2c] = 0,\qquad
b''(x)= 0,$$
$$\xi''(x)[3u_1c^3 - cw_c + 3w + c] = 0, \qquad \xi'''(x) = 0.$$
Дифференцируя первое уравнение по $x$ и складывая со вторым, получаем (поскольку $u_1\ne 0$), что $b'(x)=0$, $b={\rm const}$; из последнего уравнения $\xi_{xx}={\rm const}$.

Если $b=\xi_{xx}=0$, то алгебра двумерна.
Если $b\ne 0$, то $w = w_1c^3 + c$, $\xi_{xx}=0$, $\xi = ax + l$. Мы получаем $T = u_1c^3\ln(|g|) + w_1c^3 + c$, $\mathcal{F}=\frac 1t c^3[u_1\ln(|tf|) + w_1]+ \frac ct$, являющуюся частным случаем функции III.3.2.

Наконец, если $b=0$, но $\xi''(x)\ne 0$, то $w = w_1c^3 - c/2 + 3u_1c^3\ln c$. Тогда $T = u_1c^3\ln(|g|c^3) + w_1c^3 - c/2$, $\mathcal{F}=\frac 1t c^3[u_1\ln |tc^3f| + w_1] - \frac c{2t}$ -- это частный случай функции III.2.2.\\

{\bf Случай 4}. Подставляя $\tau = at + mb$, $\xi = ax + b\ln t + l$ в уравнение для $w$, получаем
$$b[mu_1e^{mc} + w_c - mw + 1]=0.$$
Если $b=0$, то алгебра трехмерна, если же $b\ne 0$, то $w = w_1e^{mc} + \frac 1{m} - u_1mce^{mc}$. Получаем $T = u_1e^{mc}\ln(ge^{-mc}) + w_1e^{mc} + \frac 1{m}$,
$\mathcal{F}=\frac 1te^{mc}[u_1\ln |te^{-mc}f| + w_1]+ \frac 1{mt}$, являющуюся частным случаем функции III.3.3.\\

{\bf Случай 5}. Здесь подстановка $\tau = \tau(t)$, $\xi = (2\tau_t - \frac \tau t)x + l(t)$ в уравнение (58) и расщепление его по $x$ приводит нас к уравнению, которое расщеплением по $x$ приводится к системе из двух уравнений:
$$\left(2\tau' - \frac \tau t\right)' w_c=t\left(2\tau' - \frac \tau t\right)'',\quad
3(\frac{\tau}{t} - \tau')u_1 + (l' + c(\tau' - \frac \tau t))w_c = t(l'' + c(3\tau'' - 2\frac{\tau'}{t} + 2\frac{\tau}{t^2})).$$

Рассмотрим два варианта: либо $w(c)\ne \lambda c+w_1$, либо $w(c)=\lambda c+w_1$.

В первом варианте мы немедленно из первого уравнения получаем, что  $\left(2\tau' - \frac \tau t\right)'=0$, $\tau = at + bt^{1/2}$, а второе уравнение приобретает вид
$$6bu_1 + 2(2t^{1/2}l_{t} - cb)w_c = 4t^{3/2}l_{tt} + cb.$$

Если $b=l_t=0$, то алгебра двумерна.

Если $b=0$, $l_t\ne 0$, то оставшееся уравнение приобретает вид $l_{t}w_c = tl_{tt}$, откуда $w = \lambda c + w_1$ и $l = pt^{\lambda + 1} + q$ (если $\lambda\ne -1$), $l = p\ln t + q$ (если $\lambda=-1$).

В первом случае $T = u_1\ln|g| + \lambda c + w_1$,
$\mathcal{F}=\frac 1t [u_1\ln|tf| + \lambda c + w_1]$, это частный случай функции III.4.2. Во втором случае --
$T = u_1\ln|g| - c + w_1$, $\mathcal{F}=\frac 1t[u_1\ln|tf| - c + w_1]$, и здесь мы имеем уже частный случай функции III.4.3.

Если же $b\ne 0$, то, продифференцировав уравнение по $t$ и по $c$, получаем, поскольку $w \ne \lambda c + w_1$, что $l = pt^{1/2} + q$. В этом случае последнее уравнение приобретает вид
$6bu_1 + 2(p - cb)w_c = -p + cb$ и можем найти $w = -\frac c2 + 3u_1\ln(cb - p) + w_1$. Тогда $T = u_1\ln|g(cb - p)^{3}| -\frac c2 + w_1$, $\mathcal{F}=\frac 1t[u_1\ln|t(cb - p)^{3}f| -\frac c2 + w_1]$. Поскольку $b\ne 0$, после замены $\bar t=t$, $\bar x=x-\frac pb t$ мы получаем частный случай функции III.3.2.

Наконец, рассмотрим вариант $w = \lambda c + w_1$. В этом случае первое уравнение приобретает вид
$$\left(2\tau' - \frac \tau t\right)' \lambda =t\left(2\tau' - \frac \tau t\right)'',$$
а второе расщепляется на два:
$$(\tau' - \frac \tau t)\lambda = t(3\tau'' - 2\frac{\tau'}{t} + 2\frac{\tau}{t^2}),\quad 3(\frac{\tau}{t} - \tau')u_1 + l' \lambda = t l'' .$$
Из полученных трех уравнений среднее явно решается: $\tau=at+bt^{\frac{\lambda+2}{3}}$ при $\lambda\ne 1$ и $\tau=at+bt\ln t$ при $\lambda=1$. Подстановка первого (при $\lambda\ne 1$) решения в первое уравнение дает равенство
$$b(2\lambda+1)(\lambda-1)(\lambda+2)=0,$$
и если $\lambda\ne 1,-2,-\frac 12$, то $b=0$, $\tau=at$, и из последнего уравнения $l(t)=pt^{\lambda+1}+q$. Мы получаем $T = u_1\ln|g| + \lambda c + w_1$, $\mathcal{F}= \frac 1t[u_1\ln|tf| + \lambda c + w_1]$ -- частный случай функции III.4.2 при $\lambda\ne -1$ или функции III.4.3 при $\lambda=-1$.

Тo же самое мы получаем и при $\lambda=1$: подстановка $\tau=at+bt\ln |t|$ в первое уравнение с $\lambda=1$ дает $b=0$, и мы снова получаем трехмерную алгебру для функции, являющейся частным случаем функции III.4.2.

Если $\lambda = -1/2$, то первое уравнение выполняется тождественно, $\tau = at + bt^{1/2}$, а последнее приобретает вид  $2tl'' + l' = 3u_1bt^{-1/2}$, откуда $l(t)=pt^{1/2}+q+3u_1bt^{1/2}\ln |t|$. Получаем $T = u_1\ln|g| - \frac c2 + w_1$,
$\mathcal{F}=\frac 1t[u_1\ln|tf| - \frac c2 + w_1]$, и мы снова получаем частный случай функции III.4.2 (впрочем, этот частный случай на самом деле сводится к функции II.5, но здесь это уже непринципиально).

Если же $\lambda=-2$, то $\tau=at+b$, первое уравнение выполняется тождественно, а последнее уравнение дает 
$l(t)=pt^{-1}+q+3u_1b\ln t$.
Получаем $T = u_1\ln|g| - 2c + w_1$, $\mathcal{F}=\frac 1t[u_1\ln|tf| - 2c + w_1]$, и эта функция является частным случаем III.4.2 (также приводящимся к II.5).

\subsubsection{Анзац III. $T = u(c)g^j + w(c)$}
Здесь нам остается, для завершения классификации, пользуясь результатами леммы 4, решить уравнение (60):
$$w_c(\xi_t + c(\xi_x - \tau_t) - c^2\tau_x) -  w(\frac{\tau}t - (3c\tau_x -\xi_x + 2\tau_t)) = t(\xi_{tt} + c(2\xi_{tx} - \tau_{tt}) + c^2(\xi_{xx} - 2\tau_{tx}) - c^3\tau_{xx}).\eqno(60)$$
При этом нам придется перебрать уже все случаи из леммы 4, поскольку $j$ здесь может принимать любые значения.

{\bf Случай 1.} Подставляя $\tau = at+b$, $\xi = ax+l$ в уравнение для $w$, получаем $wb = 0$. Если $b=0$, то алгебра двумерна. Если же $b\ne 0$, то $w=0$ и мы получаем функцию $\mathcal{F}=u(c)f$, представляющую семейство III.2.3.\\

{\bf Случай 2.}
Подстановка $\tau=2atx+bt$, $\xi=a(x^2-nt^2)+bx+l$ в уравнение для $w$ дает
$$a[w_c(n + c^2) -3wc -(n+ c^2)]=0.$$
Если $a=0$, то алгебра двумерна. Если же нет -- $w=w_1(c^2+n)^{3/2}+ \frac{c(c^2+n)}{n}$, и мы получаем $\mathcal{F}=\frac{(c^2+n)^{3/2}}t[((c^2+n)^{3/2}tf)^j+w_1]+ \frac{c(c^2+n)}{nt}$, то есть частный случай функции III.5.

{\bf Случай 2.1.} Подстановка $\tau_t = \xi_x$, $n\tau_x + \xi_t= 0$ в уравнение для $w$ дает

$$w_c(n+c^2)\tau_x +  w(\frac{\tau}t - 3c\tau_x  - \tau_t) = (n+c^2)t(\tau_{tx} + c\tau_{xx}).$$

Если $\tau_x = 0$, $\tau-t\tau_t = 0$, то $\xi_t = 0$, $\tau_{tt} = 0$, и алгебра оказывается двумерной.

Если $\tau_x = 0$, но $\tau-t\tau_t\ne 0$, то $\xi_t = 0$ и $w(\tau - t\tau_t) = 0$. В силу нашего предположения второй множитель -- не нуль, поэтому $w = 0$, $\xi_x=\tau_t={\rm const}$, и мы получаем $T = u_1(c^2 + n)^3 g$, $\mathcal{F}=u_1(c^2 + n)^3 f$ то есть функция II.2 (при $n\ne 0$) или II.1 (при $n=0$).

Если $\tau_x \ne 0$, $\tau-t\tau_t\ne 0$, то поделим уравнение на $\tau_x$ и продифференцируем по $x$ и по $t$. Получим
$$w\left(\frac{\frac{\tau}t - \tau_t}{\tau_x} \right)_t = (c^2+n)\left(\frac{ct\tau_{xx} + t\tau_{tx}}{\tau_x}\right)_t,\qquad w\left(\frac{\frac{\tau}t - \tau_t}{\tau_x} \right)_x = (c^2+n)\left(\frac{ct\tau_{xx} + t\tau_{tx}}{\tau_x}\right)_x.$$

Если и в том, и в другом равенстве коэффициенты при $w$ равны нулю, то
$\frac{\tau}t - \tau_t = \alpha \tau_x$, $t\tau_{tx} = \beta\tau_x$, $t\tau_{xx} = \gamma\tau_x$, причем $\alpha\ne 0$ в силу предположения $\frac{\tau}t - \tau_t \ne 0$. Последовательно решая три уравнения получаем, $\tau = k_1t(x-\alpha t) + k_2t$. Из совместности вторых производных $\xi$: $k_1\alpha = 0$, и, поскольку $\alpha\ne 0$, получаем $k_1 = 0$, и алгебра двумерна.

Если же один из коэффициентов при $w$ ненулевой, тогда $w$ имеет вид $w = (c^2 + n)(Ac+B)$. Подстановкой этого анзаца в уравнение и расщеплением по $c$ получим систему
$$t\tau_{xx}=A(\frac{\tau}t - \tau_t) - B\tau_x, \quad t\tau_{tx}=An\tau_x + B(\frac{\tau}t - \tau_t). $$
Добавляя сюда следующее из равенств $\tau_t=\xi_x$, $\tau_x+n\xi_t=0$ уравнение
$$t\tau_{tt}=-An(\frac{\tau}t - \tau_t) + Bn\tau_x,$$
мы получаем линейную систему второго порядка в нормальной форме, условия согласования этой системы дают равенства
$$2A(\frac{\tau}t - \tau_t)=2B(\frac{\tau}t - \tau_t) =0,$$
из которых получаем следующие варианты:

а) $A=B=0$, в этом случае $w=0$, $\tau_{tx}=\tau_{xx}=\tau_{tt}=0$, $\tau=at+bx+d$, $\xi=ax-bnt+l$, и мы снова получаем функцию II.2 или II.1.

б) $A^2+B^2\ne 0$, тогда $\tau_t=\frac \tau t$, $\tau=t\cdot h(x)$, где $h''(x)=0$. Если $h'(x)=0$, то алгебра двумерна. Если $h'(x)\ne 0$, то $B=0$, $A=\frac 1n$, и мы получаем $\mathcal{F}=(c^2+n)^{3}f+ \frac{c(c^2+n)}{nt}$, то есть частный случай функции III.5.\\

{\bf Случай 3.} Подставляя $\tau = pt^{\frac{j-1}{m - (j+2)}} + at$, $\xi = ax+b$ в уравнение для $w$, получаем
$$\left[-cw_c(j-1) -  w(m-3j) + c(j-1)\frac{2j+1 - m}{m - (j+2)}\right] p = 0.$$
При $p=0$ алгебра двумерна, при $p\ne 0$ мы, в силу
$j\ne 1$, получаем $w= w_1c^{\frac{3j-m}{j-1}} + \frac{j-1}{j+2-m}c$, $T = u_1c^{m}g^j + w_1c^{\frac{3j-m}{j-1}} + \frac{j-1}{j+2-m}c$, $\mathcal{F}=\frac 1tc^{\frac{3j-m}{j-1}}\left[u_1(c^{\frac{m-3}{j-1}}tf)^j + w_1\right] + \frac{j-1}{j+2-m}\frac ct$. Этo частный случай функции III.3.2.\\

{\bf Случай 3.1.}
Подстановка $\tau = at$, $\xi = \xi(x)$ в уравнение для $w$ дает
$$(\xi_x - a)(cw_c - w) = tc^2\xi_{xx}.$$
Тогда $\xi = px+b$. Если $p=a$, то алгебра двумерна, если же $p\ne a$, то $w = w_1c$,  $T = u_1c^{2j+1}g^{j} + w_1c$, $\mathcal{F}=\frac 1t[u_1c^{2j+1}(tf)^{j} + w_1c]$, это частный вид функции III.2.2.
\\

{\bf Случай 3.2.}
Подставляя $\xi = \xi(x)$, $\tau =\xi_xt+ p(x)t^{\frac{j-1}{2j+1}}$ в уравнение для $w$, получаем
$$-w_c\left[c\frac{j-1}{2j+1}pt^{\frac{j-1}{2j+1}-1}+ c^2p_xt^{\frac{j-1}{2j+1}} + c^2\xi_{xx}t\right] + w\left[3cp_xt^{\frac{j-1}{2j+1}} + 3c\xi_{xx}t - \frac{3}{2j+1}pt^{\frac{j-1}{2j+1}-1}\right] =$$
$$= -c\frac{j-1}{2j+1}(\frac{j-1}{2j+1}-1)pt^{\frac{j-1}{2j+1}-1} - c^22\frac{j-1}{2j+1}p_xt^{\frac{j-1}{2j+1}} -c^2 t\xi_{xx} - c^3p_{xx}t^{\frac{j-1}{2j+1}+1} - c^3\xi_{xxx}t^2.$$
И в левой, и в правой части этого равенства стоят степени переменной $t$ с показателями $\frac{j-1}{2j+1}-1$, $\frac{j-1}{2j+1}$, $\frac{j-1}{2j+1}+1$, $1$ и $2$. Поэтому расщепление уравнений по переменной $t$ зависит от того, совпадает ли величина $\frac{j-1}{2j+1}$ с одним из чисел $0$, $1$, $2$ и $3$. То есть от того, совпадает ли $j$ с одним из чисел $1$, $-2$, $-1$  и $-\frac 45$. Поскольку варианты $j=-2$ и $j=-1$ в рассматриваемом варианте исключаются, у нас остается три версии: $j=-\frac 45$, $j=1$ и $j$ не совпадает ни с одним из приведенных чисел.

а) если $j = -4/5$, то расщепление по $t$ дает четыре уравнения:
$$[- cw_c +  3w + 6c]p_x = 0,\quad
[- 3cw_c + 5w + 6c]p + c^3\xi_{xxx} = 0,\quad
[- cw_c +  3w + c]\xi_{xx} = 0,\quad
p_{xx} = 0.$$

Из последнего находим, что $p = p_1x + p_2$.

Если $\xi_{xx} = 0$ и $p=0$, то алгебра  двумерна. Если $\xi_{xx} = 0$ и $p\ne0$, то из второго уравнения $w = w_1c^{5/3} - 3c$, а тогда из первого следует, что либо $p=const$, либо $w_1=0$. В первом случае получаем $T = u_1c^{3/5}g^{-4/5} + w_1c^{5/3} - 3c$, $\mathcal{F}=\frac 1t[u_1c^{3/5}(tf)^{-4/5} + w_1c^{5/3} - 3c]=\frac {c^{5/3}}t[u_1(c^{4/3}tf)^{-4/5} + w_1] - \frac{3c}t$, то есть частный случаю функции III.3.2. Во втором случае  мы получаем функцию $\mathcal{F}$, определяемую той же формулой с $w_1=0$, которая принадлежит уже семейству  II.1 (при $j=-\frac 45$).

б) если $j = 1$, то мы получаем уже другой квартет уравнений:
$$p_x[-cw_c +  3w] = 0, \quad [-cw_c + 3w + c]\xi_{xx} + c^2p_{xx} = 0, \quad wp =0, \quad \xi_{xxx} = 0.$$

Если $p\ne 0$, то $w = 0$, $\xi_{xx} = 0$, $p_{xx} = 0$, и мы получаем $T = u_1c^{6}g$, $\mathcal{F}=u_1c^{6}f$, эта функция также принадлежит семейству II.1.

Если же $p = 0$, то если $\xi_{xx}=0$, то алгебра двумерна, а если $\xi_{xx}\ne 0$, то $T =u_1c^{6}g + w_1c^3 - c/2$, $\mathcal{F}=u_1c^{6}f + \frac {w_1c^3}t - \frac c{2t}$. Это частный случай III.2.2.

в) $j\ne -2,-1,-4/5, -1/2, 1$, в этом случае равенств уже пять:
$$p[cw_c(j-1) +  3w + c(j-1)\frac{j+2}{2j+1}] = 0,\qquad
p_x[-cw_c +  3w + 2c\frac{j-1}{2j+1}] = 0,$$
$$[- cw_c +  3w  + c]\xi_{xx} = 0,\qquad \xi_{xxx} = 0, \qquad p_{xx} = 0.$$

Если $p = 0$, то если $\xi_{xx}=0$, то алгебра двумерна. Если же $\xi_{xx}\ne 0$, то $T =u_1c^{3(j+1)}g^{j} + w_1c^3 - c/2$, $\mathcal{F}=\frac 1t[u_1c^{3(j+1)}(tf)^{j} + w_1c^3 - c/2]$ и мы снова получаем частный случай функции III.2.2.

Если $p\ne 0$, то $w = w_1c^{\frac{3}{1-j}} + c\frac{1-j}{2j+1}$ и остальные уравнения дают нам
$$p_x[w_1\frac{3j}{1-j}c^{\frac{3}{1-j}}] = 0,\qquad
[- \frac{3j}{1-j}w_1c^{\frac{3}{1-j}} + c\frac{3}{2j+1}]\xi_{xx} = 0,\qquad \xi_{xxx} = 0, \qquad p_{xx} = 0.$$
Поскольку $j\ne -2$, во втором уравнении степени $c$ разные, а значит, $\xi_{xx}=0$. В первом же уравнении, в силу $j\ne 0$, мы получаем  $p_xw_1=0$, что дает два варианта.

Если $w_1 = 0$, то $\xi = ax+b$, $p=p_1x+p_2$, $T = u_1c^{3(j+1)}g^{j} + c\frac{1-j}{2j+1}$, $\mathcal{F}=\frac 1t[u_1c^{3(j+1)}(tf)^{j} + c\frac{1-j}{2j+1}]$ -- это функция из семейства II.1. Если же $w_1\ne 0$, то $p = const$, тогда $T = u_1c^{3(j+1)}g^{j} + w_1c^{\frac{3}{1-j}} + \frac{1-j}{2j+1}c$,
$\mathcal{F}=\frac 1t[u_1c^{3(j+1)}(tf)^{j} + w_1c^{\frac{3}{1-j}} + \frac{1-j}{2j+1}c]$, и мы получаем частный случай функции III.3.2.\\

{\bf Случай 3.3.}
Подстановкой $\tau = \tau(t)$, $\xi = \xi(x)$ в уравнение для $w(c)$ получаем уравнение
$$cw_c(\xi_x - \tau_t) -  w(\frac{\tau}t + \xi_x - 2\tau_t) = t(- c\tau_{tt} + c^2\xi_{xx}).$$

Если $\xi_x=\tau_t$, то $\tau = at+b$, $\xi = ax+l$ и $wb = 0$. Если $b=0$, то алгебра двумерна, иначе $w=0$,  $T = u_1c^{3}g$, $\mathcal{F}=u_1c^{3}f$ и мы получаем функцию из семейства II.2.

Если же $\xi_x\ne\tau_t$, то, зафиксировав на время соответствующее значение $(t,x)$, можно увидеть, что $w(c)$ является решением линейного обыкновенного дифференциального уравнения $cw_c-kw=\alpha c^2+\beta c$, которое, в зависимости от значений коэффициентов $k$, $\alpha$, $\beta$ может иметь один из следующих видов: $w = pc^k + qc^2 + rc$ ($k\ne 1,2$), $w = pc^2\ln |c| + qc^2 + rc$, $w = pc\ln|c| + qc^2 + rc$. Нам удобно предполагать во всех трех вариантах $p\ne 0$, выделив в качестве четвертого $w=qc^2+rc$. Рассмотрим эти четыре варианта.

I. Вариант $w = pc^k + qc^2 + rc$, $k\ne 1,2$, $p\ne 0$. Подстановка этой функции в уравнение дает
$$(pkc^k + 2qc^2 + rc)(\xi_x - \tau_t) - (pc^k + qc^2 + rc)(\frac{\tau}t + \xi_x - 2\tau_t) = t(- c\tau_{tt} + c^2\xi_{xx}).$$
Приравнивание коэффициентов при одинаковых степенях $c$ расщепляет это уравнение на три:
$$p((k-1)\xi_x + (2-k)\tau_t - \frac{\tau}t) = 0, \quad q(\xi_x - \frac{\tau}t) = t\xi_{xx}, \quad r(\tau_t - \frac{\tau}t) = - t\tau_{tt}.$$
Поскольку $p\ne 0$, мы из первого уравнения получаем $(k-1)\xi_x=a$, $(2-k)\tau_t - \frac{\tau}t=-a$, $\xi=\frac a{k-1}x+l$, $\tau=\frac a{k-1}t+bt^{1/(2-k)}$. Подстановка этих функций в остальные уравнения приводит их к виду $qb=(r+\frac 1{2-k})b=0$. Если $b=0$, то алгебра двумерна, в противном случае $q=0$, $r=\frac 1{k-2}$ и мы получаем $T = u_1c^{3}g+pc^k+\frac c{k-2}$, $\mathcal{F}=u_1c^{3}f+\frac pt c^k+\frac c{t(k-2)}$, эта функция является частным случаем функции III.3.2.

II. Вариант $w = pc^2\ln c + qc^2 + rc$, $p\ne 0$. Подстановка в уравнение приводит его к виду
$$(2pc^2\ln c + pc^2 + 2qc^2 + rc)(\xi_x - \tau_t) - (pc^2\ln c + qc^2 + rc)(\frac{\tau}t + \xi_x - 2\tau_t) = t(- c\tau_{tt} + c^2\xi_{xx}),$$
откуда, ввиду линейной независимости $c$, $c^2$ и $c^2\ln c$ получаем три уравнения: $$p(\xi_x - \frac{\tau}t) = 0,\qquad
(p+2q)(\xi_x-\tau_t)-q(\frac{\tau}t + \xi_x - 2\tau_t) = t\xi_{xx},\qquad
r(\tau_t - \frac{\tau}t)=-t\tau_{tt}.$$
Из первого, ввиду $p\ne 0$, получаем $\xi=ax+l$, $\tau=at$, то есть алгебра может быть только двумерной.

III. Вариант $w = pc\ln c + qc^2 + rc$, $p\ne 0$. Подстановка этой функции в уравнение приводит его к виду
$$(pc\ln c +pc+ 2qc^2+rc)(\xi_x - \tau_t) - (pc\ln c + qc^2 +rc)(\frac{\tau}t + \xi_x - 2\tau_t) = t(- c\tau_{tt} + c^2\xi_{xx}),$$
и расщепление по $c$ дает три равенства:
$$p(\tau_t - \frac{\tau}t) = 0,\qquad q(\xi_x - \frac{\tau}t) = t\xi_{xx},\qquad p\xi_x - (p - r)\tau_t - r\frac{\tau}t = - t\tau_{tt}.$$
Из первого, в силу $p\ne 0$, следует $\tau=at$, тогда остальные два равенства приобретают вид
$q(\xi_x - a) = t\xi_{xx}$, $p(\xi_x-a)=0$, и опять же, в силу $p\ne 0$, получаем $\xi=ax+l$, то есть двумерную алгебру.

IV. Вариант $w =qc^2 + rc$. Подстановка этой функции в уравнение приводит его к виду
$$(2qc^2+rc)(\xi_x-\tau_t) - (qc^2 +rc)(\frac{\tau}t + \xi_x - 2\tau_t) = t(- c\tau_{tt} + c^2\xi_{xx}),$$
и расщепление по $c$ дает два равенства:
$$q(\xi_x - \frac{\tau}t) = t \xi_{xx},\qquad r(\tau_t - \frac{\tau}t)=-t\tau_{tt}.$$
Второе из них имеет решение $\tau=at+bt^{-r}$ при $r\ne -1$ и $\tau=at+bt\ln |t|$ при $r=-1$. Подставляя эти функции в первое равенство, получаем соответственно
$$q(\xi_x - a-bt^{-r-1}) = t \xi_{xx}\qquad \mbox{ или }\qquad
q(\xi_x - a-b\ln |t|) = t \xi_{xx}.$$
Из первого равенства при $r\ne -2$, так же, как и из второго, получаем $q(\xi_x-a)=qb=\xi_{xx}=0$. Если
$q\ne 0$, то $b=0$, $\xi_x=a$, алгебра двумерна. Если же $q=0$, то $\xi=kx+l$, $T=u_1c^3f+rc$, $\mathcal{F}=u_1c^3f+\frac rt c$, это снова частный случай функции из семейства III.3.2.

В случае же $r=-2$ первое равенство дает $q(\xi_x - a)=0$, $-qb =\xi_{xx}$, и снова если $q\ne 0$, то мы получаем двумерную алгебру: $\xi=ax+l$, $b=0$, а если $q=0$, то мы приходим к той же функции, что и в случае $r\ne -2$, то есть опять попадаем в частный случай III.3.2.\\

{\bf Случай 3.4.}
Подстановка $\tau = \xi_xt$, $\xi = \xi(x)$ в уравнение для $w$ и расщепление по $t$ дает  $\xi_{xxx} = 0$ и $[cw_c - 3w - c]\xi_{xx} = 0$. Если $\xi = ax + b$, то  алгебра двумерная. В противном случае $w = w_1c^3 - c/2$, $T = u_1c^{3/2}g^{-1/2} + w_1c^3 - \frac{c}{2}$, $\mathcal{F}=\frac 1t[u_1c^{3/2}(tf)^{-1/2} + w_1c^3 - \frac{c}{2}]$, это частный случай функции III.2.2.\\

{\bf Случай 3.5.}
Подставим $\tau = a(x)t$, $\xi= \xi(x)$ в уравнение относительно $w$. Расщепление этого уравнения по $t$ дает три уравнения:
$$a_{xx} = 0,\qquad
[-cw_c + 3w + 2c]a_x = c\xi_{xx},\qquad
[cw_c -  w](\xi_x - a) = 0.$$
Из первого уравнения получаем $a(x)=a_1x+a_2$, $\tau = (a_1x+a_2)t$. Третье уравнение предоставляет нам альтернативу.

Если $cw_c-w=0$, то есть $w=w_1c$, то $\xi = (w_1+1)a_1x^2 + b_1x + b_2$, $T =u_1c^{-3}g^{-2} + w_1c$, $\mathcal{F}=\frac 1t[u_1c^{-3}(tf)^{-2} + w_1c]$, эта функция является частным случаем III.3.2 (приводящимся, впрочем, к II.3).

Если $w\ne w_1c$, то $\xi = a_1x^2/2 + a_2x + b_2$ и $[-cw_c + 3w + c]a_1 = 0$. Если $a_1=0$, то алгебра двумерна, если же $a_1\ne 0$, то $T =u_1c^{-3}g^{-2} + w_1c^3 -c/2$, $\mathcal{F}=\frac 1t[u_1c^{-3}(tf)^{-2} + w_1c^3 -c/2]$,  эта функция является частным случаем функции III.2.2.

{\bf Случай 4.}
Подставим $\tau = at + p$, $\xi = ax - (j-1)\frac{p}{m}\ln t + b$, в уравнение для $w$, получим
$$[w_c(j-1) + mw + (j-1)]p = 0.$$
Если $p=0$, то алгебра двумерна. Если же $p\ne 0$, то
поскольку $j\ne 1$ (это, напомним, ограничение случая 4), то
$w = w_1e^{\frac{cm}{1-j}} + \frac{1-j}m$,
$T = u_1e^{mc}g^j + w_1 e^{\frac{cm}{1-j}} + \frac{1-j}m$, $\mathcal{F}=\frac 1t[u_1e^{mc}(tf)^j + w_1 e^{\frac{cm}{1-j}} + \frac{1-j}m]$, и эта функция оказывается частным случаем функции III.3.3.\\

{\bf Случай 5.}
Подставляя $\tau = \tau(t)$, $\xi = (\frac{j-1}{2j+1}\frac \tau t + \frac{j+2}{2j+1}\tau_t)x + l(t)$ в уравнение для $w$, нам будет удобно использовать обозначение $A = -\frac{j+2}{2j+1}$. Тогда алгебра имеет вид $\tau = \tau(t)$, $\xi = ((A+1)\frac \tau t - A\tau_t)x + l(t)$, и подстановка ее в уравнение и расщепление по $x$ дает два соотношения:
$$w_c((A+1)\frac \tau t - A\tau_t)_t = t((A+1)\frac \tau t - A\tau_t)_{tt},$$
$$w_c(l_t + c(A+1)(\frac \tau t - \tau_t)) -  w(A+2)(\frac \tau t - \tau_t) = tl_{tt} + c(2(A+1)(\tau_t - \frac \tau{t}) - (2A+1)t\tau_{tt}).$$
Рассмотрим два случая: когда $w(c)$ линейна и когда она не является линейной.

1. Если $w=\lambda c+w_1$, то уравнения, после расщепления второго уравнения по $c$, принимают вид
$$\lambda((A+1)\frac \tau t - A\tau_t)_t = t((A+1)\frac \tau t - A\tau_t)_{tt},$$
$$ (\lambda- 2(A+1)) (\frac \tau t - \tau_t) = (2A+1)t\tau_{tt},$$
$$\lambda l_t -  w_1(A+2)(\frac \tau t - \tau_t) = tl_{tt}.$$
Поскольку $2A+1=-\frac 1{2j+1}\ne 0$, второе уравнение является регулярным, и имеет решением $\tau=at+bt^{1+\frac{1-\lambda}{2A+1}}$ при $\lambda\ne 1$ и $\tau=at+bt\ln|t|$ при $\lambda=1$. Подстановка этих функций в первое уравнение дает
$$b\left (\lambda+1-\frac{1-\lambda}{2A+1}\right) \frac{A(1+\lambda)+1}{2A+1}\frac{1-\lambda}{2A+1} =0,\qquad (\lambda\ne 1),$$
$$b\frac 1t = bt(-\frac 1{t^2}) \qquad (\lambda=1), $$
откуда следует, что $b=0$ всегда, кроме случаев $\lambda=-\frac{A}{A+1}=\frac {j+2}{j-1}$ и $\lambda =-\frac{A+1}{A}=\frac{j-1}{j+2}$ при $j\ne -2$ и $\lambda=0$ при $j=-2$. Решение третьего уравнения дает нам двумерное множество, поэтому, за исключением указанных случаев, алгебра получается трехмерной. Соответствующая функция $\mathcal{F}=\frac 1t[u_1(tf)^j+\lambda c+w_1]$ оказывается частным случаем функции III.4.2 или III.4.3.

Исключительные случаи отличаются только тем, что алгебра четырехмерна, а соответствующие частные случаи функции III.4.2 сводятся к функции II.1.

2. Если $w\ne \lambda c+w_1$, то уже из первого уравнения мы получаем
$((A+1)\frac \tau t - A\tau_t)_t = 0$, откуда $\tau=at+bt^{1+\frac 1A}$ при $A\ne 0$ и $\tau=at$ при $A=0$.

Подстановка этой функции при $A=0$ во второе уравнение приводит его к виду $w_cl_t=tl_{tt}$, откуда, в силу предположения о нелинейности $w(c)$, следует $l_t=0$, и алгебра получается двумерной.

В случае же $A\ne 0$ подстановка во второе уравнение дает
$$w_c(l_t - c\frac {A+1}{A}bt^{\frac 1A}) +  w\frac {A+2}{A}bt^{\frac 1A} = tl_{tt} - c\frac {A+1}{A^2}bt^{\frac 1A}.$$

Если $b=0$, то мы снова приходим к уравнению $w_cl_t=tl_{tt}$, и предположение о нелинейности $w(c)$ приводит к двумерной алгебре.

Если же $b\ne 0$, то нам необходимо рассмотреть два случая: $A=-1$ и $A\ne -1$.

а) если $A = -1$ (то есть $j=1$), то $w_cl_t - wbt^{-1} = tl_{tt}$, предположение о нелинейности $w$ исключает, в силу $b\ne 0$, вырождение коэффициента $l_t$, а значит, мы получаем для $w$ анзац $w=pe^{\lambda c}+q$, где $p\ne 0$ и $q\ne 0$. Подстановка этой функции в наше уравнение и расщепление по $c$ дает равенства $p[\lambda l_t - bt^{-1}]=0$, $-qbt^{-1} = tl_{tt}$; из первого, поскольку $p\ne 0$, получаем $l=\frac b\lambda \ln|t|+l_1$, из второго, поскольку $b\ne 0$,
$q=\frac 1\lambda $. Таким образом, мы получили трехмерную алгебру для функции $\mathcal{F}=u_1f+\frac pte^{\lambda c}+\frac{1}{\lambda t}$. Это частный случай функции III.3.3.

б) $A\ne -1$. В этом случае, зафиксировав некоторое значение $t$ и значения соответствующих функций, мы получим для $w(c)$ дифференциальное уравнение, решение которого имеет вид $w=p(c-k)^{\frac{A+2}{A+1}}+qc+r$, исключений здесь нет, поскольку $\frac {A+2}{A+1}$ не может равняться ни единице (что очевидно), ни нулю (так как $A+2=\frac{3j}{2j+1}$, а вариант $j=0$ исключается условиями теоремы, так как $\mathcal{F}$ оказывается независящей от $f$).

Подстановка этой $w(c)$ в уравнение  и расщепление по $c$ приводит к равенствам
$$l_t+k\frac {A+1}{A}bt^{\frac 1A}=0,\qquad
q= - \frac {A+1}{A},\qquad
ql_t +r\frac {A+2}{A}bt^{\frac 1A} = tl_{tt},$$
из которых немедленно получаем $l(t)=-kbt^{\frac{A+1}{A}}+l_1$, $q=-\frac {A+1}{A}$,
$r=k\frac{A+1}A$, $\mathcal{F}=\frac 1t[u_1(tf)^j+p(c-k)^{\frac{3j}{j-1}}+\frac {j-1}{j+2}(c-k)]$ с трехмерной алгеброй.
$\mathcal{F}$ -- частный случай функции III.3.2.\\

{\bf Случай 5.1.}
Подстановка $\xi_x = 2\tau/t - \tau_t$ для $\tau=\tau(t,x)$, $\xi=\xi(t,x)$ приводит к равенству
$$w_c(\xi_t + 2c(\frac{\tau}{t} - \tau_t) - c^2\tau_x) -  3w(\frac{\tau}t- \tau_t - c\tau_x) = t\xi_{tt} + c(4(\tau_t - \frac{\tau}{t}) - 3t\tau_{tt}) + c^2(2\tau_x - 3t\tau_{tx}) - c^3t\tau_{xx}.$$
Зафиксируем $t$ и $x$ и воспользуемся леммой 3, позволяющей нам, используя замены вида $\bar t=t$, $\bar x=x-kt$, свести рассмотрение к следующим четырем случаям:\\
1. $w=w_1(c^2+n)^{3/2}+pc^3+qc^2+rc+s$ ($n\ne 0$),\\
2. $w=w_1 c^{3/2}+pc^3+qc^2+rc+s$,\\
3. $w=w_1 c^3+pc^3\ln c+qc^2+rc+s$,\\
4. $w=w_1c^4+pc^3+qc^2+rc+s$,\\
причем в первых трех можно без ограничения общности предполагать, что $w_1\ne 0$. Рассмотрим каждый из этих случаев.

1. Подстановка $w=w_1(c^2+n)^{3/2}+pc^3+qc^2+rc+s$ ($n\ne 0$, $w_1\ne 0$) в уравнение после расщепления приводит к системе
$$\frac{\tau}t- \tau_t=0,\quad   \xi_t+n\tau_x = 0,\quad
p(\frac{\tau}{t}-\tau_t)+q\tau_x=-t\tau_{xx},\quad
3p\xi_t+q(\frac{\tau}{t} - \tau_t)+2r\tau_x=(2\tau_x - 3t\tau_{tx}),$$
$$2q\xi_t-r(\frac{\tau}{t} - \tau_t)+3s\tau_x=4(\tau_t - \frac{\tau}{t}) - 3t\tau_{tt},\qquad
r\xi_t+3s(\frac{\tau}t- \tau_t)=t\xi_{tt}.$$
Из первого уравнения следует $\tau=th(x)$, и тогда остальные уравнения приобретают вид
$$\xi_t=-nth'(x),\,\,\,
qth'(x)=-t^2h''(x),\,\,\,
(r-pn)h'(x)=0,\,\,\,
(3s-2qn)h'(x)=0,\,\,\,
(1-r)nh'(x)=0,$$
из которых следует, во-первых, что $h''(x)=0$, а во-вторых -- что либо  $h'(x)=0$, и тогда $\tau=at$, $\xi=ax+l$ и алгебра двумерна, либо $q=s=0$, $r=1$, $p=\frac 1n$, и мы получаем трехмерную алгебру для функции $\mathcal{F}=\frac 1t[u_1(tf)^{-1}+w_1(c^2+n)^{3/2}+\frac{c(c^2+n)}n]$, что дает нам частный случай функции III.5.

2. Подстановка $w=w_1 c^{3/2}+pc^3+qc^2+rc+s$ (где $w_1\ne 0$) приводит, после расщепления по $c$, к $p=q=s=0$ и системе


$$\xi_t= \tau_x  =0,\qquad (r-4)(\tau_t - \frac{\tau}{t}) + 3t\tau_{tt}=0.$$
Условие совместности $\xi_x=2\tau/t-\tau_t$ и $\xi_t=0$ имеет вид
$t\tau_{tt}=2(\tau_t-\tau/t)$. Сопоставление этого соотношения с последним уравнением дает $(r+2)(\tau/t-\tau_t)=0$. Поэтому при $r\ne -2$ мы получаем
$\xi_t=\tau_x = \frac{\tau}t- \tau_t=\tau_{tt}=0$,
то есть двумерную алгебру. А при $r=-2$ -- трехмерную алгебру $\tau=at+bt^2$, $\xi=ax+l$ для функции $\mathcal{F}=\frac{u_1}{t^2f}+\frac{w_1c^{3/2}-2c}{t}$, являющейся частным случаем III.3.2.

3. Подстановка $w=w_1 c^3+pc^3\ln c+qc^2+rc+s$ (где $p\ne 0$) приводит, после расщепления по $c$, к той же системе
$\xi_t =\frac{\tau}{t} - \tau_t = \tau_x=0$
и той же двумерой алгебре.

4. Наконец, подстановка $w=w_1c^4+pc^3+qc^2+rc+s$ дает после расщепления по $c$ уравнения
$$w_1\tau_x= 0,\qquad 5w_1(\frac{\tau}{t} - \tau_t)=0,\qquad 4w_1\xi_t+3p(\frac{\tau}{t} - \tau_t)+q\tau_x=-t\tau_{xx},$$
$$3p\xi_t+q(\frac{\tau}{t} - \tau_t)+2(r-1)\tau_x=- 3t\tau_{tx},\qquad 2q\xi_t-(r-4)(\frac{\tau}{t} - \tau_t)+3s\tau_x= - 3t\tau_{tt},$$
$$r\xi_t-3s(\frac{\tau}{t} - \tau_t)=t\xi_{tt}.$$
Из первой тройки следует, что если $w_1\ne 0$, то мы снова получаем $\xi_t=\frac{\tau}{t} - \tau_t=\tau_x=0$ и двумерную алгебру.

Если $w_1=0$, то из оставшихся уравнений
$$3p(\frac{\tau}{t} - \tau_t)+q\tau_x=-t\tau_{xx},\qquad 3p\xi_t+q(\frac{\tau}{t} - \tau_t)+2(r-1)\tau_x=- 3t\tau_{tx},$$
$$2q\xi_t-(r-4)(\frac{\tau}{t} - \tau_t)+3s\tau_x= - 3t\tau_{tt},\qquad r\xi_t-3s(\frac{\tau}{t} - \tau_t)=t\xi_{tt}$$
первые три составляют полную систему дифференциальных уравнений второго порядка для $\tau$. Условия совместности этой системы имеют вид
$$(3p(2r-5)-2q^2)(\frac{\tau}{t} - \tau_t)+(9ps-q(r+5))\tau_x=0,\qquad (3p(2r-5)-2q^2)\xi_t+(6qs-(2r+1)(r+2))\tau_x=0.$$

В дальнейших рассуждениях мы будем, спускаясь по степени полинома, попутно производить некоторые упрощения этих систем.

а) Если $p\ne 0$, то заменой $\bar t=t$, $\bar x=x-\frac q{3p}t$ можно произвести сдвиг переменной $c$ в функции $\mathcal{F}$, а значит, и в функции $w(c)$, приведя ее к виду $w(c)=pc^3+\bar rc+\bar s$ (то есть без квадратичного члена). В этом случае система уравнений приобретет вид
$$9p(\frac{\tau}t- \tau_t) = - 3t\tau_{xx},\qquad 3p\xi_t +2(r-1)\tau_x=- 3t\tau_{tx},$$
$$-(r-4)(\frac{\tau}t- \tau_t)+3s\tau_x  = - 3t\tau_{tt},\qquad r\xi_t-3s(\frac{\tau}{t} - \tau_t)=t\xi_{tt},$$
а условия совместности -- вид
$$3p(2r-5)(\frac{\tau}{t} - \tau_t)+9ps\tau_x=0,\qquad 3p(2r-5)\xi_t-(2r+1)(r+2)\tau_x=0.$$
Из первого условия совместности, ввиду предположения $p\ne 0$, следует, что при $r\ne \frac 52$ выполнено $\tau=t\cdot h(x-kt)$, где $k=-\frac{3s}{2r-5}$. Подстановка этой формулы в первое уравнение системы второго порядка дает соотношение $9pkh'(x-kt)=3th''(x-kt)$, откуда следует, что при $k\ne 0$ (то есть при $s\ne 0$) мы получаем $h'=0$, что приводит к двумерной алгебре, а при $k=s=0$ -- только $h''=0$. Далее, подстановка $\tau=at+btx$ в уравнения дает два тождества и два равенства $3p\xi_t=-(2r+1)bt$, $t\xi_{tt}=r\xi_t$, из которых следует, что при $r\ne 1, -\frac 12$ мы получаем $b=0$ и, как следствие, снова двумерную алгебру, при $r=1$ -- трехмерную алгебру для функции $\mathcal{F}=\frac{P}{t^2f}+\frac{pc^3+c}{t}$, являющейся частным случаем функции III.5, а при $r=-\frac 12$ -- трехмерную алгебру $\tau=at+btx$, $\xi=ax+\frac 12 bx^2+l$ для функции $\mathcal{F}=\frac{P}{t^2f}+\frac{pc^3-\frac 12c}{t}$, являющейся частным случаем функции III.2.2.

Рассмотрим теперь оставшийся вариант $r=\frac 52$. В этом случае из условий совместности мы получаем $\tau_x=0$, из первого уравнения системы $\tau=at$, и решение остальных, вместе с $\xi_x=2\tau/t-\tau_t$, дает нам снова только двумерную алгебру.

б) Если теперь $p=0$, но $q\ne 0$, то уравнения принимают вид
$$q\tau_x=-t\tau_{xx},\qquad q(\frac{\tau}{t} - \tau_t)+2(r-1)\tau_x=- 3t\tau_{tx},$$ $$2q\xi_t-(r-4)(\frac{\tau}{t} - \tau_t)+3s\tau_x= - 3t\tau_{tt},\qquad
r\xi_t-3s(\frac{\tau}{t} - \tau_t)=t\xi_{tt},$$
а условия совместности -- вид
$$2q(\frac{\tau}{t} - \tau_t)+(r+5)\tau_x=0,\qquad -2q^2\xi_t+(6qs-(2r+1)(r+2))\tau_x=0.$$
Как и в предыдущем случае, из первого условия совместности следует $\tau=t h(x-kt)$, только теперь $k=\frac{r+5}{2q}$. Подстановка этой функции в первое уравнение приводит к соотношению $qh'(x-kt)=-th''(x-kt)$, откуда следует $h'=0$, $\tau=at$ и мы приходим к двумерной алгебре.

в) Если $p=q=0$, но $r\ne 0$, то, заменой $\bar t=t$, $\bar x=x+\frac sr t$ можно свести этот случай к $w=rc$, для которой, соответственно, уравнения имеют вид
$$\tau_{xx}=0,\qquad 2(r-1)\tau_x=- 3t\tau_{tx},\qquad -(r-4)(\frac{\tau}{t} - \tau_t)= - 3t\tau_{tt},\qquad r\xi_t=t\xi_{tt},$$
а условия совместности редуцируются к единственному соотношению $(2r+1)(r+2)\tau_x=0$.

Если $r\ne -2,-\frac 12$, то из этого условия мы получаем $\tau_x=0$, тогда первое и второе уравнения выполнены тождественно, из третьего при $r\ne 1$ получаем $\tau=at+bt^{\frac{4-r}{3}}$, четвертое и $\xi_x=2\tau/t-\tau_t$ дают
$\xi=x[a+b\frac{r+2}3t^{\frac{1-r}{3}}]+l(t)$,
$$b\frac{1-r}{3}\frac{2(2r+1)}{3}\frac{r+2}3=0,\qquad rl'(t)=tl''(t).$$
Поскольку мы предположили, что $r\ne -2,-1/2$, то из первого соотношения следует, что при $r\ne 1$ мы получаем $b=0$, а из второго -- $l(t)=kt^{r+1}+m$ при $r\ne -1$ и $l(t)=k\ln|t|+m$.  Алгебра оказывается трехмерной, а функция $\mathcal{F}=\frac {u_1}{t^2f}+\frac{rc}t$ -- частным случаем функции III.4.2 или III.4.3 соответственно.

При $r=1$ мы из условия согласования по-прежнему получаем $\tau_x=0$, решение третьего уравнения имеет уже вид $\tau=at+bt\ln|t|$, для второй компоненты мы получаем формулу $\xi=(a-b)x+bx\ln t+l(t)$, и подстановка этой функции в последнее уравнение приводит к
$b=0$, $l'(t)=tl''(t)$, $l(t)=kt^2+m$. Мы снова получаем трехмерную алгебру для функции $\mathcal{F}=\frac {u_1}{t^2f}+\frac{c}t$ -- частного случая функции III.4.2.

При $r=-2$ условия совместности выполнены автоматически, и решение системы из первых трех уравнений дает $\tau=at+bt^2x+dt^2$, а последнее уравнение и $\xi_x=2\tau/t-\tau_t$ -- функцию $\xi=ax+kt^{-1}+l$, мы получили пятимерную алгебру для функции $\mathcal{F}=\frac {u_1}{t^2f}-2\frac{c}t$, которая заменой $\bar t=1/t$, $\bar x=x$ превращается в функцию I.

При $r=-\frac 12$ условия совместности опять же выполнены автоматически, решение системы относительно $\tau$ имеет вид $\tau=at+bt^{3/2}+dtx$, соответственно $\xi=ax+\frac 12 bxt^{1/2}+\frac 12 dx^2+kt^{1/2}+l$, и мы получаем пятимерную алгебру для функции $\mathcal{F}=\frac {u_1}{t^2f}-\frac{c}{2t}$, которая уже заменой $\bar t=t^{1/2}$, $\bar x=x$ превращается в функцию I.

г) Наконец, если $p=q=r=0$, то заменой $\bar t=t$, $\bar x=x-st(\ln|t|-1)$ мы избавляемся от $s$, приходим к случаю $s=0$, то есть функции $\mathcal{F}=\frac {u_1}{t^2f}$, для которой условия приобретают вид
$$\tau_{xx}=0,\quad -2\tau_x=- 3t\tau_{tx},\quad
4(\frac{\tau}t- \tau_t)= - 3t\tau_{tt},\quad \xi_{tt}=0,$$
а условия совместности сводятся к $\tau_x=0$. Эти уравнения дают $\tau=at+bt^{4/3}$, из $\xi_x=2\tau/t-\tau_t$ получаем $\xi=x(a+\frac 23bt^{1/3})+l(t)$, и подстановка этой функции в $\xi_{tt}=0$ дает $b=0$, $\xi=ax+kt+m$. Алгебра получается трехмерной, а функция -- частным случаем функции III.4.2.\\

{\bf Случай 5.2.} Подставим

$\tau = at$, $\xi = \xi(t,x)$ в уравнение для $w$, получим
$$w_c(\xi_t + c(\xi_x - a)) -  w(\xi_x - a) = t(\xi_{tt} + 2c\xi_{tx} + c^2\xi_{xx}).$$

Если $\xi_t + c(\xi_x - a) = 0$, то получаем двумерную алгебру.

Если $\xi_x = a$, $\xi_t\ne 0$, то $\xi = ax+b(t)$, где $b'(t)\ne 0$. Тогда из $w_c b_t = tb_{tt}$ следует $w=w_1c + w_2$ и $b = b_1t^{w_1 + 1} + b_2$ (если $w_1\ne -1$) или $b=b_1\ln |t|+b_2$ если $w_1=-1$. Соответственно для $\mathcal{F} = u_1t^{-3/2}f^{-1/2} + \frac{w_1c + w_2}t$, частного случая функции III.4.2 или III.4.3 алгебра получается трехмерной.


Если $\xi_x \ne a$, то находим анзац
$w= mc^2 + w_1c + p + n(\alpha + c)\ln(\alpha + c)$.
Заменой $\bar t=t$, $\bar x=x+\alpha t^2/2$ мы избавляемся от слагаемого $\alpha$. Подставляя функцию
$w= w_1 c\ln|c|+ pc^2 + qc + r$.
в уравнение, получим систему:
$$w_1 \xi_t = 0, \quad p(\xi_x - a) = t\xi_{xx}, \quad 2p\xi_t + w_1(\xi_x - a) = 2t\xi_{tx}, \quad (w_1 + q)\xi_t - r(\xi_x - a) = t\xi_{tt}.$$
Из первого и третьего уравнений следует, что если $w_1\ne 0$, то $\xi_x=a$, что противоречит предположению. Значит, $w_1=0$.

Условия совместности системы имеют вид $p\xi_t=p(\xi_x-a)=0$, откуда при $p\ne 0$ мы снова получаем $\xi_x=a$, противоречащее предположению. Значит, и $p=0$.

В этом случае при $q\ne 0$ мы заменой $\bar t=t$, $\bar x=x+\frac rq t$ избавляемся от последнего слагаемого (сводя к случаю $w=qc$), и тогда
мы получаем из первых двух уравнений $\xi_x=\lambda={\rm const}$, а из третьего уравнения -- $\xi=\lambda x+bt^{q+1}+l$ при $q\ne -1$ или $\xi=\lambda x+b\ln|t|+l$ при $q=-1$. Алгебра получается четырехмерной, соответствующая функция $\mathcal{F}=u_1t^{-3/2}f^{-1/2}+\frac {qc}t$ оказывается частным случаем III.4.2, который, впрочем, приводится к II.4.

При $q=0$ мы получаем четырехмерную алгебру $\tau=at$, $\xi=\lambda x-r(\lambda -a)t\ln |t|+bt+l$ для функции $\mathcal{F}=u_1t^{-3/2}f^{-1/2}+\frac {r}t$, являющейся другим частным случаем III.4.2, тоже сводящемся к  II.4.

\subsection{Алгебры симметрий уравнений с функцией $\mathcal{F} = \mathcal{F}(t,f)$}
К такому виду, как мы уже говорили, приводятся функции, для которых соответствующее уравнение обладает как минимум двумерной алгеброй, которая является коммутативной, но не транзитивной и приводится к $\partial_x$, $t\partial_x$.

Напомним, что группа эквивалентности для этого семейства функций, помимо сдвигов и растяжений $t$ и $x$, действие которых очевидно, содержит две подгруппы, действие которых менее очевидно. Это одномерная подгруппа проективных преобразований, из которых мы будем пользоваться одним, это $\bar t=-\frac 1t$, $\bar x=\frac xt$, для которой $\bar c=ct-x$, $\bar{\mathcal{F}}=t^3\mathcal{F}$, $\bar f=f$ и бесконечномерная группа замен вида $\bar t=t$, $\bar x=x-\phi(t)$, для которой $\bar c=c-\dot\phi(t)$, $\bar{\mathcal{F}}=\mathcal{F}-\ddot\phi(t)$, $\bar f=f$.

Для функций $\mathcal{F} = \mathcal{F}(t,f)$ классифицирующее уравнение приобретает вид
$$\tau\mathcal{F}_t + f(3c\tau_x - 2\xi_x + \tau_t)\mathcal{F}_f +(3c\tau_x -\xi_x + 2\tau_t)\mathcal{F}= \xi_{tt} + c(2\xi_{tx} - \tau_{tt}) + c^2(\xi_{xx} - 2\tau_{tx}) - c^3\tau_{xx}.$$
В силу независимости $\mathcal{F}$ от $c$ мы можем рассматривать левую и правую части равенства как многочлены по $c$, и расщепить наше уравнение на четыре:
$$\tau\mathcal{F}_t + f(\tau_t- 2\xi_x)\mathcal{F}_f +(2\tau_t -\xi_x)\mathcal{F} = \xi_{tt}, \quad 3\tau_x[f\mathcal{F}_f + \mathcal{F}] = 2\xi_{tx} - \tau_{tt}, \quad \xi_{xx} - 2\tau_{tx}= 0, \quad \tau_{xx}=0.$$

Из последних двух уравнений $\tau = a(t)x + b(t)$, $\xi = a'(t)x^2 + l(t)x + n(t)$. Второе уравнение продифференцируем по $f$, получим $3\tau_x[f\mathcal{F}_f + \mathcal{F}]_f=0$. Рассмотрим два случая:\\

1. $\tau_x \ne 0$, тогда $\mathcal{F} = F_1(t)/f + F_2(t)$. С помощью замены $\bar{t} = t$, $\bar{x}= x - \varphi(t)$ мы можем упростить вид $\mathcal{F}$: $\mathcal{F} = F_1(t)/f$. Для такой функции из наших соотношений следует
$$\tau F_1'(t)+(\tau_t +\xi_x)F_1(t) =0,\qquad
\xi_{tt}=2\xi_{tx} - \tau_{tt}=\xi_{xx} - 2\tau_{tx}=\tau_{xx}=0.$$
Вторая серия равенств влечет $\tau=At^2+Btx+Ct+Ex+D$, $\xi=Atx+Bx^2+Ft+Gx+H$, и подстановка этих функций в первое равенство
после расщепления по $x$ приводит к двум равенствам:
$$[Bt+E]F_1'(t)+3B F_1(t)=0, \qquad [At^2+Ct+D]F_1'(t)+(3At+C+G)F_1(t)=0.$$
Из первого равенства следует, поскольку $F_1\ne 0$, что либо $F_1=(Bt+E)^{-3}$ (если $B\ne 0$), и после сдвига переменной $t$ мы приходим к функции $\mathcal{F}=\frac{P}{ft^3}$ из семейтсва I. Либо $F_1={\rm const}$ (если $B=0$, $E\ne 0$), и мы приходим к той же функции I напрямую. Случай $B=E=0$ дает нам $\tau_x=0$, что противоречит предположению.

2. $\tau_x = 0$, тогда $a(t)=0$, из первого соотношения получаем, что $\xi_{tt}$ не зависит от $x$. а значит $l = l_1t + l_2$, а из $2\xi_{tx} = \tau_{tt}$ получаем $\tau = l_1t^2 + b_1t + b_2$. У нас остается одно уравнение
$$(l_1t^2 + b_1t + b_2)\mathcal{F}_t + f(b_1-2l_2)\mathcal{F}_f + (3l_1t + 2b_1 - l_2)\mathcal{F} = n_{tt}.$$

Уравнение является вырожденным (то есть все коэффициенты при производных нулевые) только если $l_1=b_1=b_2=l_2=0$, $n_{tt}=0$, что дает двумерную группу. Если же не все они равны нулю, то для вида $\mathcal{F}(t,c)$ возможно несколько вариантов.

а) если $l_1\ne 0$, то решение уравнения имеет вид
$$\mathcal{F} = (l_1t^2+b_1t+b_2)^{-\frac 32}e^{-h(t)}G(fe^{2h(t)}) + N(t),$$
где $h(t)=\exp\left(-\frac 12\int\frac{b_1-2l_2}{l_1t^2+b_1t+b_2}\right)$, $N(t)$ -- некоторая функция, $G(\cdot)$ -- произвольная функция.
За счет сдвига переменной $t$ и подходящей замены $\bar t=t$, $\bar x=x-\phi(t)$ функцию $\mathcal{F}$ можно преобразовать к виду
$\mathcal{F} = (t^2\pm k^2)^{-\frac 32}e^{h(t)}G(fe^{2h(t)})$,
то есть к функции III.4.1;

б) если $l_1=0$, но $b_1\ne 0$, то
$\mathcal{F} = (b_1t+b_2)^{k-2}G(f(b_1t+b_2)^{1-2k}) + N(t)$,
что после сдвига переменной $t$ и подходящей замены $\bar t=t$, $\bar x=x-\phi(t)$ дает $\mathcal{F} = t^{k-2}G(ft^{2k-1})$,
то есть функцию III.4.2;

в) если $l_1=b_1=0$, но $b_2\ne 0$, то $\mathcal{F} = e^{-kt}G(fe^{2kt}) + N(t)$,
где $k=\frac{l_2}{b_2}$, что сводится к
$\mathcal{F} = e^{kt}G(fe^{2kt})$, то есть функции III.4.3;

г) наконец, если $l_1=b_1=b_2=0$, $l_2\ne 0$, то $\mathcal{F}=f^{-1/2}G(t)+N(t)$, что сводится к
$\mathcal{F}=f^{-1/2}G(t)$, то есть III.2.4.


\subsection{Алгебры симметрий уравнений с функцией $\mathcal{F}= cT(t,c^2f)$}
Напомним, что здесь мы имеем дело с уравнениями, алгебра симметрий которых содержит двумерную некоммутативную и нетранзитивную подалгебру, которая может быть приведена к виду $\partial_x$, $x\partial_x$.

Для этого семейства группа эквивалентности, как было показано в лемме 1, содержит, помимо указанных симметрий, бесконечномерную алгебру замен переменной $t$: $\bar t=\phi(t)$, $\bar x=x$, $\bar c=\frac c{\phi'(t)}$, $\bar f=\phi'(t)f$, $\bar {\mathcal{F}}=\frac 1{{\phi'}^2(t)}\mathcal{F}-c\frac{\phi''(t)}{\phi'(t)^3}$.

Подставим нашу функцию в уравнение (27), обозначив $g = c^2f$. Получим для функции $T = T(t,g)$ уравнение:
$$c\tau T_t + [2\xi_t + c^2\tau_x - c\tau_t]gT_g + [\xi_t + 2c^2\tau_x + c\tau_t]T = \xi_{tt} + c(2\xi_{tx} - \tau_{tt}) + c^2(\xi_{xx} - 2\tau_{tx}) - c^3\tau_{xx}.$$
Расщепление по $c$, дает нам четыре соотношения
$$\tau_{xx} = 0, \quad \tau_x[gT_g + 2T] = \xi_{xx} - 2\tau_{tx}, \quad \tau T_t - \tau_t[gT_g - T] = 2\xi_{tx} - \tau_{tt}, \quad \xi_t[2gT_g + T] = \xi_{tt}.$$

Из первого получаем $\tau = a(t)x + b(t)$, из второго и четвертого получаем, что если $T(t,g)\ne T_1(t)g^{-1/2}+T_2(t)$ и $T(t,g)\ne T_1(t)g^{-2}+T_2(t)$, то $\tau_x=\xi_t=0$, $\xi_{xx}=0$ и у нас остается единственное уравнение
$$\tau(t) T_t - \tau'(t)[gT_g - T] = - \tau''(t),$$
из которого следует, что либо уравнение вырождено, и тогда $\tau=0$ и мы имеем только двумерную группу симметрий, либо $T=\frac{G(g\tau(t))-\tau'(t)}{\tau(t)}$. Это дает нам анзац для функции $\mathcal{F}=\frac{c}{h(t)}[G(c^2fh(t))-h'(t)]$. Замена $\bar t=\varphi(t)$, $\bar x=x$ с $\varphi'(t)=1/h(t)$ приводит нашу функцию к виду
$\bar{\mathcal{F}}=\bar c G(\bar c^2 \bar f )$, принадлежащая семейству III.3.2.

Если $T(t,g)=T_1(t)g^{-1/2}+T_2(t)$, то $\mathcal{F}=T_1(t)f^{-1/2}+cT_2(t)$, и заменой $\bar t=\phi(t)$, $\bar x=x$ с $\phi''(t)=T_2(t)\phi'(t)$ мы приводим ее к виду III.2.4.

Если же $T(t,g)=T_1(t)g^{-2}+T_2(t)$, то $\mathcal{F}=T_1(t)c^{-3}f^{-2}+cT_2(t)$, той же заменой функция приводится к виду $\mathcal{F}=G(t)c^{-3}f^{-2}$, что после подстановки в уравнения дает
$$\tau_{xx} = \xi_{xx} - 2\tau_{tx}=\tau_{tt}=\xi_t= 0, \qquad \tau G_t + 3\tau_t G =0.$$
Из первой серии получаем, что $\tau=Atx+Bt+Cx+D$, $\xi=Ax^2+Ex+F$, а последнее уравнение означает, что $(\tau G^{1/3})_t=0$, и значит $\tau=H(x)G^{-1/3}(t)$. Разделение переменных в $\tau$ возможно в виде $\tau=(Ax+B)(t+k)$, так что либо $G(t)=\frac{P}{(t+k)^3}$, либо $G(t)={\rm const}$. Во втором случае мы получаем функцию $\mathcal{F}=\frac P{c^3f^2}$, в первом -- $\mathcal{F}=\frac P{(t+k)^3c^3f^2}$, которая заменой $\bar t=\ln(t+k)$, $\bar x=x$ приводится к виду $\mathcal{F}=\frac P{\bar c^3\bar f^2}+\bar c$. Обе эти функции принадлежат семейству II.3.

Таким образом, и для этого семейства алгебры симметрий размерности больше двух получаются только для функций семейств I-III, перечисленных в условиях теоремы.

На этом доказательство теоремы 4 полностью завершено.
\\

{\bf Литература.}\\

1. Платонова К. С. Групповой анализ одномерного уравнения Больцмана. I. Группы симметрий // Дифференциальные уравнения, 2017, №4, С. 538-546

2. Платонова К. С. Групповой анализ одномерного уравнения Больцмана. II. Группы эквивалентности и группы симметрий в специальном случае // Дифференциальные уравнения, 2017, №6, С. 801-813

3. Платонова К. С., Боровских А.В. Групповой анализ одномерного уравнения Больцмана. Условия сохранения физического смысла моментных величин // Теоретическая и математическая физика. 2018. Т. 195, № 3. С. 452-483.

4. Максвелл Д. К. Труды по кинетической теории. Пер.  с англ. под ред. В.В. Веденяпина и Ю.Н. Орлова – 2011.

5. Mueller I., Ruggeri T. Extended thermodynamics // Springer Tracts in Natural Philosophy, 1998, v. 37, P. 1 - 222

6.
Чепмен С., Каулинг Т., Математическая теория неоднородных газов / пер.  с англ. под ред. Н.Н.Боголюбова. М.: ИИЛ, 1960. с. 510

7. Веденяпин В. В. Кинетические уравнения Больцмана и Власова. – М. : Физматлит, 2001.

8. Черчиньяни К., Теория и приложения уравнения Больцмана / пер.  с англ. под ред. Р.Г. Баранцева. М.: "Мир", 1978, с.20

9. Bobylev A., Dorodnitsyn V. Symmetries of evolution equations with non-local operators and applications to the Boltzmann equation //Discrete and Continuous Dynamical Systems. – 2009. – Т. 24. – №. 1. – С. 35-57.

10. Grigoriev Yu.~N., Ibragimov N.~H., Kovalev V.~F., Meleshko S.~V. Symmetries of integro-differential equations: with applications in mechanics and plasma physics. – Springer Science \& Business Media, 2010. – Т. 806.

11. Овсянников Л. В. Групповой анализ дифференциальных уравнений. – " Наукa"\,, Главная редакция физико-математической литературы, 1978.

12. Олвер П. Приложения групп Ли к дифференциальным уравнениям. – М.: Мир, 1989. – Т. 639.

13. Ibragimov N.H. CRC Handbook of Lie Group Analysis of Differential Equations, vol. 1, CRC Press, 1994.

14. Ли С. Симметрии дифференциальных уравнений. в 3-х тт. Т. 2. Лекции о непрерывных группах с геометрическими и другими приложениями /пер. с нем. М.-Ижевск: НИЦ "Регулярная и хаотическая динамика\", 2011. -- 840 с.

15. Gonzalez-Lopez, A., Kamran, N., Olver, P. J., Lie algebras of vector fields in the real plane. Proceedings of the London Mathematical Society, 1992, s3-64(2), P. 339-368.

\end{document}